\documentclass[10pt,article]{article}
\usepackage{amsmath,amsfonts,amsthm,amssymb,graphicx,color}
\usepackage[english]{babel}

\makeatletter

\newcommand{\Rmnum}[1]{\expandafter\@slowromancap\romannumeral #1@}

\newtheorem{notation}{Notation}[section]
\newtheorem{notations}{Notations}[section]
\newtheorem{proposition}{Proposition}[section]
\newtheorem{property}{Property}[section]
\newtheorem{theorem}{Theorem}[section]
\newtheorem{corollary}{Corollary}[section]
\newtheorem{lemma}{Lemma}[section]
\newtheorem{definition}{Definition}[section]
\newtheorem{remark}{Remark}[section]
\newtheorem{remarks}{Remarks}[section]
\newtheorem{fact}{Fact}[section]

\def\a{\alpha}
\def\b{\beta}

\def\ga{\gamma}
\def\In{\infty}
\def\Im{\text{Im}}
\def\l{\lambda}

\def\C{\mathbb{C}}
\def\E{\mathbb{E}}
\def\P{\mathbb{P}}
\def\R{\mathbb{R}}
\def\Re{\text{Re}}
\def\N{\mathbb{N}}
\def\d{\text{d}}
\def\de{\delta}
\def\e{\text{\textit{e}}}

\def\L{\mathfrak{L}}
\def\o{\text{o}}

\def\s{\sqrt}
\def\daP{\mathcal{P}}
\def\t{\tau}
\def\th{\theta}
\def\daN{\mathcal{N}}
\def\daI{\mathcal{I}}

\def\F{\mathfrak{F}}
\def\var{\varepsilon}
\def\u{\widetilde{u}}

\begin{document}
\title{ A local limit theorem for the minimum of a random walk with markovian increasements}
\date{}
\maketitle
\centerline{Yinna YE
 $^($\footnote{Address: Department of Mathematical Sciences, Xi'an Jiaotong-Liverpool University, SIP, Suzhou, JiangSu, 215123, P. R. China.
 Email: yinna.ye@xjtlu.edu.cn}$^)$}
\vspace{2mm}

\bigskip{

 \small {\bf Abstract. }Let $(\Omega,\mathcal{F}, \P)$ be a probability space and $E$ be a finite set. Assume that $X=(X_n)$ is an irreducible and aperiodic Markov chain, defined on $(\Omega,\mathcal{F}, \P)$, with values in $E$ and with transition probability $\displaystyle P=\Big(p_{i,j}\Big)_{i,j}$. Let $(F(i,j,\d x))_{i,j\in E}$ be a family of probability measures on $\mathbb{R}$. Consider
a semi-markovian chain $(Y_n,X_n)$ on $\R\times E$ with
transition probability $\widetilde{P}$, defined by $\displaystyle \widetilde{P}\Big((u,i),A\times \{j\}\Big)=\P(Y_{n+1}\in A,X_{n+1}=j|\,Y_n= u,\,X_n=i)=p_{i,j}\,F(i,j,A)$, for any $(u,i)\in\R\times E$, any Borel set $A\subset\R$ and any $j\in E$.
We study the asymptotic behavior of the sequence of Laplace transforms
of $(X_n,m_n)$, where $m_n=\min(S_0,S_1,\cdots,S_n)$ and
$S_n=Y_0+\cdots+Y_{n-1}$. Under quite general assumptions on
$F(i,j,\d x)$, we prove that for all $(i,j)\in E\times E$,
$\sqrt{n}\;\E_i[\exp(\l m_n), X_n=j]$ converges to a
positive function $H_{i,j}(\l)$ and we obtain further information
on this limit function as $\l\rightarrow0^+$.

}

\vspace{2mm}

\small{\textit{This is the second version of "A local limit theorem for the minimum of a random walk with markovian increasements" (Apr. 2011, arXiv:1104.1554v1). In this version, author's present address is updated, typos are corrected and some notations are unified.}}

\section{Introduction and main results} \label{SecIntroLLT}

Let $(\Omega,\mathcal{F}, \P)$ be a probability space and $E$ be a finite set with $N$ elements. Assume that $X=(X_n)_{n\geq0}$ is an irreducible and aperiodic Markov chain, defined on $(\Omega,\mathcal{F}, \P)$, with values in $E$ and with transition probability $P=\Big(p_{i,j}\Big)_{i,j\in E}$. The chain $X$ admits a unique invariant
probability denoted by $\nu$. Let $\displaystyle (F(i,j,\d t))_{i,j\in E}$ be a family of probability measures on $\R$. Consider a sequence of $\R$-valued random variables $(Y_n)_{n\geq0}$ defined on $(\Omega,\mathcal{F}, \P)$, such that $(Y_n,X_n)_{n\geq0}$ is a Markov chain on $\R\times E$ with transition probability $\widetilde{P}$, defined by:\\
for any $(x,i)\in\R\times E$, any Borel set $A\subset\R$ and $j\in E$,
$$\widetilde{P}\Big((u,i),A\times\{j\}\Big)=\P(Y_{n+1}\in A,X_{n+1}=j|\,Y_n= u,\,X_n=i)=p_{i,j}\,F(i,j,A).$$
Such a chain $(Y_n,X_n)_{n\geq0}$ is called a \emph{semi-markovian chain}: once the family $\displaystyle (F(i,j,\cdot))_{i,j\in E}$ is fixed, the transitions of this chain is controlled by $(X_n)_{n\geq 0}$.
We thus consider the canonical probability space $\displaystyle \Big((\R\times E)^{\N}, \Bigl(\mathcal{B}(\R)\otimes\mathcal{P}(E)\Big)^{\otimes\N}, (\P_{(u,i)})_{(u,i)\in\R\times E}\Big)$ associated with $(Y_n,X_n)_{n\geq0}$ and, for any $(u,i)\in\R\times E$, we denoted by $\E_{(u,i)}$ the expectation with respect to $\P_{(u,i)}$. To simplify our notations, we will denote $\P_{(0,i)}$ by $\P_i$ and $E_{(0,i)}$ by $\E_i$.\\

Set
$S_0=0,  S_n=S_0+Y_1+\cdots+Y_{n}$ and $m_n=\min(S_0,S_1,\cdots,S_n).$ In the case when $E$ reduces to one point, the random variable $S_n$ is
the sum of $n$ independent and identically distributed random variables on $\R$. In this case, if
$(S_n)_{n\geq0}$ is supposed to be centered, aperiodic with a
finite variance, then for all continuous
functions with compact support on $\R_-$, one gets
$$\lim_{n\rightarrow+\In}\sqrt{n}\,\E(\varphi(m_n))=C>0,$$
with $C$  a constant depending on $\varphi$ (see \cite{koz} for instance).\\
The first goal of this paper is to extend the so-called local limit
theorem for the process $(m_n,X_n)_{n\geq0}$ associated with the semi-markovian chain $(S_n,X_n)_{n\geq0}$ defined
above. We assume once and for all the following \textbf{hypotheses H}:
\textit{\begin{description}
\item[H1]\label{hypH1} there exists $\a>0$, such that for all $\l\in\C$ with $|\Re\ \l|\leq\a$, we have
$$\sup_{(i,j)\in E\times\E}|\widehat{F}(i,j,\l)|<+\In,\quad\text{where }
\widehat{F}(i,j,\l)=\int_\R\e^{\l t}F(i,j,\d t);$$
\item[H2]\label{hypH2} there exist $n_0\geq1$ and $(i_0,j_0)\in E\times E$, such that the measure $\P_{i_0}(X_{n_0}=j_0,S_{n_0}\in\d x)$ has an absolutely continuous component with respect to the Lebesgue measure $\d x$ on $\R$;
\item[H3]\label{hypH3} $\E_\nu(S_n)=\sum_{(i,j)\in E\times E}\nu_ip_{i,j}\int_\R tF(i,j,\d t)=0.$
\end{description}}
In the case when $(S_n)_{n\geq0}$ is a random walk on $\R$ with i.i.d increasements $(Y_i)_{i\geq1}$, the hypothesis H2 becomes the `Cramer's condition', i.e. $\displaystyle \limsup_{t\rightarrow+\In}\mid\widehat{\mu}(t)\mid<1$,
where $\widehat{\mu}$ is the characteristic function of the common probability law $\mu$ of $(Y_i)_{i\geq1}$.\\
We have
\begin{theorem}\label{theoloc-min-laplace}
Under the hypotheses {\rm H}, there exists a constant $\sigma^2>0$, such
that for all $(i,j)\in E\times E$,
\begin{equation}\label{eq00}
\sqrt{n}\;\E_i(\e^{\l m_n},X_n=j)\stackrel{n\rightarrow+\In}{\longrightarrow}\frac{H_{i,j}(\l)}{\sqrt{\pi}},
\end{equation}
where $H_{i,j}(\l)>0$ for all $\l>0$ and
\begin{equation}\label{eq000}
 \lim_{\l\rightarrow0^+}\l H_{i,j}(\l)=\sqrt{\frac{2}{\pi \sigma^2}}\;\nu_j.
\end{equation}
\end{theorem}
It will be also convenient to state this result under the following form:

\begin{theorem}\label{theoloc-min-repartition}
For all $(i,j)\in E\times E$, one gets
\begin{equation}\label{eq0*}
\lim_{n\rightarrow+\In}\sqrt{n}\;\P_i(m_n\geq-x,X_n=j)=h_{i,j}(x),
\end{equation}
where the functions $(x,i)\mapsto h_{i,j}(x)$ are harmonic for $(S_n,X_n)_{n\geq0}$ and satisfy
\begin{description}
 \item [$\bullet$] for any $i,j\in E$, $x\mapsto h_{i,j}$ is increasing;
  \item[$\bullet$] $h_{i,j}(x)>0$ for $x\geq0$.
\end{description}
Furthermore,
$$h_{i,j}(x)\sim x\sqrt{\frac{2}{\sigma^2}}\;\nu_j,\quad\text{as } x\rightarrow+\In.$$
\end{theorem}
As a corollary, we obtain the following recurrence property for the
process $(m_n)_{n\geq0}$:
$$\forall x>0\text{ },\forall i\in E\text{ },\text{ }\sum_{n\geq0}\P_i(m_n\geq-x)=+\In.$$
With similar arguments, we can also precise the asymptotic behavior, as $n\rightarrow+\In$, of
the sequence
$$
\Bigl(\E_i(\e^{\l m_n-\mu S_n},X_n=j)\Bigr)_{n\geq 0}$$
for any $\l>\mu>0$ ;
in the case when the $(Y_n)$ are i.i.d (that is the case when $E$ is reduced to one point), we know that
$\displaystyle \lim_{n\rightarrow+\In}n^{3/2}\E(\e^{\l m_n-\varepsilon S_n},X_n=j)$ does exist and is $>0$. In the markovian situation we study here, a
similar result should hold with the same exponent $3/2$ which appears after a derivation; unfortunately, as far as we understand, we are not able to decide whether or not this limit does not vanish. Nevertheless, the tools used to prove Theorem \ref{theoloc-min-laplace}
and Theorem \ref{theoloc-min-repartition} allow us to state the following
``transitional result'':

\begin{theorem}\label{thm3}
For $0<\varepsilon<\l$ small enough and for all $(i,j)\in E\times E$,
$$\sum_{n=0}^{+\In}\E_i[\e^{\l m_n-\varepsilon S_n},X_n=j]<+\In.$$
\end{theorem}

The local limit theorems \ref{theoloc-min-laplace}, \ref{theoloc-min-repartition} and Theorem
\ref{thm3} have several simple consequences, which are of interest. These are natural generalizations of classical local limit theorems for $(m_n)_{n\geq0}$,
 in the case when $(S_n)_{n\geq0}$ is a random walk on $\R$ with i.i.d increments (\cite{koz}, \cite{emile}).
 A typical such application is to study the asymptotic behavior of the survival probability of a critical branching process in an i.i.d random
 environment (\cite{geig1}, \cite{quiv}). Analogous results, under appropriate conditions, hold therefore for a branching process in a markovian environment (\cite{ye}).

\section{On the spectrum of the semi-markovian chain} \label{SecSpec}

For any $\l\in\C$, consider an $\C$-valued $N\times N$  matrix $P(\l)$ defined by
$$P(\l)=\Big(P(\l)_{i,j}\Big)_{i,j\in E}, \text{ with }P(\l)_{i,j}=p_{i,j}\widehat{F}(i,j,\l)=p_{i,j}\int_{\R}\e^{\l t}F(i,j,\d t).$$
 It is easy to verify that for any $n\geq1$, $|\Re\ \l|<\a$,
$$P^{(n)}(\l)=\Big(P^{(n)}(\l)_{i,j}\Big)_{i,j}=\Big(\E_i[\e^{\l S_n},X_n=j]\Big)_{i,j}.$$

In particular, $P(0)$ is equal to the transition matrix $P$ of the Markov
chain $(X_n)_{n\geq 0}$ (and $\displaystyle P^{(n)}(0)=P^{(n)}=\Big(p_{i,j}^{(n)}\Big)_{i, j\in E}$). Its spectral radius $^(\footnote{to define the spectral radius, we first need to choose a norm on the space of  $N\times N$ matrices $A=\Big(A_{i,j}\Big)_{1\leq i,j\leq N}$ with complex coefficients ; we will set
$\|A\|:=\sup_{1\leq i,j\leq N}|A_{i,j}|.$
}^)$ is equal to $1$ since $P$ is stochastic; furthermore, since $P(0)$ is aperiodic, the eigenvalue $1$ is the unique simple eigenvalue with modulus $1$ and its associated eigenvector is $\e=\left( \begin{array}{c}
1\\
\vdots \\
1
\end{array} \right)$. According to Perron-Frobenius theorem, there thus exists a unique vector $\nu=\left( \begin{array}{c}
\nu_1\\
\vdots \\
\nu_N
\end{array} \right)$ with positive coefficients such that $\sum_{i=1}^{N}\nu_i=1$ and ${}^t\!\nu P(0)={}^t\!\nu$
 (the vector ${}^t\!\nu$ may be identified as a probability measure on $E$). So we have
 $$P=\Pi+R,$$
 where \begin{description}
         \item[$\bullet$] $\Pi$ is a matrix of rank $1$ given by
          $$\Pi=\Big(\Pi_{i,j}\Big)_{i,j\in E}=\begin{pmatrix} \nu_1&\nu_2&\cdots&\nu_N\\
                         \vdots&\vdots&\quad&\vdots\\
                         \nu_1&\nu_2&\cdots&\nu_N\end{pmatrix},$$
         \item[$\bullet$] $R$ is a matrix with spectral radius $<1$,
         \item[$\bullet$] $\Pi$ and $R$ satisfy the relation $\Pi R=R\Pi=0$.
       \end{description}

 According to the analytical perturbation theory, for $|\l|$ small enough,  $P(\l)$ has a unique eigenvalue $k(\l)$ of modulus equal to the spectral radius of
$P(\l)$ and this eigenvalue is simple. Therefore, there exists a unique vector $\nu(\l)=\left( \begin{array}{c}
\nu_1(\l)\\
\vdots \\
\nu_N(\l)
\end{array} \right)$ such that
$$\sum_{i=1}^N{\nu_i(\l)}=1$$
and ${}^t\!\nu(\l)P(\l)=k(\l){}^t\!\nu(\l)$; we can thus also define a unique vector $\e(\l)=\left( \begin{array}{c}
\e_1(\l)\\
\vdots \\
\e_N(\l)
\end{array} \right)$ such that $P(\l)\e(\l)=k(\l)\e(\l)$ and ${}^t\!\nu(\l)\e(\l)=1$. More precisely, we have the following theorem:

\begin{theorem}\label{thm55}
 Under hypotheses {\rm H1} and {\rm H2}, there exist $\ga_0<{1\over 3}$ and $0<\a_0\leq\a$
 such that
 \begin{enumerate}
  \item If $\l\in\Delta_{\a_0}:=\{\l\in\C; |\Re\ \l|,|\Im\ \l|\leq\a_0\}$, then
  \begin{equation}\label{eq2.3.1}
  P(\l)=k(\l)\Pi({\l})+R(\l),
  \end{equation}
  where
  \begin{description}
    \item[$\bullet$] $k(\l)\in\C$ is the dominant eigenvalue of $P(\l)$, and satisfies
     $$|1-k(\l)|\leq \ga_0 ;$$
    \item[$\bullet$] $\Pi(\l)$ is a rank 1 matrix, which corresponds to the projector on the 1-dimensional eigenspace associated with $k(\l)$ and is given by
    $$\Pi(\l)=\Big(\e_i(\l)\nu_j(\l)\Big)_{i,j\in E};$$

    \item[$\bullet$] $R(\l)$ is a matrix with
  spectral radius $r(R(\l))< 1-2\ga_0$.

    \item[$\bullet$] The matrices $\Pi(\l)$and $R(\l)$ satisfy the following relation:
     \begin{equation}\label{eq2.3.1.1}
      \Pi(\l)R(\l)=R(\l)\Pi(\l)=0.
     \end{equation}
  \end{description}
  Furthermore, the maps $\l\longmapsto k(\l)$,
  $\l\longmapsto\Pi(\l)$ and $\l\longmapsto
  R(\l)$ are analytic on the set $\Delta_{\a_0}$.

 \item There exists $\a'_0\leq\a_0$  and $\chi \in ]0, 1[$ such that if $|\Re\ \l|\ \leq\a'_0$ and $|\Im\ \l|\ \geq\a_0$, the
  spectral radius of $P(\l)$ satisfies the inequality
     \begin{equation}\label{spectralradiusP(lambda)}
     r(P(\l))\leq \chi <1.
     \end{equation}
 \end{enumerate}
\end{theorem}

The proof of this theorem will be stated in Appendix \ref{appendixB}.

\begin{remark}\label{remark2.2.1} {\bf From now on and for all we will assume $\alpha-0=\alpha'_0$}; by (\ref{eq2.3.1}), for $\l\in \mathbb C$ s.t. $\vert
\Re\, \l\vert \leq \alpha_0$, one gets

$\bullet$  if $\vert \Im\,\l\vert \leq \alpha_0$ (i.e. $\l \in \Delta_{\alpha_0}$) then
 \begin{equation}\label{eq2.3.2}
  (I-zP(\l))^{-1}=\frac{zk(\l)}{1-zk(\l)}\Pi(\l)+\sum_{n=0}^{+\In}z^nR^n(\l).
\end{equation}

$\bullet$ if $\vert \Im\,\l\vert \geq \alpha_0$ then
  \begin{equation}\label{eq2.3.2'}
  r(P(\l)) \leq \chi
\end{equation}
for some $\chi\in ]0, 1[$.
 \end{remark}
In this expression, one can see that, for any fixed $\l\in\Delta_{\a_0}$, the function $z\mapsto(I-zP(\l))^{-1}$ is analytic on the set of all complex numbers $\C$, excepted the points $z$ satisfying the equation $zk(\l)=1$. In the following subsection, we will give an explicit expression of the solutions of this equation, in order to give some more information of the singular points of the holomorphic function $z\mapsto(I-zP(\l))^{-1}$.\\
The hypotheses H particularly allow us to control the local expansion at $0$ of the eigenvalue $k(\l)$.

\subsection{Local expansion of the spectral radius $k(\l)$ of $P(\l)$}
In this section, for any $F:E\times E\rightarrow\mathcal{P}(\R)$ and $\l\in\C$, we set
$$P(\l,F):=\Big(P(\l,F)_{i,j}\Big)_{i,j}, \text{ with } P(\l,F)_{i,j}:=p_{i,j}\int_\R\e^{\l t}F(i,j,\d t),$$

where the matrix $\displaystyle \Big(p_{i,j}\Big)_{i,j\in E}$ is the transition probability of an irreducible and aperiodic Markov chain $X=(X_n)_{n\geq0}$
as defined at the beginning of Section \ref{SecIntroLLT}.

When there is no risk of confusion about the function $F$, we can omit the sign $F$ in this formula. (We will assume that $F$ satisfies H1, i.e. for some $\a>0$ and for all $\l\in\C$ such that $|\Re\ \l|\leq\a$, $\sup_{(i,j)\in E\times E}|\widehat{F}(i,j,\l)|<+\In$, where $\widehat{F}(i,j,\l)=\int_\R\e^{\l t}F(i,j,\d t)$.)\\
According to Rellich's analytic perturbation theory of linear operators (see
 N. Dunford and J. Schwartz 1958, VII.6, \cite{Schwartz}), we have for $\l\in\Delta_{\a_0}$,
 $$P(\l,F)=k(\l,F)\Pi(\l,F)+R(\l,F),$$
 where
  \begin{description}
    \item[$\bullet$] $k(\l,F)\in\C$ is the dominant eigenvalue of $P(\l,F)$, and satisfies $|1-k(\l,F)|\leq \ga_0$  for $0<\ga_0<{1\over 3}$;   in the particular case when $\l=0$, we get $k(0,F)=1$;
    \item[$\bullet$] $\Pi(\l,F)$ is a projection ( i.e.
 $\Pi^2(\l,F)=\Pi(\l,F)$ ) on the 1-dimensional eigenspace associated with $k(\l,F)$, and in the particular case when $\l=0$,
 $$\Pi(0,F)=\Big(\Pi_{i,j}\Big)_{i,j\in E}=\begin{pmatrix} \nu_1&\nu_2&\cdots&\nu_N\\
                         \vdots&\vdots&\quad&\vdots\\
                         \nu_1&\nu_2&\cdots&\nu_N\end{pmatrix},$$
                          with $\sum_{i\in E}\nu_i=1$ and $\forall i\in E,$ $\nu_i>0$.
    \item[$\bullet$] $R(\l,F)$ is a matrix with spectral radius $<1$ and satisfies the relation
 $$\Pi(\l,F)R(\l,F)=R(\l,F)\Pi(\l,F)=0.$$
 \end{description}
In particular, the function $\l\mapsto k(\l,F)$ is analytic on $\Delta_{\a_0}$; we now compute the first term of its local expansion.\\
We introduce the \textit{mean matrix} $M(F)$ associated with $F$ which is defined by
$$M(F)=\Big(M(F)_{i,j}\Big)_{i,j},\text{ with }M(F)_{i,j}=p_{i,j}\int_\R tF(i,j,\d\l).$$
We have the
\begin{lemma}\label{lem2.2.1}
 $k'(0,F)={}^t\!\nu M(F)\e=\sum_{i,j\in E}\nu_ip_{i,j}\int_\R tF(i,j,\d t)$.\\
 In the sequel, we will denote
 $$\ga(F):={}^t\!\nu M(F)\e=\sum_{i,j\in E}\nu_ip_{i,j}\int_\R tF(i,j,\d t).$$
\end{lemma}
\begin{proof}
 Since $P(\l,F)=k(\l,F)\Pi(\l,F)+R(\l,F)$, with $\Pi(\l,F)R(\l,F)=R(\l,F)\Pi(\l,F)=0$ and $\Pi(\l,F)^2=\Pi(\l,F)$, we have $\Pi(\l,F)P(\l,F)=k(\l,F)\Pi(\l,F)$. Using the fact that $k(0,F)=1$, the derivation of the quantities in the two hand-sides of this equality at the point $\l=0$ leads to
$$\Pi'(0,F)P(0,F)+\Pi(0,F)P'(0,F)=k'(0,F)\Pi(0,F)+\Pi'(0,F).$$
Using thus the equality $P(0,F)\e=\e$, one gets
\begin{equation}\label{eq2.1.11}
 \begin{split}
 \Pi(0,F)P'(0,F)\e&=k'(0,F)\Pi(0,F)\e\\
             &=k'(0,F)\e.
\end{split}
\end{equation}
As $P'(0,F)_{i,j}=p_{i,j}\int_\R tF(i,j,\d t)$, the equality (\ref{eq2.1.11}) implies that
$$\sum_{i,j\in E}\nu_i p_{i,j}\int_\R tF(i,j,\d t)=k'(0,F).$$
\end{proof}
\begin{corollary}
 Under the hypotheses {\rm H1} and {\rm H3}, we have $k'(0)=0$.
\end{corollary}
\begin{proof}
 This is a direct consequence of Lemma \ref{lem2.2.1}, since we suppose here that
$${}^t\!\nu M(F)\e=\sum_{i,j\E}\nu_ip_{i,j}\int_\R tF(i,j,\d t)=0.$$
\end{proof}

To compute $k''(0,F)$, we need first to ``center'' the function $F$ in the following sense:
\begin{definition} Suppose that $F=(F(i,j,\cdot))_{i,j\in E}$ and $F'=(F'(i,j,\cdot))_{i,j\in E}$ are two finite families of probability measures on $\R$. One says that $F'$ is
 \textbf{a-equivalent} to $F$, if there exists a vector $u=(u_i)_{i\in E}$, such
 that for any $i,j\in E$ satisfying $p_{i,j}\neq0$, one has
 $$F'(i,j,\cdot)=\delta_{u_j-u_i}\ast F(i,j,\cdot).$$
 \end{definition}
 This notion of equivalence is relevant since we have the
 \begin{property}\label{proA1}
  \begin{enumerate}
   \item If $F$ and $F'$ are a-equivalent and satisfy hypothesis {\rm H1}, then $k(\cdot,F)=k(\cdot,F')$ on $\Delta_{\a_0}$.
   \item For any $F:E\times E\rightarrow\mathcal{P}(\R)$ satisfying H1, there exists a function $\mathfrak{F}:E\times E\rightarrow\mathcal{P}(\R)$ which is a-equivalent to $F$ and such that $ M(\mathfrak{F})\e=\ga(F)\e=\ga(\F)\e$.
  \end{enumerate}
 \end{property}

 \begin{proof}
1. By the equality $F'(i,j,\cdot)=\delta_{u_j-u_i}\ast F(i,j,\cdot)$, for any $\l\in\Delta_{\a_0}$ and any $i,j\in E$, we have
        $$P(\l,F')_{i,j}=\e^{\l(u_j-u_i)}P(\l,F)_{i,j}.$$
        Therefore,
        \begin{equation}\label{eq2.1.111}
         \begin{split}
          P^{(n)}(\l,F')_{i,j}&=\e^{\l(u_j-u_i)}P^{(n)}(\l,F)_{i,j}\\
                                       &=\e^{\l(u_j-u_i)}\Big(k^n(\l,F)\Pi(\l,F)_{i,j}+R^{(n)}(\l,F)_{i,j}\Big).
          \end{split}
        \end{equation}
   Set $\Pi(\l,F,u):=\Big(\Pi(\l,F,u)_{i,j}\Big)_{i,j}$ with $\Pi(\l,F,u)_{i,j}:=\e^{\l(u_j-u_i)}\Pi(\l,F)_{i,j}$.\\
   According to (\ref{eq2.1.111}), for any $\l\in\Delta_{\a_0}$,
    $$\frac{P^{(n)}(\l,F')}{k^n(\l,F)}\longrightarrow\Pi(\l,F,u)\neq0,\text{ as }n\rightarrow+\In.$$
    So for any $\l\in\Delta_{\a_0}$, $|k(\l,F)|$ is equal to the spectral radius $|k(\l,F')|$ of $P(\l,F')$; there thus exists $\theta=\theta(\l)$ in $[0, 2\pi[$ such that
\begin{equation}\label{eqApp33t}
k(\l,F)=\e^{i \theta}k(\l,F').
\end{equation}
Let $\e(\l,F')$ be a non-null eigenfunction of the matrix $P(\l,F')$, corresponding to the eigenvalue $k(\l,F')$:
 \begin{equation}\label{eqApp4}
  P^{(n)}(\l,F')\e(\l,F')=k^n(\l,F')\e(\l,F').
 \end{equation}
 Using (\ref{eq2.1.111}), (\ref{eqApp33t}) and (\ref{eqApp4}), one gets for any $i\in E$,
\begin{equation}
 \begin{split}
 & k^n(\l,F')\e(\l,F')_i=\\
                    &\qquad\qquad\e^{-\l u_i}\Big[k^n(\l,F')\e^{in\th}\sum_j\e^{\l u_j}\Pi(\l,F)_{i,j}\e(\l,F')_j+\sum_j\e^{\l u_j}R^{(n)}(\l,F)_{i,j}\e(\l,F')_j\Big].
 \end{split}
\end{equation}
Let $i_\l\in E$ such that $\e(\l,F')_{i_\l}\neq0$, then
$$0\neq\e(\l,F')_{i_\l}=\e^{in\theta}{a(\l)}_{i_\l}+b(\l,n)_{i_\l},$$
where
\begin{description}
 \item[$\bullet$] $a(\l):=\Big(a(\l)_i\Big)_i$ with $a(\l)_i =\e^{-\l u_i}\sum_j\e^{\l u_j}\Pi(\l,F)_{i,j}\e(\l,F')_j$;
  \item[$\bullet$] $b(\l,n):=\Big(b(\l,n)_i\Big)_i$ with $b(\l,n)_i=\e^{-\l u_i}k(\l,F')^{-n}\sum_j\e^{\l u_j}R^{(n)}(\l,F)_{i,j}\e(\l,F')_j$.
\end{description}
 Note that $\forall i\in E,\lim_{n\rightarrow+\In}b(\l,n)_i=0$, so that
$$\lim_{n\rightarrow+\In}\e^{in\theta}=\frac{\e(\l,F')_{i_\l}}{a(\l)_{i_\l}}\neq0.$$
 We can thus conclude that $\th=0$, and so $k(\l,F)=k(\l,F')$ for any $\l\in\Delta_{\a_0}$.

2. Set $v(F):=M(F)\e-\ga(F)\e=M(F)\e-({}^t\!\nu M(F)\e)\e$. Since ${}^t\!\nu v(F)$ is null, the vector $\widetilde{u}:=\sum_{n=0}^{+\In}P^n v(F)$ exists and satisfies
       \begin{equation}\label{eqApp1t}
        \u-P\u=v(F)=M(F)\e-\ga(F)\e.
       \end{equation}
       For any $i,j\in E$, let's define a function $\F:E\times E\rightarrow\mathcal{P}(\R)$ by
       $$\F(i,j,\cdot)=\delta_{\u_j-\u_i}\ast F(i,j,\cdot).$$
       Then one obtains
       \begin{equation}\label{eqApp2}
           M(\F)\e=M(F)\e+P\u-\u.
       \end{equation}
       Using (\ref{eqApp1t}) and (\ref{eqApp2}), one has $M(\F)\e=\ga(F)\e$ and $\ga(\F)={}^t\!\nu M(\F)\e={}^t\!\nu M(F)\e=\ga(F)$.

 \end{proof}

Thank to this property, we are now able to compute $k''(0)$. We first introduce the \textit{inertial matrix} $\Sigma(F)$ associated with $F$, defined by
$$\Sigma(F):=\Big(\Sigma(F)_{i,j}\Big)_{i,j},\text{ with }\Sigma(F)_{i,j}:=p_{i,j}\int_\R t^2F(i,j,\d t).$$

\begin{property}\label{pro2.2.2}
 Let $\F:E\times E\rightarrow\mathcal{P}(\R)$ such that $\F$ is a-equivalent to $F$ and
$$M(\F)\e=\ga(F)\e.$$ Then
 $$k''(0,F)=k''(0,\F)={}^t\!\nu\Sigma(\F)\e.$$
\end{property}

\begin{proof}
 We have
       \begin{equation}\label{eq1}
        \Pi(\l,\F)P(\l,\F)=k(\l,\F)\Pi(\l,\F),
       \end{equation}
     where $k(\l,\F)$ is the unique eigenvalue of $P(\l,\F)$ of maximum absolute value with
    $$k(0,\F)=1$$ and
     $\Pi(\l,\F)$ is the corresponding eigenvector.\\
     Consider the following Taylor's formula:

     $$ k(\l,\F)=1+\l k'(0,\F)+\frac{\l ^2}{2}k''(0,\F)+\o(\l ^2),$$
     $$\Pi(\l,\F)=\Pi(0,\F)+\l \Pi'(0,\F)+\frac{\l ^2}{2} {\Pi}''(0,\F)+\o(\l ^2),$$
     $$P(\l,\F)=P(0,\F)+\l M(\F)+\frac{\l ^2}{2}\Sigma(\F)+\o(\l ^2).$$
     By identification of the coefficients of order $\l^2$ (\ref{eq1}), we get
     $$\Pi(0,\F)\Sigma(\F)+2\Pi'(0,\F)M(\F)+\Pi''(0,\F)P(0,\F)=\Pi''(0,\F)+2k'(0,\F)\Pi'(0,\F)+k''(0,\F)\Pi(0,\F).$$
     Multiplying the matrices in the two sides of this equation with $\e$ and using the facts $P(0,\F)\e=\e$, $M(\F)\e=k'(0,\F)\e$ and $\Pi(0,\F)\e=\e$, one gets
      $$k''(0,\F)={}^t\!\nu\Sigma(\F)\e.$$
      And $k''(0,F)=k''(0,\F)$ is a direct consequence of the fact that $k'(\cdot,F)=k'(\cdot,\F)$ on $\Delta_{\a_0}$.
\end{proof}
\begin{corollary}\label{thm1}
     For any $F:E\times E\rightarrow\mathcal{P}(\R)$ satisfying H1, we have $k''(0,F)=0$ if and only if $F$ is a-equivalent to $\delta_{\{0\}}$.
 \end{corollary}
 \begin{proof}
  Suppose that $F:E\times E\rightarrow\mathcal{P}(\R)$ satisfies H1, from Property \ref{pro2.2.2}, there exists $\F:E\times E\rightarrow\mathcal{P}(\R)$ such that
  $$k''(0,F)=k''(0,\F)={}^t\!\nu\Sigma(\F)\e=\sum_{i,j\in E}\nu_ip_{i,j}\int t^2\F(i,j,\d t).$$
  So that $k''(0,F)=0$ if and only if $\F=\delta_{\{0\}}$.
 \end{proof}
 \begin{corollary}\label{cor2.2.3}
 Under the hypotheses {\rm H}, we have
$$\sigma^2:=k''(0)>0.$$
 \end{corollary}
 \begin{proof}
  Suppose that $k''(0)=0$. By the definition of the semi-Markovian chain $(S_n,X_n)_{n\geq0}$, we have for a fixed $i_0\in E$, and any $n\geq1$,
 \begin{equation}\label{Appeq22}
 \P_{i_0}(S_n\in\d x)=\sum_{(i_1,\cdots,i_n)\in E^n}\left[\prod_{k=0}^{n-1}\P(i_k,i_{k+1})\right]F(i_0,i_1,\d x)\ast F(i_1,i_2,\d x)\ast\cdots\ast F(i_{n-1},i_n,\d x).
 \end{equation}
 According to Corollary \ref{thm1} and the fact that the support of $\nu$ is $E$, the measures $F(i,j,\d x)$ is a Dirac measure for any $(i,j)\in E\times E$ such that $p_{i,j}>0$. So by Formula (\ref{Appeq22}), for every $i_0\in E$ and every $n\geq1$, the law $P_{i_0}(S_n\in\d x)$ is discrete. However, the hypothesis (H2 implies that $\P_{i_0}(S_{n_0}\in\d x)$ has an absolutely component with respect to the Lebesgue measure on $\R$. This leads to a contradiction. The proof is complete.
 \end{proof}

 \subsection{The equation $zk(\l)=1$ for $z\in\C$ and $|\Re\ \l|\leq\a_0$} 
We consider here the equation
\begin{equation}\label{eq2.1.1}
 zk(\l)=1,\quad \text{for }z\in\C \text{ and }|\Re\ \l|\leq\a_0.
\end{equation}

It is shown in the previous section that $k''(0)>0$ under our conditions (H). Since $\l\mapsto k(\l)$ is analytic on the open set $\Delta_{\a_0}$, one may assume that $k''(\l)>0$ for any $\l\in]-\a_0,\a_0[$. By the implicit function theorem, for $z\in\R$, the equation
(\ref{eq2.1.1}) has at most two roots in a sub interval of
$[-\a_0,\a_0]$ ( still denoted by $[-\a_0,\a_0]$ in order to
simplify the notation). Since $k'(0)=0$, one gets $\min_{-\a_0\leq\l\leq\a_0}k(\l)=k(0)=1$. Set $q=1/\inf(k(-\a_0),k(\a_0))$, then when $z\in[q,1]$, the equation (\ref{eq2.1.1}) has exactly two solutions: one is $\l_-(z)\in [-\a_0,0]$ and another is $\l_+(z)\in [0,\a_0]$; furthermore, these two solutions coincide if and only if $z=1$, and $\l_-(1)=\l_+(1)=0$.\\
For any $\de_1,\de_2>0$ such that $q+\de_1<1$, set
$$K(\de_1,\de_2):=\{z:q+\de_1<|z|<1+\de_2,\Re\  z>0,|\Im\
 z|<\de_1\}.$$

  We will   describe  in the following sections the local behavior of  some functions of the complex variable $z\in K(\de_1,\de_2)$ but with respect to the variable $t:= \sqrt{1-z}$. In order to fix a principal determination of the function $\sqrt{\ \ } $, we introduce the subset $K^*(\de_1,\de_2)\subset K(\de_1,\de_2)$ defined by
 $$K^*(\de_1,\de_2):=\{z,q+\de_1<|z|<1+\de_2,\Re\  z>0,|\Im\  z|<\de_1,z\notin[1,1+\de_2[\}. $$ Note that the map $z \mapsto \sqrt{1-z}$ is well defined on $K^*(\de_1,\de_2)$.

By the local inversion theorem, since $k'(0)=0$ and $k''(0)>0$, one may choose $\de_1\in]0,1-q[$ and $\de_2>0$ in such a way that the two functions $z\mapsto\l_+(z)$ and $z\mapsto\l_-(z)$, defined a priori on $]q+\de_1,1+\de_2[$, admit an analytic expansion to the region $K(\de_1,\de_2)\setminus\{1\}$ and these functions remain to be the solutions of (\ref{eq2.1.1}) for $z\in K(\de_1,\de_2)\setminus\{1\}$ and $\mid\Re\ \l\mid\leq\a_0$.\\

By the above, the functions $z\mapsto\l_+(z)$ and $z\mapsto\l_-(z)$ can be decomposed on $K^*(\de_1,\de_2)$ as
\begin{equation}\label{eq2.1.2}
\l_\pm(z)=\sum_{n=1}^{+\In}(\pm1)^n \a_n(1-z)^{n/2},
\end{equation}
where $\a_n\in\C$ for any $n\geq1$. On the other hand, for any $\l$ in a neighborhood of $0$, one has
\begin{equation}\label{eq2.2.1}
 k(\l)=1+\frac{k''(0)}{2!}\l^2+\frac{k^{(3)}(0)}{3!}\l^3+\cdots.
\end{equation}
By identification of the coefficients of the terms $(1-z)$ and $(1-z)^{3/2}$ in the two sides of the equality,
\begin{equation}\label{eq2.2.2}
 k(\l_+(z))=\frac{1}{z}=\sum_{n=1}^{+\In}(1-z)^n,
\end{equation}
one obtains
$$\a_1=\sqrt{\frac{2}{k''(0)}}\text{ and }\a_2=-\frac{k^{(3)}(0)}{3(k''(0))^2}.$$
We can thus conclude that for any $z\in \overline{K}(\de_1,\de_2)$, the two solutions $\l_-(z)$ and $\l_+(z)$ of the equation (\ref{eq2.1.1}) satisfy
\begin{equation}\label{eq2.4.9}
 \l_\pm(z)=\pm\sqrt{\frac{2}{k''(0)}}(1-z)^{1/2}-\frac{k^{(3)}(0)}{3(k''(0))^2}(1-z)+O((1-z)^{3/2}).
 \end{equation}

\subsection{On the spread-out property of the transition probability}
We first introduce the
\begin{notations}\label{notation2.2.1} For any integer $N\geq 1$, let $V_N$ denote the set of $N\times N$ matrices whose coefficients are complex valued Radon measures on $\R$.

 The set  $(V_N, +, \bullet)$ is an algebraic ring, when endowed with the  sum $+$ of Radon measures
and the law $\bullet$ defined by : for any $B =\Big(B_{i,j} \Big)_{i,j\in E}$ and $C =\Big(C_{i,j} \Big)_{i,j\in E}$ in $V_N$
$$B\bullet C:=\Big((B\bullet C)_{i,j}\Big)_{i,j\in E},$$
with $(B\bullet C)_{i,j}(\d x):=\sum_{k\in E}B_{i,k}\ast C_{k,j}(\d x)$, where $\ast$ denotes the convolution of measures.

For any $n \geq 1$ we will set $B^{\bullet n}=\underbrace{B\bullet\cdots\bullet B}_{\text{n times}}=\Big(B_{i,j}^{\bullet n}\Big)_{i,j}$.

For any $[a, b]\subset \R$, we denote by $V_N[a, b]$  the subset of $V_N$  of matrices whose coefficients $\sigma$  are such that
$$
\forall \lambda \in [a, b]\quad \int _\R \exp(\lambda x) \d \vert \sigma\vert (x)<+\infty.
$$
\end{notations}

 Set $M(\d x)=\Big(p_{i,j}F(i,j,\d x)\Big)_{i,j}$, for any $i,j\in E$. Since the Markov chain $X=(X_n)_{n\geq0}$ is irreducible and $(F(i,j,\d t))_{i,j\in E}$ are probability measures on $\R$ , one gets  $M_{i,j}^{\bullet k}(\R)>0$ for any $i,j\in E$ and $k$ large enough.
 The hypothesis H2 implies that $M_{i_0,j_0}^{\bullet n_0}(\d x)$ has an absolutely continuous component. By Lemma \ref{lem2.7} of Appendix \ref{appendixA}, there exists $k_1\geq1$  such that all the terms of
 $M^{\bullet k_1}(\d x)$ have absolutely continuous components. So one gets
\begin{equation}\label{eq355}
\forall k\geq k_1,\text{ }M^{\bullet k}_{i,j}(\d x)=\varphi_{k,i,j}(x)\d x+\theta_{k,i,j}(\d x),
\end{equation}
where for any $(i,j)\in E\times E$,
\begin{description}
  \item[$\bullet$] the function $\varphi_{k,i,j}$ is
positive,   belongs to $\mathbb{L}^1(\R,\d x)$ and satisfies
$0< \int\varphi_{k,i,j}(x)\d x\leq1$;
 \item[$\bullet$] $\theta_{k,i,j}(\d
x)$ is a singulary measure with respect to the Lebesgue measure such
that $0\leq \theta_{k,i,j}(\R)<1$.
\end{description}
For $|\Re\ \l|\leq{\a}_0$ and any $k\geq1$, set
\begin{gather*}\Phi_k(\d x)=\Big(\Phi_{k,i,j}(\d x)\Big)_{i,j}=\Big(\varphi_{k,i,j}(x)\d x\Big)_{i,j},\quad \Theta_k(\d x)=\Big(\Theta_{k,i,j}(\d x)\Big)_{i,j}=\Big(\theta_{k,i,j}(\d x)\Big)_{i,j};\\
 \mathfrak{L}(\Phi_k)(\l)=\int_\R\e^{\l u}\Phi_k( u)=\Big(\widehat{\varphi}_{k,i,j}(\l)\Big)_{i,j},\quad \mathfrak{L}(\Theta_k)(\l)=\int_{\R}\e^{\l u}\Theta_k(\d
  u)=\Big(\widehat{\theta}_{k,i,j}(\l)\Big)_{i,j}.
\end{gather*}
For every $(i,j)\in E\times E$, the measure $\Phi_{k,i,j}(\d x)$ is the
absolutely continuous component of $M_{i,j}^{\bullet k}(\d x)$ and $\Theta_{k,i,j}(\d
x)$ is its orthogonal component with respect to the Lebesgue
measure; the functions $\mathfrak{L}(\Phi_k)(\l)$ and
$\mathfrak{L}(\Theta_k)(\l)$ are their respective Laplace transforms (recall that the Laplace transform of $M$ is $\mathfrak{L}(M)(\l)=P(\l)$).\\
By (\ref{eq355}) and the above notations, we have for any $p\geq 1$ and $k\geq k_1$,
\begin{equation}\label{eq356}
M^{\bullet kp}(\d x)=\left(\Phi_{k}(\d x)+\Theta_{k}(\d x)\right)^{\bullet p}=\Phi_{kp}(\d x)+\Theta_{kp}(\d x),
\end{equation}
so that
\begin{equation}\label{eq357}
\Theta_{kp}(\d x)\leq \Theta^{\bullet p}_{k}(\d x).
\end{equation}
We have the following lemma:

\begin{lemma}\label{lem2.8}
 Let $k_1\geq1$ such that (\ref{eq355}) holds. There  exists $m_1\geq 1$, such that, for $q\leq z\leq1$,
\begin{equation}\label{eq35*}
 \|\mathfrak{L}(\Theta_{k_1}^{\bullet m_1})(\l_+(z))\|<z^{-k_1 m_1}.
\end{equation}
\end{lemma}

\begin{proof}
 For any $n\geq1$, one gets
 \begin{equation*}
   \left\|\mathfrak{L}(\Theta_{k_1}^{\bullet n})(\l_+(z))\right\|\leq\left\|P^{nk_1}(\l_+(z))\right\|,
 \end{equation*}
 which readily implies
 $$\rho_{\Theta_{k_1}}(\l_+(z)):=\lim_{n\rightarrow+\In}\left\|\mathfrak{L}(\Theta_{k_1}^n)(\l_+(z))\right\|^{1/n}\leq\lim_{n\rightarrow+\In}\left\|P^{nk_1}(\l_+(z))\right\|^{1/n}=k^{k_1}(\l_+(z)),$$
 where $\rho_{\Theta_{k_1}}(\l)$ denotes the spectral radius of $\mathfrak{L}(\Theta_{k_1})(\l)$ for any $\l\in\C$. The equality $zk(\l_+(z))=1$ thus leads to
 \begin{equation}\label{eq34}
 \rho_{\Theta_{k_1}}(\l_+(z))\leq z^{-k_1}.
 \end{equation}
 Let us now prove that this inequality is strict.
 Otherwise, one should have
 $$\rho_{\Theta_{k_1}}(\l_+(z))=z^{-k_1}=k^{k_1}(\l_+(z)),$$
 which should give $1=z^{k_1}k^{k_1}(\l_+(z))=z^{k_1}\rho_{\Theta_{k_1}}(\l_+(z)).$
 Since $\rho_{\Theta_{k_1}}({\l}_+(z))$ is an eigenvalue of $\mathfrak{L}(\Theta_{k_1})({\l}_+(z))$, there would exist
 a non negative vector $\a_+(z)$, such that
 $$\mathfrak{L}(\Theta_{k_1})({\l}_+(z))\a_+(z)=\rho_{\Theta_{k_1}}({\l}_+(z))\a_+(z)=z^{-k_1}\a_+(z).$$
By the definition of $\mathfrak{L}(\Theta_{k_1})$, one gets $\mathfrak{L}(\Theta_{k_1})({\l}_+(z))=P^{k_1}({\l}_+(z))-\mathfrak{L}(\Phi_{k_1})({\l}_+(z))$, so we would get
 \begin{equation}\label{eqstar}
  \begin{split}
   0&=\Pi(\l_+(z))\left(I-z^{k_1}\mathfrak{L}(\Theta_{k_1})(\l_+(z))\right)\a_+(z)\\
    &=\Pi(\l_+(z))\left[I-z^{k_1}P^{k_1}(\l_+(z))+z^{k_1}\mathfrak{L}(\Phi_{k_1})({\l}_+(z))\right]\a_+(z).\\
  \end{split}
 \end{equation}
  The equalities (\ref{eq2.3.1}), (\ref{eq2.3.1.1}) and the fact that $zk(\l_+(z))=1$ give
 $$\Pi(\l_+(z))\left[I-z^{k_1}P^{k_1}({\l}_+(z))\right]= [1-z^{k_1}k ^{k_1}({\l}_+(z))]\Pi(\l_+(z))=0.$$
Consequently, (\ref{eqstar}) leads to the equality
\begin{equation}\label{eq358}
0=z^{k_1}\Pi(\l_+(z))\left[\mathfrak{L}(\Phi_{k_1})({\l}_+(z))\right]\a_+(z),\text{ for } q\leq z\leq1.
\end{equation}
However, since all the terms of matrix
$\mathfrak{L}(\Phi_{k_1})({\l}_+(z))$ are strictly positive, the vector
$\mathfrak{L}(\Phi_{k_1})({\l}_+(z))\,\a_+(z)$
is strictly positive and the non-negative matrix
$\Pi(\l_+(z))$  has rank 1. We hence obtain
$$\Pi(\l_+(z))\left[\mathfrak{L}(\Phi_{k_1})({\l}_+(z))\right]\a_+(z)\neq0.$$
This contradicts (\ref{eq358}). So if we take $m_1$ large enough, we can thus obtain (\ref{eq35*}).
\end{proof}

\textbf{From now on, we fix $k_1$, $m_1\geq1$ such that (\ref{eq35*}) holds and we set $n_1:=k_1m_1$.} We now fix $\kappa>0$ and denote $\varphi_{\kappa}$ the density function of the
$\Gamma(2,\kappa)$-distribution defined by $\varphi_{\kappa}(x)=\kappa^2 x\e^{-\kappa
x}1_{]0,+\In[}$; for any $s\in\C$ such that $\Re \ s<\kappa$, the Laplace transform $\widehat{\varphi}_{\kappa}$ of $\varphi_{\kappa}$ exists and one gets $\widehat{\varphi}_{\kappa}(s)=\frac{\kappa^2}{(s-\kappa)^2}$. Consider the following matrice
\begin{gather*}
 \Phi_{n_1,\kappa}(\d x):=\Phi_{n_1}*\varphi_{\kappa}(\d x),
 \end{gather*}
 and $\mathfrak{L}(\Phi_{n_1,\kappa})$ its Laplace transform defined for $\mid\Re \ \l\mid\leq\a_0$. One gets the
\begin{property}\label{prop2.2}
 There exist $\delta_1,\delta_2,\varepsilon_1>0$, $0<\gamma<1$ and $\kappa>0$, such that for all $s\in[-\varepsilon_1,\varepsilon_1]$, $z\in\overline{K}(\delta_1,\delta_2)$,
\begin{equation}
\left\|\Phi_{n_1}(]-\In,-x]\cup[x,+\In[)\right\|=O(\e^{-\a_0 x}),\text{ for } x>0;
\end{equation}
\begin{equation}\label{eq36**}
   |z|^{n_1}\|\mathfrak{L}(\Theta_{k_1}^{\bullet m_1})(s)\|\leq\ga;
   \end{equation}
 \begin{equation}\label{eq36*}
   |z|^{n_1}\left\|\mathfrak{L}(\Phi_{n_1})(s)-\mathfrak{L}(\Phi_{n_1,\kappa})(s)\right\|\leq\frac{1-\ga}{2},\quad \text{for } -\a_0\leq s\leq\a_0.
 \end{equation}
\end{property}

\begin{proof}
 \begin{description}
  \item[1)] The first equality is derived from the fact that
  $$\left\|\int_0^{+\In}\e^{\a_0 x}\Phi_{n_1}(\d x)\right\|\leq\left\|\int_0^{+\In}\e^{\a_0 x}M^{\bullet n_1}(\d x)\right\|<+\In$$
  (resp. $\left\|\int_{-\In}^0\e^{-\a_0 x}\Phi_{n_1}(\d x)\right\|<+\In$).\\
  Therefore, for $x>0$,
  $$\left\|\Phi_{n_1}[x,+\In[\right\|\leq\left\|\e^{-\a_0 x}\int_{x}^{+\In}\e^{\a_0 t}\Phi_{n_1}(\d t)\right\|\leq C\e^{-\a_0 x}$$
 (and $\left\|\Phi_{n_1}[-\In,x[\right\|\leq C\e^{-\a_0 x},$ for $x<0$).
  \item[2)] The equalities (\ref{eq356}), (\ref{eq357}) and Lemma \ref{lem2.8} give, for $q\leq z\leq1$,
  $$z^{n_1}\Big\|P^{n_1}(\l_+(z))-\mathfrak{L}(\Phi_{n_1})(\l_+(z))\Big\|=z^{n_1}\Big\|\mathfrak{L}(\Theta_{n_1})(\l_+(z))\Big\|\leq z^{n_1}\Big\|\mathfrak{L}(\Theta_{k_1}^{\bullet m_1})(\l_+(z))\Big\|<1.$$
  Recall that $z\mapsto\l_+(z)$ is continuous on $[q,1]$ and $s\mapsto\|P^{n_1}(s)-\mathfrak{L}(\Phi_{n_1})(s)\|$
  is continuous on a neigborhood of $0$, we can then choose some suitable $\delta_1,\delta_2$, $\varepsilon_1>0$ and $0<\gamma<1$, such that (\ref{eq36**}) holds.
  \item[3)] The inequality (\ref{eq36*}) is an immediate consequence of the
  following  lemma, applied to the densities $\varphi_{n_1,i,j}(x)$ of $M^{\bullet n_1}_{i,j}$ for any $i,j\in E$.
 \end{description}
 \end{proof}

\begin{lemma}\label{prop4.2}
Fix $a<0<b$ and let $h:\R\rightarrow\R$ be a Borel function, such
that $\forall s\in[a,b]$,
$$\int_\R\e^{sx}|h(x)|\d x<+\In.$$
Set $h_\kappa=h*\varphi_\kappa$, where $h\ast\varphi_\kappa(x)=\int_0^{+\In} h(x+y)\varphi_\kappa(y)\d y$. Then
\begin{equation}\label{eqhphi}
 \lim_{\kappa\rightarrow+\In}\sup_{a\leq s\leq b}\int_\R\e^{sx}|h(x)-h_\kappa(x)|\d x=0.
\end{equation}
\end{lemma}
\begin{proof}
 We first prove that
 \begin{equation}\label{eqplus}
\lim_{y\rightarrow0}\sup_{a\leq s\leq b}\int_\R\e^{sx}|h(x+y)-h(x)|\d x=0.
\end{equation}

Indeed, fix $\varepsilon>0$ and choose a continuous function $\psi_{\varepsilon}$ with compact support $[\a,\b]$ such that
 \begin{equation}\label{star}
   \int(\e^{a t}+\e^{bt})|h(t)-\psi_{\varepsilon}(t)|\d t<\varepsilon.
 \end{equation}
For $a\leq s\leq b$ and $\mid y\mid\leq1$, one thus gets
 \begin{equation*}
  \int_\R\e^{sx}|h(x+y)-\psi_{\varepsilon}(x+y)|\d x\leq\e^{-ys}\int_{\R}(\e^{a t}+e^{bt})\mid h(t)-\psi_{\varepsilon}(t)\mid\d t\leq\e^{-ys}\varepsilon\leq(\e^{-a}+\e^b)\varepsilon.
 \end{equation*}
 Therefore,
 \begin{equation*}
  \begin{split}
   \int_\R\e^{sx}|h(x+y)-h(x)|\d x&\leq\int_\R\e^{sx}|h(x+y)-\psi_{\varepsilon}(x+y)|\d x+\int_\R\e^{sx}|\psi_{\varepsilon}(x+y)-\psi_{\varepsilon}(x)|\d x\\
   &\qquad\qquad\qquad+\int_\R\e^{sx}|\psi_{\varepsilon}(x)-h(x)|\d
x\\
&\leq 2(\e^{-a}+\e^{b})\var +\int_{\a-1}^{\b+1}(\e^{a x}+\e^{bx})|\psi_{\varepsilon}(x+y)-\psi_{\varepsilon}(x)|\d
x.
  \end{split}
 \end{equation*}
By the uniform continuity of $\psi_{\varepsilon}$ on $\R$, one gets
$|\psi_{\varepsilon}(x+y)-\psi_{\varepsilon}(x)|\stackrel{y\rightarrow0}{\longrightarrow}0$ uniformly on $\mathbb R$  and by
the dominated convergence theorem
$$\limsup_{y\rightarrow+\In}\sup_{a\leq s\leq b}\int_\R\e^{sx}|h(x+y)-h(x)|\d x\leq 2(\e^{-a}+\e^b)\varepsilon.$$
One can conclude since $\var$ is arbitrary.\\
We are now able to prove (\ref{eqhphi}). Since $\varphi_\kappa$ is a density, one gets
$$\int_\R\e^{sx}|h(x)-h_\kappa(x)|\d x\leq I_r(s,\kappa)+J_r(s,\kappa),$$ with
 $$\displaystyle I_r(s,\kappa):=\int_0^r\varphi_\kappa(y)\left(\int_\R\e^{sx}|h(x+y)-h(x)|\d x \right)\d y$$
and
$$\displaystyle  J_r(s,\kappa):=\int_r^{+\In}\varphi_\kappa(y)\left(\int_\R\e^{sx}|h(x+y)-h(x)|\d x \right)\d y.$$
Fix $\varepsilon>0$. By (\ref{eqplus}), one may choose $r$ small enough in such a way that, for $\mid y\mid\leq r$ and any $s \in [a, b]$
$$  \int_\R \e^{sx}\mid h(x+y)-h(x)\mid\d x\leq\varepsilon,$$
and since $\varphi_\kappa$ is a density of probability, one gets $\forall s\in[a,b], \forall\kappa>0,\quad I_r(s,\kappa)\leq\varepsilon$.

On the other hand,
\begin{eqnarray*}J_r(s,\kappa)&\leq&\int_r^{+\In}\e^{s  y}\varphi_\kappa(y)\Big(\int_\R\e^{st}\mid h(t)-h(t-y)\mid\d t\Big)\d y
\\
&\leq&\left[\int_r^{+\In}(1+\e^{\mid a\mid y})\varphi_\kappa(y)\d
y\right]\times \sup_{a\leq s\leq b}\left(\int_\R\e^{st}|h(t)|\d
 t\right).
 \end{eqnarray*}
Setting $u=\kappa y$, one obtains $\displaystyle \int_r^{+\In}\e^{\mid a\mid y}\varphi_\kappa(y)\d
y=\int_{r\kappa}^{+\In}u\e^{u(\frac{\mid a\mid}{\kappa}-1)}\d u$,
and so,  for $\kappa>2|a|$,
$$\int_r^{+\In}\e^{\mid a\mid y}\varphi_\kappa(y)\d y\leq\int_{r\kappa}^{+\In} u\e^{-\frac{u}{2}}\d u ; $$
then $\displaystyle \limsup_{\kappa\to +\infty}\sup_{s\in[a,b]}J_r(s,\kappa)=0$.
\end{proof}

We now introduce the following matrices,
\begin{gather*}
 B(z,\d x):=z^{n_1}\Big(M^{\bullet n_1}(\d x)-\Phi_{n_1,\kappa}(\d x)\Big),\\
 \widetilde{B}(z,\d x):=\sum_{k=1}^{+\In}B^{\bullet k}(z,\d x)
 \end{gather*}
and denote  $\mathfrak{L}(B)$ and $\mathfrak{L}(\widetilde{B})$ their Laplace transforms defined for $\mid\Re \ \l\mid\leq\a_0$.
\begin{lemma}\label{lem2.9}
 There exist $\delta_1,\delta_2$ and $\varepsilon>0$  such that
 \begin{enumerate}
   \item\label{as1} $\sup_{z\in \overline{K}(\delta_1,\delta_2)\atop|s|\leq\varepsilon}\left\|\int_\R\e^{su}\widetilde{B}(z,\d u)\right\|<+\In;$
   \item\label{as2} for $z\in \overline{K}(\delta_1,\delta_2)$, $|s|\leq\varepsilon$, $\theta\in\R$, the matrix $I-\mathfrak{L}(B)(z,s+i\theta)$ is invertible and
       $$\left(I-\mathfrak{L}(B)(z,s+i\theta)\right)^{-1}=I+\mathfrak{L}(\widetilde{B})(z,s+i\theta).$$
 \end{enumerate}
\end{lemma}
\begin{proof}
 \begin{description}
  \item[1)] For $z\in \overline{K}(\delta_1,\delta_2)$ and
  $|s|\leq\varepsilon$, we have
      \begin{equation*}
        \Big\|\L(B)(z,s)\Big\|\leq|z|^{n_1}\Big\|\L(\Theta_{n_1})(s)\Big\|+|z|^{n_1}\Big\|\L (\Phi_{n_1})(s)-\L(\Phi_{n_1,\kappa})(s)\Big\|.
      \end{equation*}
      From (\ref{eq357}), (\ref{eq36**}) and (\ref{eq36*}), there exist $\delta_1,\delta_2,\varepsilon>0$ and
      $0<\ga<1$ such that
$$\left\|\int_\R\e^{su}B(z,\d u)\right\|\leq\frac{1+\ga}{2}<1.$$
      Therefore, for any $z\in \overline{K}(\delta_1,\delta_2)$, $|s|\leq\varepsilon$,
     $$\Big\|\L(\widetilde{B})(z,s)\Big\|\leq\sum_{k\geq0}\Big\|\L(B)(z,s)\Big\|^k\leq\sum_{k\geq0}(\frac{1+\ga}{2})^k<+\In.$$
  \item[2)] By the first assertion, for any
  $z\in\overline{K}(\delta_1,\delta_2)$, $|s|\leq\varepsilon$ and
  $\theta\in\R$, the matrix
$$I-\mathfrak{L}(B)(z,s+i\theta)$$
 is invertible,
  with inverse
  $$\Big(I-\mathfrak{L}(B)(z,s+i\theta)\Big)^{-1}=\sum_{k=0}^{+\In}\mathfrak{L}(B^k)(z,s+i\theta)=I+\mathfrak{L}(\widetilde{B})(z, s+i \theta).$$
 \end{description}
\end{proof}

\subsection{The resolvent of $P(\l)$}\label{subsec2.3}
We denote by $V_N[-\a_0,\a_0]$ the algebra of $N\times N$ matrices whose terms are Laplace transforms of Radon measures $\sigma$ on $\R$, satisfying
$$\int_\R\e^{\l x}\d |\sigma|(x)<+\In,\quad\text{for } |\Re\ \l|\leq\a_0.$$


\begin{theorem}\label{thm2.4.2}
There exist $\delta_1,\delta_2$ and $\varepsilon>0$  such that
 \begin{description}
  \item[1)] The function $A(z,\l)$ defined by
   \begin{equation}\label{eq2.4.2}
      A(z,\l):=(I-zP(\l))^{-1}+\frac{\Pi_+(z)}{(\l-\l_+(z))\b_+(z)}+\frac{\Pi_-(z)}{(\l-\l_-(z))\b_-(z)}
       \end{equation}
       is analytic for $(z,\l)$ in the open set
       $$E(\de_1,\de_2,\varepsilon):=\{(z,\l);z\in K(\de_1,\de_2),\l\in S_z(\varepsilon)\},$$
        with $S_z(\varepsilon):=\{\l:\l_-(z,\varepsilon)<\Re\ \l<\l_+(z,\varepsilon)\}$, where $\b_{\pm}(z):=zk'(\l_\pm(z))$, $\Pi_\pm(z):=\Pi(\l_\pm(z))$, $\l_-(z,\varepsilon)=\Re\ \l_+(z)-\varepsilon$ and $\l_+(z,\varepsilon)=\Re\ \l_+(z)+\varepsilon$
  \item[2)] For $(z,\l)\in E(\de_1,\de_2,\varepsilon)$, one gets
   \begin{equation}\label{eq2.4.3}
    \begin{split}(I-zP(\l))^{-1}=I-&\frac{\Pi_+(z)}{(\l-\l_+(z))\b_+(z)}-\frac{\Pi_-(z)}{(\l-\l_-(z))\b_-(z)}+\\
                                   &\int1_{]-\In,0[}\,\e^{\l x}\d a_-(z,x)+\int1_{[0,+\In[}\,\e^{\l x}\d a_+(z,x),
     \end{split}
   \end{equation}
   where $a_+(z,\cdot)$ (resp. $a_-(z,\cdot)$) is a Radon measure on $ \R_+$ (resp. $\R_-$), with values in $V_{N\times N}[-\a_0,\a_0]$.\\
   Furthermore, for $x\geq0$ (resp. $x<0$), the function $z\mapsto a_+(z,x)$ (resp. $z\mapsto a_-(z,x)$) is analytic on
   $K(\de_1,\de_2)$, and satisfy :  for any  $ z\in\overline{K}(\de_1,\de_2)$
  :
   \begin{equation}\label{eq2.4.4}
    \|a_+(z,+\In)-a_+(z,x)\|\leq
    C\e^{-(\Re\ \l_+(z)\ +\varepsilon)x},\quad x\geq0,
   \end{equation}
   \begin{equation}\label{eq2.4.5}
    \|a_-(z,-\In)-a_-(z,x)\|\leq
    C\e^{-(\Re\ \l_-(z)\ -\varepsilon)x},\quad x<0.
   \end{equation}
 \end{description}
\end{theorem}

\begin{proof}
Throughout the present proof,  the parameters $\delta_1, \delta_2$ and $\varepsilon$ will satisfy the conclusions of Lemma \ref{lem2.9}.
\begin{description}
\item[1)]
As we mentioned in Remark \ref{remark2.2.1}, for $(z,\l)$ such hat $1-zk(\l)\neq0$, $|\Re\ \l|\leq\a_0$ and $|\Im\ \l|\leq\a_0$ (i.e; $\l \in \Delta_{\a_0}$), the operator $I-zP(\l)$ is invertible with inverse
$$(I-zP(\l))^{-1}=\frac{zk(\l)}{1-zk(\l)}\Pi(\l)+\sum_{n=0}^{+\In}z^nR^n(\l).$$
By the implicit function theorem, there exists real numbers $\de_1,\de_2>0$ such that when $z\in K(\de_1,\de_2)$,
the equation $1-zk(\l)=0$ has two distinct roots $\l_-(z)$ and $\l_+(z)$, given by
 \begin{equation}
 \l_\pm(z)=\pm \sqrt{\frac{2}{k''(0)}}\sqrt{1-z}\pm\frac{k^{(3)}(0)}{3(k''(0))^2}(1-z)+\sum_{k=3}^{+\In}(\pm1)^k\a_k( 1-z)^{k/2}.
 \end{equation}
 So we can choose $\de_1, \de_2$ and $\varepsilon$ such that
$\Re\ \l_-(z)\ -\varepsilon<\Re\ \l_+(z)\ +\varepsilon$ for any $z \in K(\delta_1, \delta_2)$. The residue of the map $\displaystyle\l\mapsto\frac{zk(\l)\Pi(\l)}{1-zk(\l)}$ at $\l_+(z)$ (resp. $\l_-(z)$) can be computed as
$$Res\left(\frac{zk(\l)\Pi(\l)}{1-zk(\l)},\l_\pm(z)\right)=-\frac{\Pi_\pm(z)}{\b_\pm(z)}.$$
Therefore, the function
$$(z,\l)\longmapsto\frac{zk(\l)\Pi(\l)}{1-zk(\l)}+\frac{\Pi_{+}(z)}{\b_+(z)(\l-\l_+(z))}+\frac{\Pi_{-}(z)}{\b_-(z)(\l-\l_-(z))}$$
is analytic for $(z,\l)\in E(\delta_1,\delta_2,\varepsilon)$.

Moreover,  $\displaystyle\sup_{|\Re\ \l|\leq\a_0}r(R(\l))<1$   ; the function $(z,\l)\longmapsto\sum_{n=0}^{+\In}z^nR^n(\l)$ is thus analytic on the domain $E(\delta_1,\delta_2,\varepsilon)$  when $\de_1, \de_2$ and $\var$ are small enough.

At last, by Theorem \ref{thm55} (2), one may choose $\a_0$ small enough in such a way
$$ \sup_{\stackrel{\vert {\rm Re}\l \vert \leq \a_0}{\vert {\rm Im}\l \vert \geq \a_0}}r(P(\l))<1
$$ which leads to the analyticity of the map $(\l, z)\mapsto  (I-zP(\l))^{-1}$ on the set $ \{(z, \l) \in E(\delta_1, \delta_2, \varepsilon)\slash \vert {\rm Im}\l \vert \geq \a_0\}$ ; the analyticity of the maps $(z, \l) \mapsto \frac{\Pi_{+}(z)}{\b_+(z)(\l-\l_+(z))}$ and $(z, \l) \mapsto \frac{\Pi_{-}(z)}{\b_-(z)(\l-\l_-(z))}$ on this domain also hold and
the proof of assertion 1) is achieved.

\item[2)] For $q\leq z<1$ and $\Re\ \l_-(z)<\Re\ \l<Re\ \l_+(z)$, one gets $zk(\Re\ \l)<1$; since $r(P(\l))\leq r(P(\Re\ \l))=k(Re\ \l)$, one thus obtains $zr(P(\l))<1$ for such a $z$ and so
\begin{equation}\label{eq2.4.10}
 (I-zP(\l))^{-1}=\sum_{n=0}^{+\In}z^nP^n(\l)=\left(\sum_{n=0}^{+\In}z^n\E_i(\e^{\l
 S_n},X_n=j)\right)_{i,j}.
\end{equation}

For every $(i,j)\in E\times E$, we consider the following distribution functions:
\begin{gather*}
\text{for }x\geq0,\qquad(a_+(z,x))_{i,j}:=\sum_{n=1}^{+\In}z^{n}\P_i(0\leq S_n<x,X_n=j)-\frac{(\Pi_+(z))_{i,j}}{\l_+(z)\b_+(z)}(1-\e^{-\l_+(z)x});\\
\text{for }x < 0,\qquad(a_-(z,x))_{i,j}:=\sum_{n=1}^{+\In}z^{n}\P_i(x\leq S_n<0,X_n=j)+\frac{(\Pi_-(z))_{i,j}}{\l_-(z)\b_-(z)}(1-\e^{-\l_-(z)x}).
\end{gather*}
The measures $a_+(z,x)$ and $a_-(z,x)$ satisfy the following identities
$$\int1_{[0,+\In[}(x)\e^{\l x}\d (a_+(z,x))_{i,j}=\sum_{n=1}^{+\In}z^n\E_i(\e^{\l S_n},S_n\geq0,X_n=j)+\frac{(\Pi_+(z))_{i,j}}{(\l-\l_+(z))\b_+(z)},$$
$$\int1_{]-\In,0[}(x)\e^{\l x}\d (a_-(z,x))_{i,j}=\sum_{n=1}^{+\In}z^n\E_i(\e^{\l S_n},S_n<0,X_n=j)+\frac{(\Pi_-(z))_{i,j}}{(\l-\l_-(z))\b_-(z)}.$$
Summing the two precedent equalities and using (\ref{eq2.4.10}), we
find the expected formula (\ref{eq2.4.3}).\\
Now we prove the analyticity of the functions $z\longmapsto a_+(z,\cdot)$ and
$z\longmapsto a_-(z,\cdot)$. By (\ref{eq2.4.2}) and (\ref{eq2.4.3}),
we get
$$A(z,\l)=I+\int1_{[0,+\In[}(x)\e^{\l x}\d a_{+}(z,x)+\int1_{]-\In,0[}(x)e^{\l x}\d a_-(z,x).$$

Observe that the function $x\mapsto a_+(z,x)$ is continuous and vanishes at $x=0$ ;
applying the inversion formula for the Laplace integral transform (\cite{Widd}), we obtain for $x\geq 0$ and $0<\delta< \Re \ \l_+(z)$,
\begin{equation}\label{eq2.4.11}
 \begin{split}
  a_+(z,+\In)-a_+(z,x)&=\sum_{n=1}^{+\In}z^n\P_i(S_n\geq x,X_n=j)-\frac{\Pi(\l_+(z))\e^{-\l_+(z)x}}{\l_+(z)\b_+(z)}\\
                     &=\frac{1}{2\pi i}\int_{\Re\ \l\ =\ \de}\e^{-\l x}\frac{A(z,\l)}{\l}\,\d\l.
 \end{split}
\end{equation}
On the other hand, the function $(z,\l)\mapsto A(\l,z)$ is analytic on the set $E(\de_1,\de_2,\varepsilon)$ and by Cauchy's theorem, one gets
\begin{equation*}
 \begin{split}
  a_+(z,+\In)-a_+(z,x)&=\frac{1}{2\pi
  i}\int_{Re\ \l\ =\ \Re\ \l_+(z)\ +\varepsilon}\e^{-\l x}\frac{A(z,\l)}{\l}
  \d\l\\
  &=\frac{1}{2\pi}\e^{-(\Re\ \l_+(z)\ +\varepsilon)x}\int_\R\e^{-ix\theta}\frac{A(z,\Re\ \l_+(z)\ +\varepsilon+i\theta)}{\Re\ \l_+(z)+\ \varepsilon+i\theta}\,\d\theta.
 \end{split}
\end{equation*}
To compute this last integral, we use the following
\begin{lemma}\label{lem2.4.3}
 Let $a\neq b$ two complex numbers such that $\Re\ a >0$ and $ \Re\ b>0$. For $x\geq 0$, one gets
$$
  \int_{-\In}^{+\In}\frac{\e^{ix\theta}}{(i\theta-a)(i\theta-b)}\d
  \theta=0.
$$
\end{lemma}

By (\ref{eq2.4.2}) and Lemma \ref{lem2.4.3}, one gets for $x\geq 0$,
\begin{equation*}
 \begin{split}
a_+(z,+\In)-a_+(z,x)&=\frac{1}{2\pi}\e^{-(z_+(z,\varepsilon))x}\int_\R\frac{\e^{-ix\theta}[I-zP(z_+(z,\varepsilon)+i\theta)]^{-1}}{z_+(z,\varepsilon)+i\theta}\,\d\theta\\
&=\frac{1}{2\pi}\e^{-(\Re\ \l_+(z)\ +\varepsilon)x}\ W_+(z,\varepsilon,x)
 \end{split}
\end{equation*}
with
$$W_+(z,\varepsilon,x):=\int_\R\frac{\e^{-ix\theta}[I-zP(z_+(z,\varepsilon)+i\theta)]^{-1}}{z_+(z,\varepsilon)+i\theta}\,\d\theta.$$

 By a similar argument, one may write for $x<0$,
$$
a_-(z,-\In)-a_-(z,x)=\frac{1}{2\pi}\e^{-(\Re\ \l_-(z)\ -\varepsilon)x}\ W_-(z,\varepsilon,x)$$
with
$$W_-(z,\varepsilon, x)=\int_\R\frac{\e^{-ix\theta}[I-zP(\Re\ \l_-(z)-\varepsilon+i\theta)]^{-1}}{\Re\ \l_-(z)-\varepsilon+i\theta}\,\d\theta. $$
Note that by definition of  $a_\pm$,  the functions $x \mapsto W_\pm(z, \var, x)$ are left-continuous, for any $z \in K(\de_1, \de_2)$.
One completes the proof  by a simple application of the following
   :
\begin{property}\label{TDFdesresolvantes}
We fix $\varepsilon>0$ and $\delta_1, \delta_2>0$ small enough in such a way   the conclusions of  Lemma \ref{lem2.9} hold  for any
$
 z \in \overline K(\delta_1, \delta_2).
$
We set
\begin{description}
\item[$\bullet$] $\l_{\pm}(z,\varepsilon)=\Re\ \l_{\pm}(z)\pm\varepsilon;$
  \item[$\bullet$]$\displaystyle W_+(z, \varepsilon, x)=\int_\R\e^{-i\theta
 x}\frac{[I-zP(\l_+(z,\varepsilon)+i\theta)]^{-1}}{\l_+(z,\varepsilon)+i\theta}\;\d
 \theta,\qquad \text{for }x\geq 0;$
 \item[$\bullet$] $\displaystyle W_-(z, \varepsilon, x )=\int_\R\e^{-i\theta
 x}\frac{[I-zP(\l_-(z,\varepsilon )+i\theta)]^{-1}}{\l_-(z,\varepsilon)+i\theta}\;\d
 \theta,\qquad \text{for }x<0$.
\end{description}
 Then, there  exists a constant $C=C(\varepsilon)>0$ such that for $x\geq0$ (resp. $x<0$), one gets
\begin{equation}\label{4}
\forall z\in\overline{K}(\de_1,\de_2),\quad\|W_+(z,x,\varepsilon)\|\leq C \quad(\text{resp. } \|W_-(z,x,\varepsilon)\|\leq C).
\end{equation}
\end{property}
\end{description}
\end{proof}
\begin{proof} Note first that by the choice of the constants $\varepsilon_1, \varepsilon_2$ and $\delta_1$, one gets $\vert \l_{\pm}(z,\varepsilon)\vert \leq \varepsilon_1$ for any $z\in\overline{K}(\de_1,\de_2)$.


 For $z\in\overline{K}(\de_1,\de_2)$, the matrices $I-z^{n_1}P^{n_1}(s)$ and $I-\mathfrak{L}(B)(z,s)$ are invertible; the identity
$$z^{n_1}P^{n_1}(s)=\mathfrak{L}(B)(z,s)+z^{n_1}\mathfrak{L}(\Phi_{n_1,\kappa})(s)$$
allows us to write
 $$
   [I-z^{n_1}P^{n_1}(s)]^{-1}=[I-\mathfrak{L}(B)(z,s)]^{-1}+[I-z^{n_1}P^{n_1}(s)]^{-1}z^{n_1}\mathfrak{L}(\Phi_{n_1,\kappa})(s)[I-\mathfrak{L}(B)(z,s)]^{-1}.\quad^(\footnotemark[1]^)
$$\footnotetext[1]{We use the classical fact that for any $N\times N$ matrices $U$  and $V$  such that  $I-U$ and $I-V$  are invertible, setting  $W=U-V$,  one has
$(I-U)^{-1}=(I-V)^{-1}+(I-U)^{-1}W(I-V)^{-1}$.  We apply this identity to  $U=z^{n_1}P^{n_1}(s)$  and  $V=\mathfrak{L}(B)(z,s)$.}
Throughout this proof, in order to simplify the notations, we set $\blacktriangle:=\l_+(z,\varepsilon)+i\theta$, so that
 \begin{equation*}
  \begin{split}
   [I-zP(\blacktriangle)]^{-1}=&\ I+[zP(\blacktriangle)+\cdots+z^{n_1}P^{n_1}(\blacktriangle)][I-z^{n_1}P^{n_1}(\blacktriangle)]^{-1}\\
                     =&\ I+[zP(\blacktriangle)+\cdots+z^{n_1}P^{n_1}(\blacktriangle)][I-\mathfrak{L}(B)(z,\blacktriangle)]^{-1}\\
                     &\quad+[zP(\blacktriangle)+\cdots+z^{n_1}P^{n_1}(\blacktriangle)][I-z^{n_1}P^{n_1}(\blacktriangle)]^{-1}z^{n_1}\mathfrak{L}(\Phi_{n_1,\kappa})(\blacktriangle)[I-\mathfrak{L}(B)(z,\blacktriangle)]^{-1}
  \end{split}
 \end{equation*}
and we may decompose $W_+(z, \var, x)$ as  $W_+(z, \var, x)=W_{+1}(z, \var, x)+W_{+2}(z, \var, x)+W_{+3}(z, \var, x)$ with
 $$W_{+1}(z,\varepsilon, x):=\int_\R\frac{\e^{-i\theta
 x}I}{\blacktriangle}\;\d\theta,$$
$$W_{+2}(z,\varepsilon, x):=\int_\R\frac{\e^{-i\theta x}[zP(\blacktriangle)+\cdots+z^{n_1}P^{n_1}(\blacktriangle)][I-\mathfrak{L}(B)(z,\blacktriangle)]^{-1}}{\blacktriangle}\;\d\theta,$$
$$W_{+3}(z,\varepsilon, x):=\int_\R\e^{-i\theta
                      x}\frac{z^{n_1}[zP(\blacktriangle)+\cdots+z^{n_1}P^{n_1}(\blacktriangle)][I-z^{n_1}P^{n_1}(\blacktriangle)]^{-1}\mathfrak{L}(\Phi_{n_1,\kappa})(\blacktriangle)[I-\mathfrak{L}(B)(z,\blacktriangle)]^{-1}}{\blacktriangle}\;\d
                      \theta.$$
The fact that $W_{+1}(z,\varepsilon, x)$ is bounded uniformly in $z \in \overline K(\delta_1, \delta_2)$  and $x\geq 0$ is a direct consequence of  the following Lemma; indeed, one gets
$\displaystyle
\int_\R {e^{-i\th x}\over \blacktriangle }\d\theta =    \pi (1-{\rm sgn}(x)) e^{-  \l_+(z, \var)x }=0 $,
 since $x \geq 0$.

\begin{lemma}\label{intcomplexe'}
For any $ a >0$ and any $x \in \R$ one gets $\displaystyle
\int_\R{e^{i\theta x}\over a+i\theta} \d\theta=  \pi e^{-ax}(1+{\rm sgn}(x)).
$
\end{lemma}

 Now, we focuse our attention on the term $W_{+2}(z,\varepsilon)$. By Lemma \ref{lem2.9}, the function
 $z\mapsto[zP(\blacktriangle)+\cdots+z^{n_1}P^{n_1}(\blacktriangle)][I-\mathfrak{L}(B)(z,\blacktriangle)]^{-1}$ is
 the Laplace transform  at point
 $\blacktriangle$ of the measure $\mu(z,\d x)=[zM(\d x)+\cdots+z^{n_1}M^{n_1}(\d x)]\bullet\widetilde{B}(z,\d x)$.
  By the definition of $P$ and Lemma \ref{lem2.9}, for $z\in [q+\delta_1,1+\delta_2]$, the term $\mu(z,\cdot)$ is a matrix of finite measures on $\R$, so we get
 $$\sup_{z\in\overline{K}(\delta_1,\delta_2)}\|[zM(\R)+\cdots+z^{n_1}M^{n_1}(\R)]\widetilde{B}(\R)\|<+\In.$$
 By the inversion formula for the Laplace integral transform, for any continuity point $x\geq 0$ of the map $t\mapsto\mu(z,[t,+\In[)$, one gets
 \begin{equation}\label{6}
  \e^{-\l_+(z,\varepsilon)x}W_{+2}(z,\var, x)=\mu(z,[x,+\In[).
 \end{equation}
 This equality holds in fact for any $x \geq 0$ since the two members are left-continous on $\R$.
  Therefore, for any $ x \geq 0$, one gets
  \begin{equation*}
    \|W_{+2}(z,\var, x)\|=\|\e^{\l_+(z,\varepsilon)x}\mu(z,[x,+\In[)\|\leq\int_{-\In}^{+\In}e^{\l_+(z,\varepsilon)t}\|zM(\d t)+\cdots+z^{n_1}M^{n_1}(\d
                                       t)]\bullet\widetilde{B}(z,\d t)\|.
  \end{equation*}
  Using Lemma \ref{lem2.9}
  and the fact that
  $\displaystyle \sup_{z\in\overline{K}(\de_1,\de_2)}\|P(\l_+(z,\varepsilon))\|<+\In$, we obtain immediately
  $$
  \sup_{\stackrel{z\in\overline{K}(\de_1,\de_2)}{x\geq 0}}\|W_{+2}(z, \var, x)\|<+\In.
  $$
 We finally  study the last term $W_{+3}(z,x)$. One gets
  $\displaystyle \|\mathfrak{L}(\Phi_{n_1,\kappa})(\blacktriangle)\|=\frac{\kappa^2}{|\blacktriangle-\kappa |^2}\|\mathfrak{L}(\Phi_{n_1})(\blacktriangle)\|,$
 with
  $$\sup_{z\in\overline{K}(\de_1,\de_2)}\|\mathfrak{L}(\Phi_{n_1})(\blacktriangle)\|\leq\|P(\l_+(z,\varepsilon))\|^{n_1}<+\In.$$
  On the other hand, by Lemma \ref{lem2.9} one gets $\displaystyle \sup_{z\in\overline{K}(\de_1,\de_2)}\|[I-\mathfrak{L}(B)(z,\blacktriangle)]^{-1}\|<\In$. Since the
  matrices
  $ [I-z^{n_1}P^{n_1}(\blacktriangle)]^{-1} $ and
  $ zP(\blacktriangle)+\cdots+z^{n_1}P^{n_1}(\blacktriangle) $ are clearly bounded in   $z\in \overline{K}(\delta_1,\delta_2)$, there finally exists a constant $C>0$ such that
   $$\forall z\in\overline{K}(\de_1,\de_2), \forall x\geq 0,\qquad \|W_{+3}(z,\varepsilon, x)\|\leq
C\sup_{z\in\overline{K}(\de_1,\de_2)}\int_\R\frac{1}{|\blacktriangle|}\times\frac{\kappa^2}{|\kappa-\blacktriangle|^2}\d
\theta<+\In.$$
\end{proof}

It remains to prove Lemmas \ref{lem2.4.3} and   \ref{intcomplexe'}.
\begin{proof} [Proof of Lemma \ref{lem2.4.3}]
 For $z\in\C$ and $x\geq 0$, set $f(x,z):=\frac{\e^{xz}}{(z-a)(z-b)}$ ; one gets
  \begin{equation}\label{eq47*}
  \int_{\ga_1\cup\ga_2\cup\ga_3\cup\ga_4}f(x,z)\d z=0,
 \end{equation}
  where $\ga_k$, $1\leq k\leq4$, are the paths defined as follows (see  Figure \ref{fig-dessin1}) : for $\a$, $A>0$
  \begin{alignat*}{2}
  \ga_1=\{z=i\theta;-A\leq\theta\leq A\},&\qquad \ga_2=\{z=-t+iA;0\leq
  t\leq\a\},\\
  \ga_3=\{z=-\a-i\theta, -A\leq\theta\leq A\},&\qquad \ga_4=\{z=t-iA;-\a\leq
  t\leq0\}.
  \end{alignat*}
  \begin{figure}
\begin{center}
  \includegraphics[width=300pt,height=240pt]{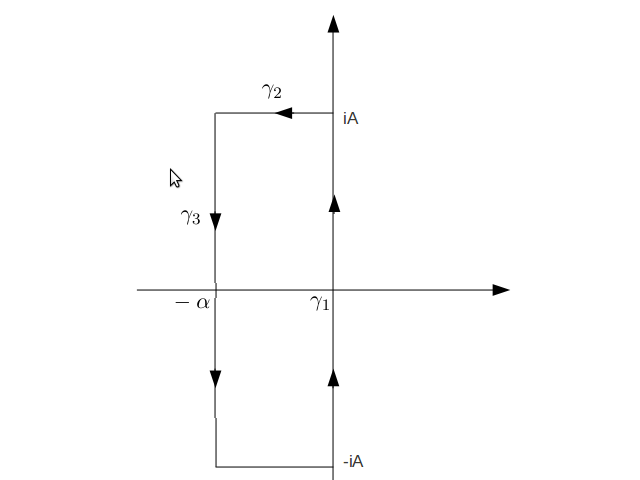}
\end{center}
 \caption{The closed path $\ga_1\cup\ga_2\cup\ga_3\cup\ga_4$ of Lemma \ref{lem2.4.3}.}\label{fig-dessin1}
\end{figure}

 In addition,
  \begin{equation*}
    \begin{split}
    \left|\int_{\ga_2}f(x,z)\,\d z\right|&\leq\int_0^{\a}\left|\frac{\e^{(-t+iA)x}}{(-t+iA-a)(-t+iA-b)}\right|\,\d
    t\\
                        &=\int_{0}^{\a}\frac{\e^{-tx}\d t}{\sqrt{(t+\Re\,a)^2+(A-\Im\,a)^2}\sqrt{(t+\Re\, b)^2+(A-\Im\,b)^2}}\\
                        &\leq\frac{\a} {\sqrt{(\Re\,a)^2+(A-\Im\,a)^2}\sqrt{(\Re\,b)^2+(A-\Im\, b)^2}}\stackrel{A\rightarrow+\In}{\longrightarrow}0.
   \end{split}
  \end{equation*}
 The same argument leads to
  $$\left|\int_{\ga_4}f(z)\,\d z\right|=\left|\int_{-\a}^{0}\frac{\e^{(t-iA)x}}{(t-iA-a)(t-iA-b)}\,\d
   t\right|\stackrel{A\rightarrow+\In}{\longrightarrow}0.$$
 On the other hand,
 \begin{equation*}
  \begin{split}
  \Big|\int_{\ga_3}f(z)\,\d z\Big|&\leq\e^{-\a x}\int_{-A}^{A}\frac{\d\th}{|\a+i\theta+a||\a+i\theta+b|}\\
                                &\leq\e^{-\a x}\int_{-\In}^{+\In}\frac{\d\th}{\sqrt{(\a+\Re\,a)^2+(\theta-\Im\,a)^2}\sqrt{(\a-\Re\,b)^2+(\theta-\Im\,b)^2}}\stackrel{\a\rightarrow+\In}{\longrightarrow}0.
 \end{split}
  \end{equation*}
Then $\displaystyle\lim_{A\rightarrow+\In}\int_{\ga_1}f(x,z)\,\d z=\int_{-\In}^{+\In}\frac{\e^{i\theta x}}{(i\theta-a)(i\theta-b)}\,\d\theta=0.$
\end{proof}

\begin{proof} [Proof of Lemma \ref{intcomplexe'}]

 For $z\in\C$ and $x\in \R$, set $g(x,z):=\frac{\e^{xz}}{z}$. For any fixed $x>0$, one gets

\begin{figure}
\begin{center}
  \includegraphics[width=450pt,height=300pt]{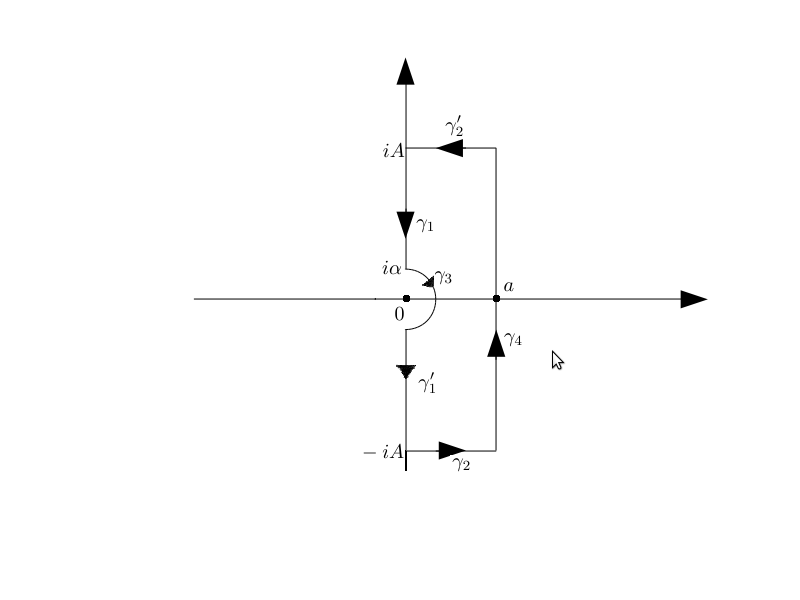}
\end{center}
 \caption{The closed path $\ga_1\cup\ga'_1\cup\ga_2\cup\ga'_2\cup\ga_3\cup\ga_4$ of Lemma \ref{intcomplexe'}.}\label{fig-dessin22}
\end{figure}
  \begin{equation}\label{eq47*}
  \int_{\ga_1\cup\ga'_1\cup \ga_2\cup \ga'_2\cup\ga_3\cup\ga_4}g(x,z)\d z=0,
 \end{equation}
  where $\ga_k$, $1\leq k\leq6$, are the paths defined as follows (see  Figure \ref{fig-dessin22}):  for $A>\a>0$
 \begin{itemize}
 \item $\gamma_1$ is the oriented segment from $iA$ to $i\a$
 \item $\gamma'_1$ is the oriented segment from $-i\a$ to $-iA$
 \item $\gamma_2$ is the oriented segment from $-iA$ to $a-iA$
 \item $\gamma'_2$ is the oriented segment from $a+iA$ to $iA$
 \item $\gamma_3$ is the clockwise oriented arc of circle from $i\a$ to $-i\a$
\item $\gamma_4$ is the oriented segment from $a-iA$ to $a+iA$
 \end{itemize}


 One gets
  \begin{enumerate}
 \item  $\displaystyle \int_{\ga_1\cup\ga'_1}g(x,z)\d z= -2i   \int_{\a }^{A} {\sin t x \over t} \d t \stackrel{\begin{array}{c}   \alpha \to 0   \\  A\to +\infty\end{array}}{\longrightarrow} -i\pi \ {\rm sgn}(x),$
 \item $ \displaystyle \Big\vert \int_{\ga_2\cup\ga'_2}g(x,z)\d z\Big\vert \leq 2{e^{ax}\over A}\stackrel{A\to +\infty}{\longrightarrow} 0,$
 \item
 $ \displaystyle  \int_{\ga_3}g(x,z)\d z =-i \int_{-{\pi \over 2}}
^ {\pi \over 2}e^{x\alpha e^{i\theta}} \d\theta
\stackrel{\a\to  0}{\longrightarrow} -i\pi$
\end{enumerate}
and equality (\ref{eq47*}) thus implies
$$
\int_{\gamma_4}g(x,z)\d z=ie^{ax} \int_{-A}^{ A}{e^{ix\theta}\over a+i\theta}\d \theta\quad \stackrel{A\to +\infty}{\longrightarrow} \quad  i\pi  (1+{\rm sgn}(x))
$$ and the Lemma follows.
\end{proof}

\section{On the factorization of $I-zP(\l)$}

\subsection{Preliminaries and motivation}
We first introduce the two following stopping times, which correspond to the first entrance time of the random walk $(S_n)_{n \geq 1}$ inside one of the semi-group $\R^+, \R^{*+}, \R^-$ and $\R^{*-}$ :
$$T_+= \inf\{n\geq1,\,S_n\geq0\};\qquad T^*_+= \inf\{n\geq1,\,S_n>0\};$$
 $$ T_-= \inf\{n\geq1,\,S_n\leq0\};\qquad T^*_-= \inf\{n\geq1,\,S_n<0\}.$$
Recall that $V_N[-\a_0,\a_0]$ is the algebra of $N\times N$ matrices whose terms are Laplace transforms of Radon measures $\sigma$ on $\R$, satisfying $\displaystyle\int_\R\e^{\l x}\d |\sigma|(x)<+\In$, for $|\Re\ \l|\leq\a_0$.
Let $G\in V_N[-\a_0,\a_0]$, defined by
$$G(\l)= \left(\int_\R\e^{\l x}\d \sigma_{i,j}(x)\right)_{1\leq i,j\leq N}.$$
For $|\Re\ \l|\leq\a_0$, we set $^($\footnote{the letter $\mathcal N$ corresponds to the restriction of the Radon measure to the {\it negative} or {\it strictly negative } half line $\R^-$ or  $\R^{*-}$ and the letter  $\mathcal P$  corresponds to the  {\it positive} or {\it strictly positive } half line $\R^+$ or  $\R^{*+}$}$^)$
\begin{gather*}
\mathcal{N}G(\l)= \Big(\int_{]-\In,0]}\e^{\l x}\d \sigma_{i,j}(x)\Big)_{i,j},\quad \mathcal{N}^*G(\l)= \Big(\int_{]-\In,0[}\e^{\l x}\d \sigma_{i,j}(x)\Big)_{i,j};\\
\daP G(\l)= \Big(\int_{[0,+\In[}\e^{\l x}\d \sigma_{i,j}(x)\Big)_{i,j},\quad
\daP^*G(\l)= \Big(\int_{]0,+\In[}\e^{\l x}\d \sigma_{i,j}(x)\Big)_{i,j}.
\end{gather*}

For $|z|<1$, we consider the following matrices of measures on $\R$:
$$ B_z(\d y)= \left(\sum_{n= 1}^{+\In}z^{n}\P_{i}\{S_1\geq S_n,S_2\geq S_n,\cdots,S_{n-1}\geq S_n,S_n\in\d y,X_n= j\}\right)_{i,j},$$
$$B^*_z(\d y)= \left(\sum_{n= 1}^{+\In}z^{n}\P_{i}\{S_1>S_n,S_2>S_n,\cdots,S_{n-1}> S_n,S_n\in\d y,X_n= j\}\right)_{i,j},$$
$$C_z(\d y)= \left(\sum_{n= 1}^{+\In}z^{n}\P_{i}\{S_1\geq0,S_2\geq0,\cdots,S_{n-1}\geq0,S_n\in\d y,X_n= j\}\right)_{i,j},$$
$$C^*_z(\d y)= \left(\sum_{n= 1}^{+\In}z^{n}\P_{i}\{S_1>0,S_2>0,\cdots,S_{n-1}>0,S_n\in\d y,X_n= j\}\right)_{i,j}.$$
For $\Re\ \l=0$, the related Laplace transforms of the above measures, denoted respectively by $B_z(\l)$, $B^*_z(\l)$, $C_z(\l)$ and $C^*_z(\l)$, are defined as following:

\begin{gather*}
 B_z(\l)=\int_{-\In}^{+\In}\e^{\l y}B_z(\d y),\quad B^*_z(\l)=\int_{-\In}^{+\In}\e^{\l y}B^*_z(\d y);\\
 C_z(\l)=\int_{-\In}^{+\In}\e^{\l y}C_z(\d y),\quad C^*_z(\l)=\int_{-\In}^{+\In}\e^{\l y}C^*_z(\d y).
\end{gather*}

 Note that the series which appear in these formulas do converge for $|z|<1$ and that the matrices$B_z(\l)$, $B^*_z(\l)$, $C_z(\l)$ and $C^*_z(\l)$ belong to $V_N[0,0]$.

 Let us now explain briefly how we will use these waiting times to prove the local limit theorem for the process  $m_n:= \min(0,S_1,\cdots,S_n)$. Indeed, the
Laplace transform of $m_n$ may be expressed in terms of the operators $\mathcal N^* $ and $  \mathcal P$ and the matrices
$ B^*_z $ and $C_z$ ; we have the

\begin{lemma}\label{lem4.1}
For $\l>0$ and $|z|<1$,
\begin{equation}\label{eq25}
 \sum_{n= 0}^{+\In}z^n\E_i(\e^{\l m_n};X_n= j)= \{[I+\daN^* B^*_z(\l)][I+\daP C_z(0)]\}_{i,j}.
\end{equation}
 \end{lemma}
\begin{proof}
Applying Markov property to the process $(X_n,S_n)$, we get
\begin{equation*}
  \begin{split}
   &\sum_{n= 0}^{+\In}z^n\E_i(\e^{\l m_n};X_n= j)\\
   = &\sum_{n= 0}^{+\In}z^n\sum_{k= 0}^n \E_i(\e^{\l S_k};S_0>S_k,\cdots,S_{k-1}>S_k,S_{k+1}\geq S_k,\cdots,S_n\geq S_k,X_n= j)\\
                                             = &\sum_{n= 0}^{+\In}z^n\sum_{k= 0}^n\sum_{l\in E}\E_i(\e^{\l S_k};S_1>S_k,\cdots,S_{k-1}>S_k,S_k<0,X_k= l)\;\E_l(S_1\geq0,\cdots,S_{n-k}\geq0,X_{n-k}= j)\\
                                             = &\sum_{l\in E}\left[\sum_{k= 0}^{+\In}z^k\E_i(\e^{\l S_k};S_1>S_k,\cdots,S_{k-1}>S_k,S_k<0,X_k= l)\right]\left[\sum_{p= 0}^{+\In}z^p\E_l(S_1\geq0,\cdots,S_p\geq0,X_p= j)\right]\\
                                             = &\Big\{[I+\daN^* B^*_z(\l)][I+\daP C_z(0)]\Big\}_{i,j}.
  \end{split}
 \end{equation*}
\end{proof}
  We will have to study the regularity with respect to $z$ and $\lambda$ of each factor $I+\daN^* B^*_z(\l)$ and $I+\daP C_z(0)$ ; to do this, we will use a classical approach based on the so-called {\it Wiener-Hopf factorization}.

 \subsection{The initial probabilistic factorization}
 We have the
\begin{proposition}\label{initfactorization}
 For $\Re\ \l= 0$ and $|z|<1$, one gets

 \begin{equation}\label{eq2.5}
  I-zP(\l)= (I-\daP B^*_z(\l))(I-\daN^* C_z(\l)),
\end{equation}

 \begin{equation}\label{eq2.6}
 (I-\daP B^*_z(\l))^{-1}= I+\daP C_z(\l),
\end{equation}

\begin{equation}\label{eq2.7}
  (I-\daN^* C_z(\l))^{-1}= I+\daN^* B^*_z(\l).
 \end{equation}
\end{proposition}

\begin{proof}
We first  check that
  \begin{equation}\label{eq2.8}
 (I-\daN^* C_z(\l))(I-zP(\l))^{-1}= I+\daP C_{z}(\l),
 \end{equation}
 and (\ref{eq2.5}) will follow by (\ref{eq2.6}). Note that, for $\Re\ \l= 0$, $r(P(\l))\leq r(P(0))= 1$. So for $|z|<1$, $(I-zP(\l))$ is invertible, with inverse
 $$(I-zP(\l))^{-1}= I+\sum_{n= 1}^{+\In}z^nP^n(\l).$$
By the definition of $P(\l)$ and the strong Markov property, we get
 \begin{equation*}
  \begin{split}
   &\delta_{i,j}+\sum_{n= 1}^{+\In}z^n(P^n(\l))_{i,j}\\
   &= \delta_{i,j}+\sum_{n= 1}^{+\In}\E_i(z^n\e^{\l S_n};X_n= j)\\
                                                  &= \delta_{i,j}+\E_i\left(\sum_{n= 1}^{T^*_--1}z^n\e^{\l S_n}; X_n= j\right)+\E_i\left(\sum_{n= T^*_-}^{+\In}z^n\e^{\l S_n}; X_n= j\right)\\
                                                  &= \delta_{i,j}+\E_i\left(\sum_{n= 1}^{+\In}z^n\e^{\l S_n};T^*_-\geq n+1; X_n= j\right)+\E_i\left\{z^{T^*_-}\e^{\l S_{T^*_-}}\left[\E_{X_{T^*_-}}\left(\sum_{n= 0}^{+\In}z^n\e^{\l S_n}; X_n= j\right)\right]\right\}\\
                                                  &= \delta_{i,j}+(\daP C_z(\l))_{i,j}+\sum_{l\in E}\left\{\left[\E_i\left(\sum_{k= 1}^{+\In}z^k\e^{\l S_k};T^*_-= k;X_k= l\right)\right]\left[\sum_{n= 0}^{+\In}\E_l\left(z^n\e^{\l S_n};X_n= j\right)\right]\right\}\\
                                                  &= \delta_{i,j}+(\daP C_z(\l))_{i,j}\\
                                                  &\qquad+\sum_{l\in E}\left\{\left[\sum_{k= 1}^{+\In}z^k\E_i(\e^{\l S_k};S_1\geq0,S_2\geq0,\cdots,S_{k-1}\geq0,S_k<0;X_k= l)\right]\left[\delta_{l,j}+\sum_{n= 1}^{+\In}\E_l(z^n\e^{\l S_n};X_n= j)\right]\right\}\\
   &= \delta_{i,j}+(\daP C_z(\l))_{i,j}+\sum_{l\in E}\left(\daN^* C_z(\l)\right)_{i,l}((I-zP(\l))^{-1})_{l,j}\\
                                                  &= \delta_{i,j}+(\daP C_z(\l))_{i,j}+\left(\daN^*
C_z(\l)(I-zP(\l))^{-1}\right)_{i,j}.
 \end{split}
\end{equation*}

We  now
 prove (\ref{eq2.6}) (and the proof of (\ref{eq2.5}) will be complete, as we claimed above). Set $F_z(\l)= (I-\daP B^*_z(\l))(I+\daP C_z(\l))$ ; we want to check that $F_z(\l)= I$. One gets
 \begin{equation}\label{eq2.9}
  \Big(F_z(\l)\Big)_{i,j}= \delta_{i,j}+\Big(\daP C_z(\l)\Big)_{i,j}-\Big(\daP B^*_z(\l)\Big)_{i,j}-\Big(\daP B^*_z(\l)\daP C_z(\l)\Big)_{i,j}.
\end{equation}
By the strong Markov property, we get
\begin{equation}
\begin{split}
&\Big(\daP B^*_z(\l)\daP C_z(\l)\Big)_{i,j}\\
= &\sum_{n= 1}^{+\In}z^n\E_i\left[\e^{\l S_n};
S_1> S_n,S_2>S_n,\cdots,S_{n-1}>
S_{n}\geq0;\E_{X_n}\left(\sum_{k= 1}^{+\In}z^k\e^{\l
S_k};S_1\geq0,\cdots,S_k\geq0,X_k= j\right)\right]\\
= &\sum_{n\geq1,k\geq1}z^{n+k}\E_i[\e^{\l S_{n+k}};S_1>
S_n,\cdots,S_{n-1}> S_n,S_{n+1}\geq S_n,\cdots,S_{n+k}\geq S_n\geq0,X_{n+k}= j]\\
= &\sum_{m= 2}^{+\In}z^m\left[\sum_{n= 1}^{m-1}\E_i(\e^{\l S_m};S_1>
S_n,\cdots,S_{n-1}> S_n,S_{n+1}\geq S_n,\cdots,S_m\geq S_n\geq0,X_m= j)\right].
\end{split}
\end{equation}
Therefore,

\begin{equation*}
 \begin{split}
  (F_z(\l))_{i,j}= \ &\delta_{i,j}+\sum_{m= 1}^{+\In}z^m\E_i(\e^{\l
  S_m};S_1\geq0,S_2\geq0,\cdots,S_m\geq0,X_m= j)\\
  &-\sum_{m= 1}^{+\In}z^m\E_i(\e^{\l S_m};S_1> S_m,S_2> S_m,\cdots, S_{m-1}> S_m
  \geq0;X_m= j)\\
  &-\sum_{m= 2}^{+\In}z^m\left[\sum_{n= 1}^{m-1}\E_i(\e^{\l S_m};S_1>
S_n,\cdots,S_{n-1}> S_n,S_{n+1}\geq S_n,\cdots,S_m\geq S_n\geq0,X_m= j)\right]\\
  = \ &\delta_{i,j}+\sum_{m= 1}^{+\In}z^m\E_i(\e^{\l
  S_m};S_1\geq0,S_2\geq0,\cdots,S_m\geq0,X_m= j)-z\E_i(\e ^{\l S_1};S_1\geq0,X_1= j)\\
  &-\sum_{m= 2}^{+\In}z^m\left[\sum_{n= 1}^{m}\E_i(\e^{\l S_m};S_1>
S_n,\cdots,S_{n-1}> S_n,S_{n+1}\geq S_n,\cdots,S_m\geq S_n\geq0,X_m= j)\right]\\
  = \ &\delta_{i,j}+\sum_{m= 1}^{+\In}z^m\E_i(\e^{\l
  S_m};S_1\geq0,S_2\geq0,\cdots,S_m\geq0,X_m=j)\\
   &-\sum_{m= 1}^{+\In}z^m\left[\sum_{n= 1}^{m}\E_i(\e^{\l S_m};S_1>
S_n,\cdots,S_{n-1}> S_n,S_{n+1}\geq S_n,\cdots,S_m\geq S_n\geq0,X_m= j)\right]
 \end{split}
\end{equation*}
To prove $F_z(\l)=I$, we have to check that, for any $m\geq1$,
\begin{equation*}
 \begin{split}
&\E_i(\e^{\l
  S_m};S_1\geq0,S_2\geq0,\cdots,S_m\geq0,X_m=j)\\
  & \qquad\qquad\qquad=\sum_{n= 1}^{m}\E_i(\e^{\l S_m};S_1>
S_n,\cdots,S_{n-1}> S_n,S_n\geq 0, S_{n+1}\geq S_n,\cdots,S_m\geq S_n,X_m= j).
 \end{split}
\end{equation*}
Let us thus  consider the random variables $T_m, m\geq 1, $ defined by
$$T_m=\inf\{1\leq n\leq m: S_n=\inf(S_1,\cdots, S_m)\}.$$
We have the following equalities
\begin{equation*}
 \begin{split}
   &\E_i(\e^{\l S_m};S_1\geq0,S_2\geq0,\cdots,S_m\geq0,X_m=j)\\
   =&\sum_{n=1}^{m}\E_i(\e^{\l S_m};S_1\geq0,S_2\geq0,\cdots,S_m\geq0,T_m=n,X_m=j)\\
   =&\sum_{n= 1}^{m}\E_i(\e^{\l S_m};S_1>
S_n,\cdots,S_{n-1}> S_n, S_n\geq 0,S_{n+1}\geq S_n,\cdots,S_m\geq S_n,X_m= j),
 \end{split}
\end{equation*}
which achieves the proof.

The proof of the equality (\ref{eq2.7}) goes along the same lines.
\end{proof}

\begin{remarks}\label{probaexp}
\begin{description}
\item [1)] When $E$ reduces to one point, the sequence $(S_n)_{n\geq 0}$ is a random walk on $\mathbb{R}$ and Proposition
\ref{initfactorization} corresponds to the classical Wiener-Hopf factorization (\cite{Fell}).
\item [2)] There is another way to express the matrices $\daN^* C_z(\l) $ and $\daP B^*_z(\l)$ ;
for $|z|<1$,  one gets
 $$\daN^* C_z(\l)= \Big\{\E_i\Big(z^{T_-^*}\e^{\l S_{T_-^*}}; X_{T_-^*}= j\Big)\Big\}_{i,j} \quad {\rm when}\quad \Re\ \l\geq 0$$
 $$\daP B^*_z(\l)= X^{-1}\Big\{\E_i\Big(z^{\widetilde{T}_+^*}\e^{\l \widetilde{S}_{\widetilde{T}_+^*}}; X_{\widetilde{T}_+^*}= j\Big)\Big\}_{i,j}^{\Large t}X \quad {\rm when}\quad \Re\ \l\leq 0
 \quad ^(\footnote{\rm  where, for any  $N\times N$ complex matrice A,  we denote by  $A^t$ the transposed matrice of  A. }^)$$
 where $X$ is the diagonal matrice
 $\displaystyle X:=   \begin{pmatrix} \nu_1&\quad&(0)\\
                         \quad&\ddots&\quad\\
                         (0)&\quad&\nu_N\end{pmatrix}.$
\end{description}

\end{remarks}
To explain (briefly) how two obtain for instance this ``new'' expression of $\daN^* C_z(\l)$, we  introduce the dual chain $(\widetilde{S}_n,\widetilde{X}_n)$ of $(S_n,X_n)$   whose transition probability  is given  by
$$\widetilde{P}_{(i,x)}(\{j\}\times A)= \frac{\nu_j}{\nu_i}p_{j, i}F(A-x,j,i).$$
We also  consider the    $N\times N$ matrice $\widetilde{C}^-_z$  defined by :

for $|z|<1$, $|\Re\ \l|\leq\a_0$
$$\widetilde{C}^-_z= \left(\sum_{n= 1}^{+\In}z^n\E_i(\e^{\l \widetilde{S}_n}, \widetilde{S}_1\leq0, \widetilde{S}_2\leq0,\cdots,\widetilde{S}_{n-1}\leq0,\widetilde{X}_n= j)\right)_{i,j}.$$
The remark (2) is a straightforward consequence of the\begin{fact}\label{lem2.5} One gets
$\widetilde{C}_z^-= X^{-1} ( B^*_z)^tX.$
\end{fact}
\begin{proof}
We have the equality
\begin{equation*}
\begin{split}
&\E_i(\e^{\l
\widetilde{S}_n},\,\widetilde{S}_1\leq0,\,\cdots,\widetilde{S}_{n-1}\leq0,\,\widetilde{X}_n= j)\\
= &\sum_{k_1,\,k_2,\cdots,k_{n-1}}\int_{\R^n}\frac{\nu_{k_1}}{\nu_i}\,\frac{\nu_{k_2}}{\nu_{k_1}}\cdots\frac{\nu_{j}}{\nu_{k_{n-1}}}\,\e^{\l(\widetilde{y}_1+\cdots\widetilde{y}_n)}1_{[\widetilde{y}_1\leq0]}\,1_{[\widetilde{y}_1+\widetilde{y}_2\leq0]}\cdots1_{[\widetilde{y}_1+\cdots+\widetilde{y}_n\leq0]}\\
&\qquad \times F(k_1,i,\d\widetilde{y}_1)P_{k_1,i}F_(k_2,k_1,\d\widetilde{y}_2)P_{k_2,k_1}\cdots
F(j,k_{n-1},\d\widetilde{y}_n)P_{j,k_{n-1}}.
\end{split}
\end{equation*}
Replacing in this equality $\widetilde{y}_k$ by $y_{n+1-k}$ and $\widetilde{X}_k$ by $X_{n-k}$ for all
$0\leq k\leq n$, we obtain
\begin{equation*}
\begin{split}
&\E_i(\e^{\l
\widetilde{S}_n},\,\widetilde{S}_1\leq0,\,\cdots,\widetilde{S}_{n-1}\leq0,\,\widetilde{X}_n= j)\\
= &\frac{\nu_j}{\nu_i}\sum_{k_1,\cdots,k_{n-1}}\E(\e^{\l
S_n},S_n\leq S_{n-1},\cdots,S_n\leq
S_1,X_0= j,
X_1= k_{n-1},\cdots,X_{n-1}= k_1,X_n= i)\\
= &\frac{\nu_j}{\nu_i}\;\E_j(\e^{\l S_n},S_1\geq S_n, S_2\geq
S_n,\cdots,S_{n-1}\geq S_n,X_n= i).
\end{split}
\end{equation*}
Therefore,
$ \widetilde{C}_z^-(\l)_{i,j}= \frac{\nu_j}{\nu_i} B^*_z(\l)_{j,i}.$
\end{proof}

In the sequel, we will extend this factorization to a larger set of parameters. We will first  prove, by  arguments of elementary type,  that this identity is valid for $\vert z\vert \leq 1$ and ${\rm Re} \ \lambda \in [-\alpha_0, \alpha_0]$. In a second step, we will   extend this identity for ${\rm Re} \ \lambda =0$ and $z$ in a neigbourhood of the unit disc, excepted the point $1$ ; this is much more delicate and it relies on a general argument of algebraic type, due to Presman (\cite{pres2}).

\subsection{General factorization theory of Presman}\label{subsec5.1}

  Let ${\frak R}$ be an arbitrary algebraic ring with unit element $e$ and
$\daI$ be the identity operator in ${\frak R}$. Let the additive operator
$\frak{N}$ be defined on a two-side ideal ${\frak R}'$ of the ring ${\frak R}$, with
\begin{equation}\label{eq2.10}
 ({\frak N} f)({\frak N} g)= {\frak N}[({\frak N} f)g+f({\frak N} g)-fg]
\end{equation}
holding for any $f,g\in{\frak R}'$. It is easy to check that the
operator $\frak P= \daI-{\frak N}$ also satisfies the relation
(\ref{eq2.10}).

\begin{definition}
We say that the element $e-a$ of a ring ${\frak R}$ admits a \textbf{left
canonical factorization} with respect to the operator ${\frak N}$ (l.c.f.
${\frak N}$) if $a\in{\frak R}'$ and if there exist $b,c\in{\frak R}'$ such that
\begin{gather}
e-a= (e-{\frak P} b)(e-{\frak N} c)\label{eq2.12}\\
(e-{\frak P} b)^{-1}= e+{\frak P} c\label{eq2.13}\\
(e-{\frak N} c)^{-1}= e+{\frak N} b\label{eq2.14}.
\end{gather}
 In this case, we say that \textbf{$b$ and $c$ provide a l.c.f.
${\frak N}$}. We call $e-{\frak P} b$ and $e-{\frak N} c$ respectively, \textbf{the
positive and negative components of the l.c.f. ${\frak N}$}.
\end{definition}
The following lemma states the uniqueness of such a factorization once it exists.
\begin{lemma}[\cite{pres2}, lemma 1.1]\label{lem5.2}
 \begin{enumerate}
   \item If $b$ and $c$ provide a l.c.f ${\frak N}$ of the element $e-a$ then
       \begin{enumerate}
        \item\label{itma1} the l.c.f. ${\frak N}$ is unique and is determined by any one of
          the elements ${\frak N} b$, ${\frak P} b$, ${\frak N} c$, ${\frak P} c$;
        \item\label{itema2} for any $\d\in{\frak R}$, the equations
       \begin{alignat}{1}\label{eq2.11}
          x-{\frak P}(xa)= d,\quad y-{\frak N}(ay)= d
        \end{alignat}
        have a unique solution, given by the formulas:
         \begin{gather}
           x= d+\{{\frak P}[da(e+{\frak N} b)]\}(e+{\frak P} c)\label{eq2.15},\\
           y= d+(e+{\frak N} b){\frak N}[(e+{\frak P} c)ad]\label{eq2.16};
         \end{gather}
        \item\label{itema3} for $d= e$, the elements $x= e+{\frak P} c$ and $y= e+{\frak N} b$ are solutions of equation
        (\ref{eq2.11});
       \item\label{itema4} $c_1= c$ (resp. $b_1= b$) is the unique solution of the equation
       $$(e+{\frak P} C_1)(e-a)= e-{\frak N} C_1\quad (\text{ resp. }(e-a)(e+{\frak N} b_1)= e-{\frak P} b).$$
      \end{enumerate}
   \item\label{itemb} If, for $d= e$, equations (\ref{eq2.11}) have solutions $x'$
   and $y'$, then $x'(e-a)y'= e$; moreover, if any two of the three
   elements $x'$, $y'$, $e-a$ are invertible, then $b'= ay'$ and
   $c'= x'a$ provide a l.c.f. ${\frak N}$ of the element $e-a$.
 \end{enumerate}
\end{lemma}

Now,  we assume that $a$ depends analytically on the complex variable $z$ in a neigbourhood of some $z_0$ and describe the regularity of the two components of the l.c.f $ \frak N$ ; namely, we get the following

\begin{lemma}[\cite{pres2}, lemma 1.2]\label{lem5.3}
 Let $a(z)$ be an analytic function in a neighborhood of the point $z_0$, taking values in an ideal ${\frak R}'$ of the Banach algebra ${\frak R}$ and suppose that $b_0$ and $c_0$ provide a l.c.f. ${\frak N}$ of the element $e-a(z_0)$. Then $e-a(z)$ admits l.c.f. ${\frak N}$ in a neighborhood of the point $z_0$, where the elements $b(z)$ and $c(z)$ which provide the l.c.f. ${\frak N}$ of the element $e-a(z)$ are analytic functions of $z$ taking values in ${\frak R}'$.
 \end{lemma}

We achieve this paragraph explaining how one will use this general result in our context.

We will consider  the algebraic ring $V_N[-\alpha_0, \alpha_0]$ of $N\times N$ matrices whose terms are Laplace transforms of Radon measures $\sigma$ on $\R$,  with exponential moment of order $\alpha_0$
The operator ${\frak N}$ will be here the operator $\daN^*$ defined above  and acting on $V_N[-\alpha_0, \alpha_0]$ and $\frak P$ will be equal to $\daP$.

If $\nu$, $\mu$ are two Radon measures on $\R$, we have the following identity :
$$\nu^{*-}\ast\mu^{*-}= (\nu^{*-}\ast\mu+\nu\ast\mu^{*-}-\nu\ast\mu)^{*-}.
\quad ^(
\footnote{where, for any Radon measure  $\gamma$   on $\R$, we have denote by $\ga^{*-}$ its restriction to $\R^{*-}$ defined by $$\ga^{*-}(\d x)= 1_{]-\In,0[}(x)\gamma(\d x).$$}^)$$
Taking into account   this equality, we obtain that  $\daN^*$ and $\daP$ both satisfy the identity (\ref{eq2.10}) for any $f$, $g\in V_N[a,b]$.

For $|z|<1$ and $|\Re\ \l|\leq\a_0$, we will consider  the following $\C$-valued $N\times N$ matrices:
$$ B^*_z(\l):= \left(\sum_{n= 1}^{+\In}z^{n}\int_{-\In}^{+\In}\e^{\l y}\d \P_{i}\{S_1>S_n,S_2>S_n,\cdots,S_{n-1}>S_n,S_n\leq y,X_n= j\}\right)_{i,j},$$
$$C_z(\l):= \left(\sum_{n= 1}^{+\In}z^{n}\int_{-\In}^{+\In}\e^{\l y}\d \P_{i}\{S_1\geq0,S_2\geq0,\cdots,S_{n-1}\geq0,S_n\leq y,X_n= j\}\right)_{i,j}.$$

Recall now that $P(\l)$  belongs to $V_N[-\a_0,\a_0]$ ; furthermore, by Proposition  \ref{initfactorization},  for any complex number $z$ with modulus $ <1$ and any $\lambda \in \C$ such that {\rm Re}$\ \lambda = 0$,  the operator $I-zP(\lambda)$ admits a l.c.f $\daN^*$   on $V_N[ 0, 0]$ provided with $ B^*_z$ and $C_z$.

The above general Presman's result are therefore applicable to $z\mapsto  A_z:= zP(\l)$ with values in $V_N[-\a_0,\a_0]$ for $\vert z\vert <1$ and analytic on the unit open disc  of the complex plane.

 In particular, the elements $B^*_z$ and $C_z$ belong to $V_N[ 0, 0]$. In fact, one may precise this last statement, with the following lemma due to Presman (Lemma 1.3 in \cite{pres2}) :

 \begin{lemma}\label{presfinal}
If $I-A_z$ is an analytic function  of $z$ in a neighbourhood of the point $z_0$, taking values in the ring  $V_N[-\a_0,\a_0]$ and if in this neighbourhood $I-A_z$, as an element of $V_N[0,0]$, admits a l.c.f.  with respect to $\mathcal N^*$ with corresponding elements $B^*_z$ and $C_z$, then $\mathcal P B^*_z$ (resp.  $\mathcal N^*C_z$) is analytic in $z$  in this neigbourhood, with values in
$   \mathcal PV_N]-\In,\a_0]$ (resp. $\mathcal N^*  V_N[-\a_0,+\In[$).
 \end{lemma}

In the sequel, we analyze the factorization of $I-zP(\l)$ in a neigbourhood of the unit disc of the complex plane for some values of $\lambda \in \C$ ; we thus introduce the


\begin{notation}We will   denote by $\mathbb D$ the closed unit ball in the complex number plane : $$\mathbb D:= \{z\in \C: \vert z\vert \leq 1\}.$$
The open unit ball will be denoted $\mathbb D^\circ$.
\end{notation}

\subsection{The factorization of $I-zP(\lambda)$ for $z \in \mathbb D^\circ$ and $\Re \ \l $ closed to $0$}

We first state the  the following

 \begin{theorem}\label{expansiondansdisque}
   There exists $\alpha_1\in ]0, \alpha _0[$ such that for any $z \in \mathbb D^\circ$, one gets
 \begin{enumerate}
 \item
  For  $-\a_1\leq\Re\ \l\leq\a_1$,
  \begin{equation}\label{eq54}
   I-zP(\l)= (I-\daP B^*_z(\l))(I-\daN^* C_z(\l)),
 \end{equation}
 \item
  For   $ \Re\ \l\leq 0$,
  \begin{equation}\label{eq55}
   (I-\daP B^*_z(\l))^{-1}= I+\daP C_z(\l),
  \end{equation}
   \item
 For   $ \Re\ \l\geq 0$,
  \begin{equation}\label{eq56}
   (I-\daN^* C_z(\l))^{-1}= I+\daN^* B^*_z(\l).
  \end{equation}
  \end{enumerate}
  Furthermore, the maps  $z\mapsto \mathcal P  B^*_z(\lambda)$ and $z\mapsto \mathcal N^*C_z (\lambda)$
  are  analytic on $\mathbb D^\circ$  with values  $\mathcal P V_N]-\infty,\a_1]$ (resp. $\mathcal N^* V_N[-\a_1,+\infty[$). \end{theorem}

\begin{proof}
By the  argument developped to establish Proposition \ref{initfactorization}, one checks easily that  (\ref{eq2.6}) (resp. (\ref{eq2.7}))  is valid for $\vert z\vert <1$ and   ${\rm Re}\ \l \leq 0$ ( resp. ${\rm Re}\ \l \geq 0$). So (\ref{eq55}) and (\ref{eq56}) are valid.

The existence of the factorization in $V_N[0,0]$ for any $z\in \mathbb D^\circ$ is given by Proposition \ref{initfactorization}. The analyticity of the different components $B^*_z(\l)$ and $C_z(\l)$ on $ \mathbb D^\circ$  for $\Re \ \l=0$  is a consequence of Lemma  \ref{lem5.3} ; we may also apply Lemma \ref{presfinal}  and  conclude that
$$
\daP B^*_z\in \daP V_N]-\infty ,\a_0]
\quad \mbox{\rm and} \quad
\daN^* C_z\in \daN^* V_N[-\a_0,+\infty[.$$

Now, for any $z \in \mathbb D^\circ$, the maps $\lambda\mapsto I-zP(\l)$ and $\lambda \mapsto (I-\daP B^*_{z}(\l))( I-\daN^* C_{z}(\l))$ are analytic on the strip $\{\vert {\rm Re}\ \l \vert \leq \alpha_1\}$ for any $\alpha_1\in ]0, \alpha_0[$ and they coincide on the line ${\rm Re}\ \l=0$ ; they thus coincide as analytic functions on the strip $\{\vert {\rm Re}\ \l \vert \leq \alpha_1\}$.  So (\ref{eq54}) holds for $-\a_1\leq\Re\ \l\leq\a_1$ ; the analyticity of the maps  $z\mapsto \mathcal P  B^*_z(\lambda)$ and $z\mapsto \mathcal N^*C_z (\lambda)$
  are  analytic on $\mathbb D^\circ$  with values  $\mathcal P V_N]-\infty,\a_1]$ (resp. $\mathcal N^* V_N[-\a_1,+\infty[$ ) is a direct consequence of  Lemma \ref{presfinal}.

  For the analyticity of these two maps when $\Re \ \l \in ]-\alpha_1, \a_1]$, one may also  use the explicit form of the functions $B^*_z$ and $C_z$ and argue as follows :\\
-  for  $\Re\  \l = 0$, it is a consequence of Lemma \ref{lem5.3} as we said a few lines above ;\\
 - when $\Re\ \l > 0$, it is a direct consequence of the identity
  $$ \daN^* C_z(\l)= \Big\{\E_i\Big(z^{T_-^*}\e^{\l S_{T_-^*}}; X_{T_-^*}= j\Big)\Big\}_{i,j};$$
\\
  - when $\Re\ \l \in [-\alpha_1, 0[$, we use  (\ref{eq54}) and (\ref{eq55}) to write
  $$
  \daN^* C_z(\l)= I- (I+\daP C_z(\l))(I-zP(\l))
  $$
with $\displaystyle \daP C_z(\l)= \Big\{\sum_{n\geq 0}z^{n}\E_i\Big(\e^{\l S_{n}}; T_{-}^*>n, X_{n}= j\Big)\Big\}_{i,j}.$ The two factors on the right hand side of this last equality  are clearly analytic in $z \in \mathbb D^\circ$ and the result follows.
The same argument holds for   $z\mapsto \mathcal P  B^*_z(\lambda)$.
\end{proof}

\subsection{Expansion of the factorization outside the unit disc and far from $z=1$}

We  study here the extension of the preceding factorization when
$\Re\ \l\  = 0$ and $z$ lives in a neighbourhood of $ \mathbb D \setminus\{1\}$. We have the
\begin{theorem}\label{expansiondehors}
There exists  a neighbourhood $U$ of $ \mathbb D \setminus\{1\}$    such that, for $\Re\ \l= 0$, the two maps $z \mapsto B^*_z(\l)$ and $z \mapsto C_z(\l)$  may be   continuously expanded on  $U$  in such a way
\begin{enumerate}
 \item for any $z\in U_{}$, the formulas (\ref{eq54}), (\ref{eq55}) and (\ref{eq56}) hold.
 \item the maps  $z\mapsto \daP B^*_z$ and $z\mapsto\daN^* C_z$ are analytic  on $U  $, with values in $V_N]-\In,\a_0]$ and  $V_N[-\a_0,+\In[$ respectively.
\end{enumerate}
\end{theorem}

\begin{proof}
We fix $\lambda $ s.t. $\Re\  \l = 0$,  $z_0 \in \mathbb C$ with $\vert z_0\vert =1,  z_0\neq 1$ and choose a sequence $(z_n)_{n\geq 1}$ of complex numbers  in $\mathbb D^\circ$ which converges to $z_0$.

By 2) of Remarks \ref{probaexp}, the two limits
$\displaystyle  B^+_{z_0}(\l):=\lim_{n \to +\infty}\daP B^*_{z_n}(\l)$ and
$\displaystyle C^{*-}_{z_0}(\l):=\lim_{n \to +\infty}\daN^* C_{z_n}(\l)$ do exist ; furthermore,   (\ref{eq54})  holds at any point $z_n$ and letting $n \to +\infty$ one gets
$$
I-z_0P(\l)= (I-B^+_{z_0}(\l))(I-C^{*-}_{z_0}(\l)).
$$
Since $z_0\neq 1$, the matrice $I-z_0P(\l)$ is invertible, so is $I-B^+_{z_0}(\l)$  ; by (\ref{eq55}), the limit
 $\displaystyle
\lim_{n \to +\infty} PC_{z_n}(\l)
 $
 does also exists (and is equal to $C_{z_0}^+(\l):= -I+ (I-B^+_{z_0}(\l))^{-1}$).

Consequently, $\displaystyle
C_{z_0}(\l):=  \lim_{n \to +\infty}  C_{z_n}(\l)  =
\lim_{n \to +\infty} \daN^*C_{z_n}(\l)
 +
\lim_{n \to +\infty} \daP C_{z_n}(\l) =C^{*-}_{z_0}(\l)+C^{+}_{z_0}(\l) $ does exist and one gets  $C^{*-}_{z_0}= \daN^*C_{z_0}(\l)$ and
 $C^{+}_{z_0}(\l)=\daP C_{z_0}(\l)$.

By the same argument, on shows that
$\displaystyle
B^*_{z_0}(\l):=\lim_{n \to +\infty}  B^*_{z_n}(\l)$ does exist  and (\ref{eq56}) holds at $z_0$.

Finally $B^*_{z_0}(\l)$ and $C_{z_0}(\l)$ provide a l.c.f $\daN^*$ of $I-z_0P(\l)$ ; since $z\mapsto
I-z_0P(\l)$ is analytic in a neigbourhood of $z_0$, so are the maps $z\mapsto {\mathcal P}B^*_{z }(\l)$ and
 $z\mapsto {\mathcal N}^*C_{z }(\l)$ by Lemma \ref{lem5.3}, with values in ${\mathcal P}V_N]-\In,\a_0]$ and  ${\mathcal N}^*V_N[-\a_0,+\In[$ respectively, by Lemma \ref{presfinal}.
\end{proof}
In the sequel we will specify the neigbourhood $U$ as follows ;   recall that
$$D_{\rho,\theta}:=  \{z;z\neq1,|\arg(z-1)|>\theta>0,|z|<\rho\}$$
 and
 $$K(\de_1,\de_2):=\{z:q+\de_1<|z|<1+\de_2,\Re\ z>0,|\Im
 z|<\de_1\}.$$
We have the

\begin{corollary}\label{Rhodelta}
 There exist $\rho>1$ and $\theta \in ]0, \pi/2[$ such that

 $\bullet$  the formulas (\ref{eq54}), (\ref{eq55}) and (\ref{eq56}) hold for $\Re \ \l \ =0$ and $z \in D_{\rho,\theta}\setminus K(\de_1,\de_2)$,

  $\bullet$ for $\vert{\rm Re}(\l)\vert \leq \a_0$, the map $z\mapsto \daP B^*_z$ (resp. $z\mapsto\daN^* C_z$ ) is analytic  on $D_{\rho,\theta}\setminus K(\de_1,\de_2)$; furthermore, $I-\daP B^*_z$ (resp. $I-\daN^*C_z$) is invertible (and their inverses are also analytic)  on this domain.

\end{corollary}

\section{On the local behavior of the factors of the Laplace transform of the minimum}\label{sec4}

We know, by Lemma \ref{lem4.1} that the Laplace transform of the minimum $m_n$ may be decomposed as follows :
for $\l>0$ and $|z|<1$,
$$
 \sum_{n= 0}^{+\In}z^n\E_i(\e^{\l m_n};X_n= j)= \{[I+\daN^* B^*_z(\l)][I+\daP C_z(0)]\}_{i,j}.$$
 In this section, we will study each the behavior of these two factors near $z=1$. More precisely, we will first consider the case when$\vert z\vert \leq 1$  and after investigate the case when $\vert z\vert >1$.
\subsection{Preliminaries}
As mentioned in the previous section, the matrices $I+\daN^* B^*_z(\l)$ and $I+\daP C_z(0)$ could be seen as the inverse of two factors for the matrix $I-zP(\l)$, we will first study the regularities of these quantities for $z\in\overline{K}(\delta_1,\delta_2)$. In the following , the constants $\delta$ and $\varepsilon$ are choosen small enough in such a way that, for $z \in \bar K(\delta, 0)$, one gets $[\lambda_-(z)-\varepsilon, \lambda_+(z)+\varepsilon]\subset ]-\alpha_0, \alpha_0[$.\\
We have the

 \begin{proposition}\label{Kdelta}
 There exist  $\delta_1>0$, for $z\in\overline{K}(\de_1,0)$, and any $\varepsilon>0$ such that Theorem \ref{thm2.4.2} is satisfied, one gets
 \begin{enumerate}
  \item for $\Re\ \l <\l_+(z,\varepsilon)$ with $\l\neq\l_+(z)$,
   \begin{equation}\label{eq01}
    (I-\daP B^*_z(\l))^{-1}=I+\daP C_z(\l)=I-\frac{\left[I-\daN^* C_z(\l_+(z))\right]\Pi_+(z)}{(\l_+(z)-\l)\b_+(z)}+\int_{0}^{+\In}\e^{\l
 x} k_+(z,\d x)
   \end{equation}
   \item for $\Re\ \l >\l_-(z,\varepsilon)$ with $\l\neq\l_-(z)$,
   \begin{equation}\label{eq02}
   (I-\daN^* C_z(\l))^{-1}=I+\daN^* B^*_z(\l)=I-\frac{\Pi_-(z)\left[I-\daP B^*_z(\l_-(z))\right]}{(\l_-(z)-\l)\b_-(z)}+\int_{-\In}^{0}\e^{\l
 x}k_-(z,\d x)
 \end{equation}
 \end{enumerate}
where $k_+(z,\cdot)$ (resp. $k_-(z,\cdot)$) is a measure on $[0,+\In[$ (resp. $]-\In, 0]$)  taking values in the vector space $M_{N\times N}(\C)$ of $N\times N$ complex matrices, such that for $z\in\overline{K}(\de,0)$, one gets
 \begin{equation}\label{eq03}
  \|k_+(z,x)\|\leq C\e^{-(\l_+(z)+\varepsilon )x}\quad for \quad
  x>0,
\end{equation}
\begin{equation}\label{eq04}
 \|k_-(z,x)\|\leq C\e^{-(\l_-(z)-\varepsilon)x} \quad for \quad
  x<0,
\end{equation}
where $k_+(z,x)=k_+(z,]x,+\In[)$ for $x>0$ and $k_-(z,x)=k_-(z,]-\In,x[)$  for $x<0$.

\noindent Furthermore,   the following limits exist :
 \begin{equation}\label{limit}
 \lim_{|z|\uparrow1}\frac{(I-\daN^* C_z(\l_+(z)))\Pi_+(z)}{\b_+(z)}=A_+ \qquad \mbox{\rm and} \quad
 \lim_{|z|\uparrow1}\frac{\Pi_-(z)(I-\daP B^*_z(\l_-(z)))}{\b_-(z)}=A_-,
  \end{equation}
 where $A_+$(resp. $A_-$) is a $N\times N$ matrix  with non positive (resp. non negative) coefficients.
 \end{proposition}

 \begin{proof}
Since the probabilistic expression of $\daN^*C_z$ is quite simple,
we first prove that (\ref{eq01}) and (\ref{eq03})  hold  when $z\in K(\delta,0)$ for any $0<\delta<\a_0$ ; then, we will establish the existence of $A_+$ in (\ref{limit})  when $\delta$ is quite small (namely $\delta \leq \delta_1$),  which will allows us to prove that (\ref{eq01}) and (\ref{eq03}) holds in fact for $z\in \overline K(\delta_1,0)$ and $\Re\ \ \l<z_+(z,\varepsilon)$, $\l\neq\l_+(z)$ .

We first prove that equality (\ref{eq01})  holds for $z\in K(\delta,0)$, $0<\delta<\a_0$ ; the same argument works to establish   (\ref{eq02}).\\
 According to Theorem \ref{expansiondansdisque} and the definition of $\daP
C_z(\l)$, for $q\leq \mid z\mid<1$ and $\Re\ \l\leq0$, one gets
\begin{equation}\label{eq47}
 \begin{split}
  (I-\daP B^*_z(\l))^{-1}=I+\daP C_z(\l)&=I+\sum_{n=1}^{+\In}z^n\,\E_i(\e^{\l S_n},T^*_->n,X_n=j)\\
                                      &=I+\int_0^{+\In}\e^{\l y}\d b_+(z,y).
 \end{split}
\end{equation}
By (\ref{eq54}) and the inversion formula of Laplace, for $\l_-(z)<-\delta<0$, one may write for $x>0$,
$$b_+(z,x)-b_+(z,-\infty)=-\frac{1}{2\pi i}\int_{\Re\ \l=-\delta}\e^{-\l x}\frac{(I-\daN^* C_z(\l))(I-zP(\l))^{-1}}{\l}\,\d\l.$$
Now we transfer the contour of integration to the straight line
$\Re\ \l=\l_+(z)+\varepsilon$; using Cauchy's formula on the convex open set $\Omega=\{-\delta<\Re\ \l<\Re\ \l_+(z)+\varepsilon,|\Im\l|<\b\}$ and the fact that $\displaystyle \lambda \mapsto{1\over \lambda} (I-\daN^* C_z(\l))(I-zP(\l))^{-1}$ is analytic in $\Omega\setminus\{0,\l_+(z)\}$, we get for $y>0$,
\begin{equation}\label{eq48}
\begin{split}
 b_+(z,y)-&b_+(z,-\In)=-(I-\daN^* C_z(0))(I-zP(0))^{-1}+\frac{\e^{-\l_+(z)y}[I-\daN^*
 C_z(\l_+(z))]\Pi_+(z)}{\b_+(z)\l_+(z)}\\
          &-\frac{\e^{-(\l_+(z)+\varepsilon)y}}{2\pi
          i}\int_{\Re\ \l=0}\e^{-\l y}\frac{[I-\daN^*
          C_z(\l+\l_+(z)+\varepsilon)][I-zP(\l+\l_+(z)+\varepsilon)]^{-1}}{\l+\l_+(z)+\varepsilon}\,\d\l.
\end{split}
\end{equation}
As in the proof of Theorem  \ref{thm2.4.2}, we set
$\l_+(z, \var):= \Re\ \l_{+}(z)+\varepsilon$     and, for $x\geq 0$
 \begin{equation}\label{eq49}
k_+(z,x)=-\frac{\e^{- \l_+(z,\varepsilon) x}}{2\pi
          }\int_{\R}\e^{-i\theta x}\frac{[I-\daN^*
          C_z(\l_+(z, \var)+i\theta][I-zP(\l_+(z, \var)+i\theta)]^{-1}}{\l_+(z, \var)+i\theta}\,\d\theta.
\end{equation}
Consequently, for $z\in K(\delta,0)$, $\Re\ \l<\l_+(z,\varepsilon)$ and $\l\neq\l_+(z)$, one gets
\begin{equation}\label{eqsansbar}
 (I-\daP B^*_z(\l))^{-1}=I-\frac{\left[I-\daN^* C_z(\l_+(z))\right]\Pi_+(z)}{(\l_+(z)-\l)\b_+(z)}+\int_{0}^{+\In}\e^{\l
 x} k_+(z,\d x).
\end{equation}
Inequality (\ref{eq03}) is thus a direct consequence of  the following result, which is the analogous in the present  context of Property \ref{TDFdesresolvantes}

\begin{property}\label{TDFdesfacteursderesolvantes}
We fix $\varepsilon>0$ and $\delta_1$ small enough in such a way that Theorem \ref{thm2.4.2} is satisfied.
We set
\begin{description}
  \item[$\bullet$]$\displaystyle W'_+(z, \varepsilon, x)=\int_\R\e^{-i\theta
 x}\frac{[I-\daN^*
          C_z(\l_+(z, \var)+i\theta)][I-zP(\l_+(z, \var)+i\theta)]^{-1}}{\l_+(z, \var)+i\theta}\;\d
 \theta,$\\
             \text{\hspace{11cm} for} $x\geq 0$;
 \item[$\bullet$] $\displaystyle W'_-(z, \varepsilon, x )=\int_\R\e^{-i\theta
 x}\frac{[I-zP(\l_-(z,\varepsilon )+i\theta)]^{-1}[I-\daP
          B^*_z(\l_+(z, \var)+i\theta)]}{\l_-(z,\varepsilon)+i\theta}\;\d
 \theta,$\\
             \text{\hspace{11cm} for} $x<0$;
\end{description}
 Then, there  exists a constant $C'=C'(\varepsilon)>0$ such that for $x\geq0$ (resp. $x<0$), one gets
\begin{equation}\label{4}
\forall z\in\overline{K}(\de_1,0),\quad\|W'_+(z,x,\varepsilon)\|\leq C \quad(\text{resp. } \|W'_-(z,x,\varepsilon)\|\leq C).
\end{equation}
\end{property}

Let us now establish (\ref{limit}). Since for any $|z|\leq1$, $$\daN^*C_z(\l_+(z))=\Big\{\E_i\Big(z^{T^*_-}\e^{\l_+(z) S_{T^*_-}};X_{T^*_-}=j\Big)\Big\}_{i,j}$$
and $\lim_{z\rightarrow1}\l_+(z)=0$, we obtain that for any $z\in\overline{K}(0,\delta)$,
 $$\|\daN^* C_z(\l_+(z))\|\leq\Big\|\Big\{\E_i\Big(z^{T^*_-}\e^{\Re\ \l_+(z) S_{T^*_-}};X_{T^*_-}=j\Big)\Big\}\Big\|<+\In.$$
Moreover, by the second assertion of Theorem \ref{thm2.4.2}, we may choose $\delta_1\leq\delta$ and $0<\varepsilon_i<\a_0$, $i=1, 2$, such that $\|(I-zP(\l))^{-1}\|<+\In$ for all $z\in \overline{K}(\delta_1,0)$ and $\varepsilon_1\leq\Re\ \l\leq\varepsilon_2$.

Therefore, for any $\varepsilon_1\leq\Re\ \l\leq\varepsilon_2$ and $|z|<1$, one gets
$$(I-\daP B^*_z(\l))^{-1}=(I-\daN^* C_z(\l))(I-zP(\l))^{-1}$$
 and the limits as $\vert z\vert \to 1$ of the two factors on the right hand side do exist ; this implies that $(I-\daP B^*_z(\l))^{-1}$ exists for $z \in  \overline{K}(\delta_1,0)$ and $\varepsilon_1\leq\Re\ \l\leq\varepsilon_2$, with the identity
\begin{equation}\label{limitofdapb}
(I-\daP B^*_z(\l))^{-1}=(I-\daN^* C_z(\l))(I-zP(\l))^{-1}.
\end{equation}

In particular, letting  $|z|\rightarrow1$ in (\ref{eqsansbar}), we obtain
\begin{equation}\label{daA}
 \lim_{|z|\uparrow1}\frac{[I-\daN^* C_z(\l_+(z))]\Pi_+(z)}{\b_+(z)}\quad \mbox{\rm exists} \quad (= A_+).
\end{equation}
It remains to prove that (\ref{eq01}) holds for $|z|=1$, $\Re\ \l<\l_+(z,\varepsilon)$ and $\l\neq\l_+(z)$. Taking  into account
(\ref{limitofdapb}) and (\ref{daA}), we can comfirm that for any $\varepsilon_1\leq\Re\ \l\leq \varepsilon_2$, as $|z|\rightarrow1$, the limits for the members in the equality (\ref{eqsansbar}) exist and (\ref{eq01}) hold for $|z|=1$. Since the different members in (\ref{eq01}) exist as Laplace transforms (of   certain measures)  for  $\Re\ \l<\l_+(z,\varepsilon)$ and $\l\neq\l_+(z)$ and any fixed $z\in\overline{K}(0,\delta_1)$, this equality (\ref{eq01})  holds in fact for  such values of $z$ and $\lambda.$

The equalities (\ref{eq02}),  (\ref{eq04}) and the existence of $A_-$  may  be proved with the same method. \end{proof}

It remains to give the main lines of the proof of Property  \ref{TDFdesfacteursderesolvantes}.
\begin{proof} [Proof of Property \ref{TDFdesfacteursderesolvantes}]  We just give the main steps of the proof for $W'_+(z, \var, x)$, which is quite similar to the one of Property \ref{TDFdesresolvantes}  ;  we also  set $\blacktriangle:=\l_+(z,\varepsilon)+i\theta$, and  decompose $W'_+(z, \var, x)$ as  $W'_+(z, \var, x)=W'_{+1}(z, \var, x)+W'_{+2}(z, \var, x)+W'_{+3}(z, \var, x)$ with
 $$W'_{+1}(z,\varepsilon, x):=\int_\R {e^{-i\theta x}\over \blacktriangle} [I-\daN^*
          C_z(\blacktriangle)] \;\d\theta,$$
$$W'_{+2}(z,\varepsilon, x):=\int_\R {e^{-i\theta x}\over \blacktriangle}   [I-\daN^*
          C_z(\blacktriangle)][zP(\blacktriangle)+\cdots+z^{n_1}P^{n_1}(\blacktriangle)][I-\mathfrak{L}(B)(z,\blacktriangle)]^{-1}  \;\d\theta,$$
\begin{eqnarray*}W'_{+3}(z,\varepsilon, x)&:=&\int_\R {e^{-i\theta x}\over \blacktriangle}
 [I-\daN^*
          C_z(\blacktriangle)]z^{n_1}[zP(\blacktriangle)+\cdots+z^{n_1}P^{n_1}(\blacktriangle)]\\
          &\ & \qquad \qquad \qquad \qquad \times
          [I-z^{n_1}P^{n_1}(\blacktriangle)]^{-1}\mathfrak{L}(\Phi_{n_1,\kappa})(\blacktriangle)[I-\mathfrak{L}(B)(z,\blacktriangle)]^{-1}
          \;\d
                      \theta.\end{eqnarray*}
To check that $W'_{+1}(z,\varepsilon, x)$ is bounded uniformly in $z \in \overline K(\delta_1, \delta_2)$  and $x\geq 0$, one first uses   Lemma  \ref{intcomplexe'} to get
\begin{eqnarray*}
 W'_{+1}(z,\varepsilon, x) &=& \int_\R {e^{-ix\theta} I\over \l_+(z, \var)+i\th}\d \th-\Bigl(\mathbb E_i\Bigl[z^{T^*_-}e^{\l_+(z, \var)S_{T^*_-}} \int_\R{e^{i\theta(S_{T^*_-}-x)}\over \blacktriangle}\d \theta ; X_{T^*_-}=j  \Bigr]\Bigr)_{i, j}=0.
\end{eqnarray*}

To control  $W'_{+2}(z,\varepsilon, x)$, one uses the fact that  the function
 $z\mapsto[I-\daN^*
          C_z(\blacktriangle)][zP(\blacktriangle)+\cdots+z^{n_1}P^{n_1}(\blacktriangle)][I-\mathfrak{L}(B)(z,\blacktriangle)]^{-1}$ is
 the Laplace transform  at point
 $\blacktriangle$ of the measure $\mu'(z,\d y)=N_z(\d y)\bullet [zM(\d y)+\cdots+z^{n_1}M^{n_1}(\d y)]\bullet\widetilde{B}(z,\d y)$, where $\displaystyle N_z(\d y):= \Big\{\mathbb \delta_{i, j}(\d y)-\E_i(z^{T^*_-},  S_{T^*_- }\in \d y, X_{T^*_-}=j)\Big\}_{i, j}$
 and one may conclude as in the proof of Property \ref{TDFdesresolvantes}.

 The control of  $W'_{+3}(z,\varepsilon, x)$ is like the one of $W_{+3}(z,\varepsilon, x)$ in Property \ref{TDFdesresolvantes}. The proof for
 $W'_{-}(z,\varepsilon, x)$ and $x<0$ goes along the same lines.
\end{proof}

 In the following Proposition, we precise the type of regularity of $(I-\daP B^*_z(\l))^{-1}$ and $(I-\daN^* C_z(\l))^{-1}$ on the domain $K^*(\delta_1, \delta_2)$ for small enough $\delta_1, \delta_2 >0$ (by Corollary \ref{Rhodelta}, we already know that they are analytic on $D_{\rho,\theta}\cap \Bigl(K(\de_1,\de_2)\Bigr)^c$ for some suitable $\rho>1$ and $\theta >0$).

  We set
$$F_\pm(z,\l)=I+\frac{\l_\pm(z)-a_\pm}{\l-\l_\pm(z)}\,\Pi_\pm(z),$$
where $a_+=\a_0+1$ and $a_-=-\a_0-1$. Recall that for $z\in K(\delta_1,\delta_2)$, the matrices $\Pi_\pm(z):=\Pi(\l_\pm(z))$ are rank $1$ and given by
$$\Pi_\pm(\l)=\Big(e_i(\l_\pm(z))\nu_j(\l_\pm(z))\Big)_{i,j\in E},$$
with ${}^t\!\nu(\l_\pm(z))\e(\l_\pm(z))=1$.

Note that $F_\pm (z, \l)$ are analytic with respect to $z \in K^*(\de_1, \de_2)$, excepted at ${1\over k(\l)}$ (so that $\l  \neq \l_\pm(z)$).

On the other hand, one gets
$$F^{-1}_+(z,\l)=I-\frac{\l_+(z)-a_+}{\l-a_+}\,\Pi_+(z)$$
(and similarly $\displaystyle F^{-1}_-(z,\l)=I-\frac{\l_-(z)-a_-}{\l-a_-}\,\Pi_-(z)$)$^($\footnote{Remark that for any column vector $a$ and row vector  $b$ , setting $ba=\beta \in \C$, then
$$\det(I-ab)=1-\b\text{ and }(I-ab)^{-1}=I+(1-\b)^{-1}ab.$$
One applies  these formulae to $a=-\frac{\l_+(z)-a_+}{\l-\l_+(z)}\left(
e_1(\l_+(z)),\cdots,e_N(\l_+(z))\right)^t$ and $b=(\nu_1(\l_+(z)),\cdots,\nu_N(\l_+(z)))$
to obtain the announced  expression of $F_+^{-1}(z, \l)$; and the expression of $F_-^{-1}(z, \l)$ can be obtained analogously.}$^)$. Let us emphasize that $F_\pm^{-1}(z, \l)$ are analytic on $K^*(\de_1, \de_2)$ (even at point ${1\over k(\l)}$ ).

We now set $ B(z,\l)=F_+(z,\l)(I-zP(\l))F_-(z,\l) ; $ by the above, the matrice $ B(z,\l)$ is invertible, we denote by $ B^{-1}(z,\l) $ its inverse; we also set $B_+(z,\l)=F_+(z,\l)(I-\daP B^*_z(\l))$ and $B_-(z,\l)=(I-\daN^* C_z(\l))F_-(z,\l)$.

For $z\in\overline{K}(\delta_1,0)$, according to the relation (\ref{eq54}), we have
 \begin{equation}\label{eq57}
  B(z,\l)=B_+(z,\l)B_-(z,\l),\quad |\Re\ \l|\leq\a_0,
 \end{equation}
 \begin{equation}\label{eq58}
  B^{-1}(z,\l)=B_-^{-1}(z,\l)B_+^{-1}(z,\l),\quad \l\in
  S_z(\varepsilon).
 \end{equation}

 The regularity of $ B(z,\l), B_\pm(z,\l)$ and $ B^{-1}(z,\l), B^{-1}_\pm (z,\l)$  is described in the following

\begin{proposition}\label{propo2}
 For  $\de_1,\de_2,\varepsilon>0$ small enough and  $z\in
 K^*(\de_1,\de_2)$, one gets
$$
 B(z,\l) \in V_N[-\a_0,+\a_0] \quad \mbox{\rm and} \quad
 B^{-1}(z,\l)= V_N[ \l_-(z,\varepsilon),\l_+(z,\varepsilon)].$$

Furthermore, the maps
 \begin{description}
 \item \indent $\bullet$ $z\longmapsto B (z,\l)$,   $z\longmapsto B_- (z,\l)$, $z\longmapsto B_+ (z,\l)$
 \item \indent $\bullet$
 $z\longmapsto B^{-1} (z,\l)$, $z\longmapsto B^{-1}_-(z,\l)$, $z\longmapsto B^{-1}_+ (z,\l)$
 \item \indent $\bullet$  $z\longmapsto \daP B^*_z(\l)$,
 $z\longmapsto \daN^*
C_z(\l)$,
\end{description} admit an analytic expansion on  $K^* (\de_1,\de_2)$   and   with respect to the variable $t=\sqrt{1-z}$ for $z \in
K^*(\de_1,\de_2)$.

 Furthermore, the maps    $z\longmapsto(I-\daP B^*_z(\l))^{-1}$ and $z\longmapsto(I-\daN^*
C_z(\l))^{-1}$
are analytic on  $K^* (\de_1,\de_2)$ excepted at point  $ {1\over k(\lambda)}$ ; in particular, they are analytic on $D_{\rho, \theta}$.\end{proposition}

\begin{proof}
 We first assume that $\delta_1$ is choosen in such a way that the conclusions of Proposition \ref{Kdelta} are valid.
 Since $B_+(z,\l)\in V[-\a_0,\a_0]$, by the formula (\ref{eq01}) in Proposition \ref{Kdelta}, we find
 $$B_+^{-1}(z,\l)=\Bigl(I+\int_{0}^{+\In}\e^{\l x}  k_+(z,dx)\Bigr)F_+^{-1}(z,\l)+\frac{(I-\daN^* C_z(\l_+(z)))\Pi_+(z)}{(\l-a_+)\b_+(z)}.$$
 The equality (\ref{limit}) thus implies that $B_+^{-1}(z,\l)$ is bounded for $z\in\overline{K}(\delta_1,0)$ and $\lambda \in S_z(\varepsilon)$. The same holds for $B_-^{-1}(z,\l)$.\\

 The relations (\ref{eq57}) and (\ref{eq58}) show that $B^{\pm1}(z,\l)$ admit a canonical factorization for all $z$ on the unit circle such that $|\Im z|<\de_1$. Since these functions are regular with respect to the variable $t=\sqrt{1-z}$ for $ z\in K^*(\delta_1,\delta_2)$, we may by  Lemma \ref{lem5.3} adapt the choice of $\delta_1$ and $\delta_2$ in such a way  that the components of factorizations (\ref{eq57}) and (\ref{eq58}), regarded as functions of $t$, admit an   analytic expansion with respect to the variable $t$. By the identity
 \begin{equation}\label{eq59'}
\daP B^*_z(\l)=I-F_+^{-1}(z,\l)B_+ (z,\l),
 \end{equation}
 one obtains the expected regularity of the functions
 $z\mapsto\daP B^*_z(\l)$.

 At last, for $z\neq 1/k(\lambda)$, one gets  by the previous equality
 \begin{equation}\label{eq59}
 (I-\daP B^*_z(\l))^{-1}=B_+^{-1}(z,\l)F_+(z,\l),
 \end{equation}
 with $F_+(z,\l)$ well  defined and analytic in $z$ since $\lambda \neq \lambda_{\pm}(z)$
 and one concludes.

 The same holds similarly for $\l\mapsto\daN^* C_z(\l)$ and $\l\mapsto(I-\daN^* C_z(\l))^{-1}$.
\end{proof}

\subsection{ On the regularity of the factors  $I+\daN^* B^*_z(\l) $ and $I+\daP C_z(\l) $  on $D_{\rho,\theta}$ for $\lambda \in \R^*$}

 In this section we fix $\rho>1$ and $\th\in ]0, \pi/2[$ such that the conclusions of Corollary \ref{Rhodelta} hold. We   prove the
 \begin{theorem}\label{NOsingularite}
 \begin{enumerate}
 \item
  For $\l>0$ (resp. $\l<0$) closed to $0$, the function $I+\daN^* B^*_z(\l)$ (resp. $I+\daP C_z(\l)$)   admits an analytic expansion on
  $D_{\rho,\theta}$.
  \item We have
\begin{equation}\label{eq29}
 \lim_{\l\rightarrow0^+}\l(I+\daN^* B^*_1(\l))=A_-,
\end{equation}
\begin{equation}\label{eq30}
 \lim_{\l\rightarrow0^-}\l(I+\daP C_1(\l))=A_+,
\end{equation}
with
\begin{equation}\label{eq28}
 -\frac{k''(0)}{2}A_-A_+=\Pi(0).
\end{equation}
\end{enumerate}
 \end{theorem}

\begin{proof}
\begin{enumerate}
\item
 First case : when $z\in D_{\rho,\theta}\setminus K(\delta_1,\delta_2)$ and $\l\in\R^*$, this is a direct consequence of Corollary  \ref{Rhodelta}.

     Second case : when $z\in K(\delta_1, 0)$, by the first assertion of Theorem \ref{expansiondansdisque}, we have
     \begin{equation}\label{eq61}
      (I-\daN^* C_z(\l))^{-1}=I+\daN^* B_z(\l),\quad \Re\ \l\geq0,
     \end{equation}
     \begin{equation}\label{eq62}
      (I-\daP B^*_z(\l))^{-1}=I+\daP C_z(\l),\quad \Re\ \l\leq0.
     \end{equation}
    Now, by Proposition \ref{propo2}, the quantities of left hand-side of the above formulae are proved to be analytic with respect to $z\in K^*(\delta_1,\delta_2)$ for some $\delta_2>0$ small enough and  for $z \neq {1\over k(\l)}$. Recall that $z \neq {1\over k(\l)}\Leftrightarrow  \lambda \neq \lambda_\pm(z)$. We hence obtain the expected  result, using the fact that $D_{\rho,\theta}\subset(D_{\rho,\theta}\setminus K(\delta_1,\delta_2))\cup K^*(\delta_1,\delta_2)$.

\item The equalities (\ref{eq29}) and (\ref{eq30}) are direct consequences of Proposition \ref{Kdelta}. Indeed, according to this Proposition,  one gets
\begin{equation}\label{eq63}
 \lim_{z\rightarrow1}\l_+(z)(I-\daP B^*_z(0))^{-1}=-A_+,
\end{equation}
\begin{equation}\label{eq64}
 \lim_{z\rightarrow1}\l_-(z)(I-\daN^* C_z(0))^{-1}=-A_-.
\end{equation}
On the other hand, for $q<z<1$, one gets
$$(1-z)(I-zP(0))^{-1}=[\sqrt{1-z}\,(I-\daN^* C_z(0))^{-1}][\sqrt{1-z}\,(I-\daP B^*_z(0))^{-1}],$$
with $(I-zP(0))^{-1}=\frac{z\Pi(0)}{1-z}+\sum_{n=0}^{+\In}z^n R^n(0) ;$
so
$$
\lim_{z\rightarrow1}[\sqrt{1-z}\,(I-\daN^* C_z(0))^{-1}][\sqrt{1-z}\,(I-\daP B^*_z(0))^{-1}]=\Pi(0).$$
Since $\lim_{z\rightarrow1}\frac{\sqrt{1-z}}{\l_-(z)}=-\sqrt{\frac{k''(0)}{2}}$ and $\lim_{z\rightarrow1}\frac{\sqrt{1-z}}{\l_+(z)}=\sqrt{\frac{k''(0)}{2}}$ (see (\ref{eq2.4.9})), we hence obtain
$$-\frac{k''(0)}{2}A_-A_+=\Pi(0)=\begin{pmatrix} \nu_1&\nu_2&\cdots&\nu_N\\
                         \nu_1&\nu_2&\cdots&\nu_N\\
                         \qquad&\cdots\\
                         \nu_1&\nu_2&\cdots&\nu_N\end{pmatrix},$$
which yields to the result.

\end{enumerate}

\end{proof}

\subsection{ On the regularity of the factors  $I+\daN^* B^*_z(0)$ and $I+\daP C_z(0)$ on $D_{\rho,\theta}$ }

We prove here the
\begin{theorem}\label{singalriteenSQRT}The functions $\sqrt{1-z}\,(I+\daN^* B_z(0))$ and  $\sqrt{1-z}\,(I+\daP C_z(0))$ admit  an analytic expansion on $D_{\rho,\theta}$ and may be continuously extended on $\overline{D_{\rho,\theta}}$. Furthermore, one gets

\begin{equation}\label{eq26}
 \lim_{z\rightarrow1}\sqrt{1-z}\,(I+\daN^* B^*_z(0))=\sqrt{\frac{k''(0)}{2}}\,A_-,
\end{equation}
\begin{equation}\label{eq27}
 \lim_{z\rightarrow1}\sqrt{1-z}\,(I+\daP C_z(0))=-\sqrt{\frac{k''(0)}{2}}\,A_+.
\end{equation}
\end{theorem}

\begin{proof}
First case : when $z\in D_{\rho,\theta}\setminus K(\delta_1,\delta_2)$, the analysis of $z\mapsto\sqrt{1-z}\,(I+\daP C_z(0))$ (resp. $z\mapsto\sqrt{1-z}\,(I+\daN^* B^*_z(0))$) is derived from Corollary \ref{Rhodelta} and the fact that $z\mapsto(I-zP(0))^{-1}$ is analytic in $D_{\rho,\theta}\setminus K(\delta_1,\delta_2)$.\\
   Second case :  the map $z \mapsto \sqrt{1-z} (I-\daP B^*_z(0))^{-1} $  is the analytic expansion on $ K^*(\delta_1,\delta_2)$ of $z\mapsto \sqrt{1-z}\,(I+\daP C_z(0))$, and,   by  (\ref{eq59}), one gets
       $$ \sqrt{1-z} (I-\daP B^*_z(0))^{-1} = \sqrt{1-z}\,B_+^{-1}(z,0)F_+(z,0).$$
    By Proposition \ref{propo2}, the map  $z\mapsto B_+^{-1}(z,0)$ is analytic on $K^*(\delta_1,\delta_2)$ and one gets
   $$\sqrt{1-z}\,F_+(z,0)=\sqrt{1-z}\,\left(I-\frac{\l_+(z)-a_+}{\l_+(z)}\,\Pi_+(z)\right),$$
so that    $\displaystyle \lim_{z\rightarrow1}\sqrt{1-z}\,F_+(z,0)$ exists
   since $ \frac{\sqrt{1-z}}{\l_+(z)}\to  \sqrt{\frac{k''(0)}{2}}$ as
   $z\to 1$. Hence,
$$z\mapsto\sqrt{1-z}\,(I+\daP C_z(0))$$
 is analytic on $D_{\rho,\theta}$ and admits an analytic expansion on the boundary of $D_{\rho,\theta}$.

\end{proof}

\section{Proofs of the local limit theorems}

This section is devoted to the proof of our local limit
theorems \ref{theoloc-min-laplace}, \ref{theoloc-min-repartition} and \ref{thm3}.
\subsection{Preliminaries}

In the previous section, we have described  the local behavior near $z=1$ of a family of analytic functions, expressed in terms of Laplace transforms ; we thus need some argument which relies the type of singularity  near $z=1$ of such a function to its behavior at infinity.  The following lemma is a classical result in the theory of complex variables functions.

\begin{lemma}[\cite {flajolet}]\label{ELP&MP}
If a function $z\mapsto G(z)$ satisfies simultaneously the following three conditions:
\begin{itemize}
  \item  $G$ is analytic on $D_{\rho,\theta}$ and can be written as $G(z)=\sum_{n=0}^{+\In}g_n z^n$;
  \item  $\sqrt{1-z}G(z)$ is bounded in $D_{\rho,\theta}$;
  \item  $\lim_{z\rightarrow1}\sqrt{1-z}G(z)=C>0,$\\
  \end{itemize}
then
  $$g_n\sim\frac{C}{\sqrt{\pi n}},\qquad n\rightarrow+\In.$$
\end{lemma}

\begin{proof} For the sake of completeness, we detail here the proof. For every $\var>0$, $r\in ]1, \rho[$ and $\theta'>\theta$, let's consider  the  arcs $\ga_0=\ga_0(\var, \th'),  \ga_1=\ga_1(\var,r'), \ga'_1=\ga'_1(\var,r')$ and $\ga_2=\ga_2(r)$   defined respectively by
\begin{alignat}{3}
 \ga_0&:=\{z=1+\var\e^{-it}; \th'\leq t\leq 2\pi-\th'\};\\
 \ga_1&:=\{z=1+t\e^{i\th'}; \varepsilon\leq t\leq r'\}\quad {\rm and} \quad \ga'_1:=\{z=1+(r'-t)\e^{i(2\pi-\th')};0\leq t\leq r'-\var\};\\
 \ga_2&:=\{z=r\e^{it}; \th'' \leq t\leq 2\pi-\th'' \},
\end{alignat}
where $r'$ and $\th''$ verify the following system of equations:

$$\left\{
            \begin{array}{ll}
             r\cos\th''=1+r'\cos\th'; \\
             r\sin\th''=r'\sin\th'.
            \end{array}
          \right.
$$

\begin{figure}
 \begin{center}
  \includegraphics[width=540pt,height=289pt]{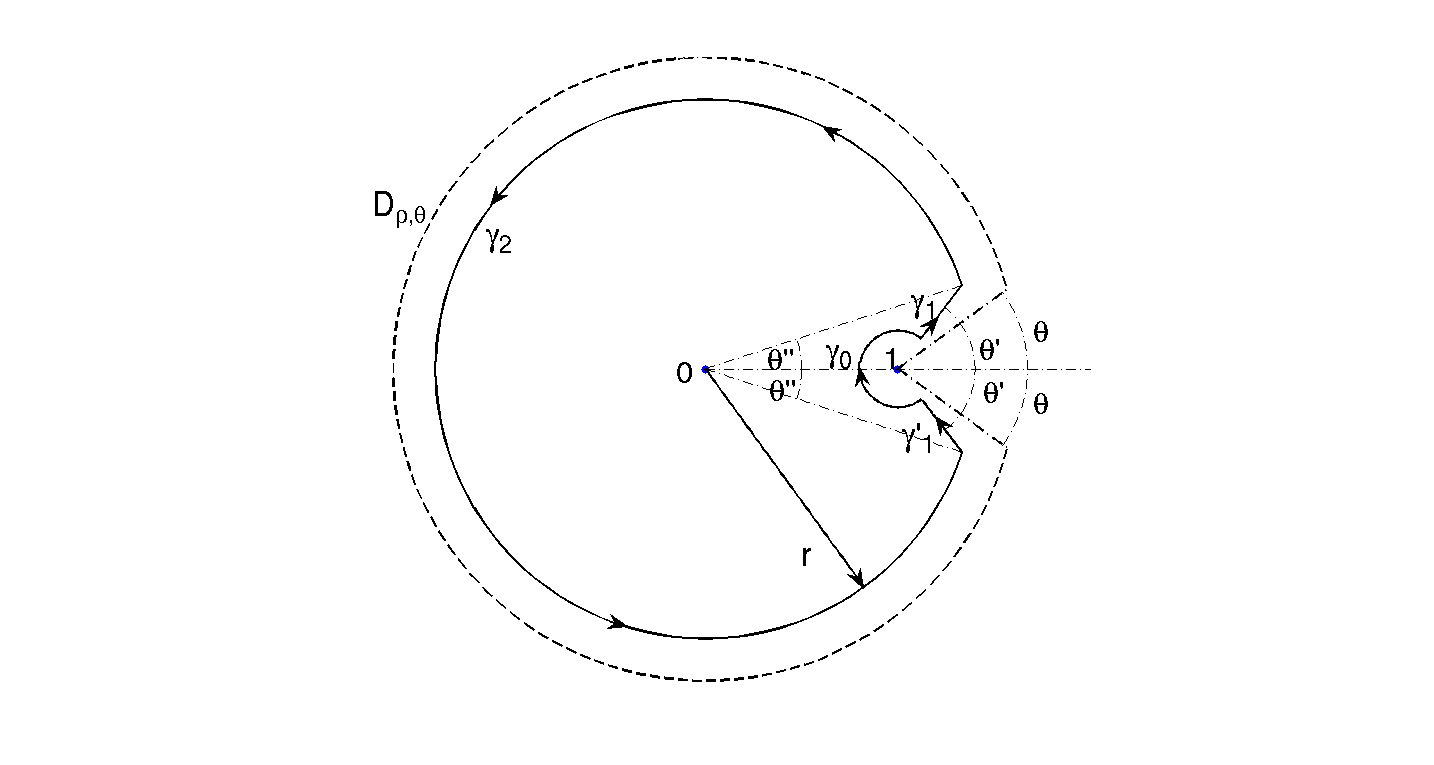}
 \end{center}
 \caption{The closed path $\ga_0\cup\ga_1\cup\ga_2\cup\ga'_1$ of Lemma \ref{ELP&MP} and the open set $D_{\rho,\theta}$.}\label{fig-dessin3}
\end{figure}

 Define a closed path $\gamma(\var, r)$, composed by the curves $\ga_0$, $\ga_1, \ga_2$ and $\ga'_1$, as showed in Figure \ref{fig-dessin3}.
We now introduce the  complex function $F(z)$ defined by
 $$F(z)=G(z)-\frac{C}{\s{1-z}}:=\frac{\delta(z)}{\s{1-z}}.$$

 Since $z\mapsto G(z)$ is analytic on $D^\circ_{\rho,\theta}$, so is $F$  on this set  and one may   write, for  $|z|<1 $
 $$F(z)=\sum_{n=0}^{+\In}f_n\,z^n  $$
 where $f_n=\frac{1}{2\pi i}\int_\ga\frac{F(z)}{z^{n+1}}\d z$ doest not depend on $\var$, $r$ and $\th$. By hypothesis, there exists some constant $M>0$ such that  $|F(z)|\leq\frac{M}{|\s{1-z}|}$ for $z\in D^\circ_{\rho,\th}$ ; one thus gets
  $$\frac{1}{2\pi}\int_{\ga_0}\frac{|F(z)|}{|z|^{n+1}}\d z\leq\frac{M}{2\pi}\int_{\th'}^{2\pi-\th'}\frac{\s{\var}}{(1+\var\e^{it})^{n+1}}\d t\leq  \frac{ M \s{\var}}{|1-\var|^{n+1}}$$
and
 $$\frac{1}{2\pi}\int_{\ga_2}\frac{|F(z)|}{|z|^{n+1}}\d z\leq \frac{M}{r^{n}\s{ r-1}}.$$
On the other hand,
$$\frac{1}{2\pi}\int_{\ga_1\cup \ga'_1}\frac{|F(z)|}{|z|^{n+1}}\d z\leq\frac{\delta(\var, r')}{\pi}\int_0^{r'}\frac{\d t}{\s{t}(1+t\cos \th')^{n+1}},$$
with $\delta(\var, r'):=\displaystyle  \sup_{z\in\ga_1 \cup \ga'_1 }|\delta(z)|$.

Set $\displaystyle I_n(r'):=\int_0^{r'} \frac{\d t}{\s{t}(1+t\cos \th')^{n+1}}$.
 Since $\ln(1+u) \geq{\ln r'\over r'} u$ as soon as $0\leq u\leq r'$, for any $t \in [0, r']$ one gets
$\displaystyle \ln (1+t \cos \th') \geq {\ln r' \over r'}t \cos \th'$, so that
$$
I_n(r')
=\int_0^ {r'}  {e^{-(n+1)\ln  (1+t\cos \th')} \over \sqrt{t}}\d t\leq
\int_0^{r'}  {e^{-(n+1){\ln  r'\over r'}t \cos \th'} \over \sqrt{t}}\d t.
$$
Setting
  $s=(n+1)t$, one obtains $\displaystyle I_n(r')\leq {1\over \sqrt{n+1}}
  \int_0^{+\infty}  {e^{-{\ln  r'\over r'}s\cos \th'} \over \sqrt{s}}\d s$, i.e. $\s{n}\ I_n(r') \leq M'$ for some constant
  $M' \in ]0, +\infty[$ ; this readily implies
  $\displaystyle
  \frac{1}{2\pi}\int_{\ga_1\cup\ga'_1}\frac{|F(z)|}{|z|^{n+1}}\d z\leq \delta(\varepsilon, r') M'.
  $
In summary one gets
$$
\sqrt{n} \vert f_n \vert \leq \frac{ M \s{\var n}}{|1-\var|^{n+1}}+ \frac{M\s{n}}{r^{n}\s{ r-1}}+\delta(\varepsilon, r') M',
$$
so that
$\displaystyle \sqrt{n} \vert f_n \vert \leq  \frac{M\s{n}}{r^{n}\s{ r-1}}+\delta(0, r') M'$ since $\varepsilon$ may be choosen arbitrarily small.
Letting  now $n\to +\infty $,  one gets,   since $r>1$
$$\limsup_{n\rightarrow+\infty} \sqrt{n} \vert f_n \vert \leq \delta(0, r')M'
$$
and one concludes that $\s{n} f_n\to 0$ as $n \to +\infty$ noticing that $ \lim_{r'\to 0} \delta(0, r')=0$.\\

One achieves the proof writing $ g_n=f_n+Ca_n$  with $a_n=\frac{2n!}{4^n(n!)^2}=\frac{1+o(n)}{\s{\pi n}}$, so that
$$g_n=f_n+Ca_n\sim\frac{C}{\s{\pi n}} \quad\text{as }n\rightarrow+\In.$$
\end{proof}

\subsection{Proof of Theorem \ref{theoloc-min-laplace}}
We fix $\lambda >0$ and set, for any $i, j \in E$ and $z \in \mathbb D^\circ$,
$$
G_{i, j}(z, \l):= \sum_{n=0}^{+\In}z^n\E_i(\e^{\l m_n};X_n=j)
$$
and
$$
H_{i,j}(z,\l):=\sqrt{1-z}\ G_{i, j}(z, \l).
$$
By lemma \ref{lem4.1}, we have
\begin{equation*}
H_{i,j}(z,\l) =\left\{[I+\daN^* B^*_z(\l)]\;\sqrt{1-z}[I+\daP C_z(0)]\right\}_{i,j}.
\end{equation*}
By  (\ref{eq27}), we get
\begin{equation}\label{eq31}
H_{i,j}(\l):=\lim_{z\rightarrow1}H_{i,j}(z,\l)=-\sqrt{\frac{k''(0)}{2}}\{(I+\daN^* B^*_1(\l))A_+\}_{i,j}.
\end{equation}
By (\ref{eq29}) and (\ref{eq28}), we obtain
\begin{equation}\label{eq500}
\lim_{\l\rightarrow0^+}\l H_{i,j}(\l)=-\sqrt{\frac{k''(0)}{2}}{\left(A_-A_+\right)}_{i,j}=\sqrt{\frac{2}{k''(0)}}\left(\Pi(0)\right)_{i,j}=\sqrt{\frac{2}{k''(0)}}\ \nu_j>0.
\end{equation}
%
From (\ref{eq30}), the coefficients of $A_+$ are $\leq 0$, the function $H_{i,j}$ is in fact the Laplace transform of a positive  measure $\mu_{i,j}$ on $\R_-$  and this measure is $\neq 0$ by (\ref{eq500}) ; in particular, there exists an interval $[a,b]\subset\R_-$ such that $\mu_{i,j}([a,b])>0$. Therefore, for all $\l>0$, one gets
\begin{equation}\label{eq32}
H_{i,j}(\l)=\int_{-\In}^0\e^{\l x}\d\mu_{i,j}(x)\geq\int_{a}^b\e^{\l x}\d\mu_{i,j}(x)\geq\e^{\l a}\mu_{i,j}([a,b])>0.
\end{equation}
Consequently, by the above, for any $\l>0$, the function $z\mapsto G_{i,j}(z,\l)$ is analytic on $D_{\rho,\theta}$, $z \mapsto \sqrt{1-z}G_{i, j}(z, \lambda)$ is  bounded on $D_{\rho,\theta}$ and $\lim_{z\rightarrow1}H_{i,j}(z,\l)>0$. By Lemma \ref{ELP&MP}, we obtain
\begin{equation}\label{eq41*}
\sqrt{n}\,\E_i(\e^{\l m_n},X_n=j)\stackrel{n\rightarrow+\In}{\longrightarrow}\frac{H_{i,j}(\l)}{\sqrt{\pi}}.
\end{equation}


\subsection{Proof of Theorem \ref{theoloc-min-repartition}}

 In this paragraph, we precise the previous statement in terms of distribution function. We thus introduce, for any
 any $(i,j)\in E\times E$,   the distribution function $h_{i,j}:\R_+\rightarrow\R$ of the measure $\mu_{i, j}$,  defined by
$$h_{i,j}(x)=\left\{
            \begin{array}{ll}
              -\sqrt{\frac{k''(0)}{2\pi}}\left\{[I+\daN^* B^*_1(1_{[-x,0]})]A_+\right\}_{i,j}, & x > 0; \\
              -\sqrt{\frac{k''(0)}{2\pi}}\;(A_+)_{i,j}, & x=0;
            \end{array}
          \right.
$$
where $\displaystyle \daN^* B^*_1(1_{[-x,0]})=\sum_{n=1}^{+\In}\P_i(S_1>S_n,S_2>S_n,\cdots, S_{n-1}>S_n,-x\leq S_n\leq0,X_n=j)$, for $x>0$.  We will decompose the `` potential '' $ \daN^* B^*_1(1_{[-x,0]})$ in terms of the ladder epochs $\{\t_j\}_{j\geq0}$ of the random walk $(S_n)_n$, defined recursively by :
$$\t_0=0\quad\text{and}\quad \t_j=\inf\{n;\text{ for all }n\geq\t_{j-1},S_n<S_{\t_{j-1}}\},\text{ for }j\geq1.$$
For any $x \in \mathbb R^{*+}$ and   $l\geq0$, we thus consider the matrix $B_l(x)$ defined by
$$B_{l}(x)=\left({B_{l}(x)}_{i,j}\right)_{i,j},$$
with ${B_{l}(x)}_{i,j}=\sum_{k\in E}\P_i(S_{\t_l}\geq-x,X_{\t_l}=k)(A_+)_{k,j}$.

One gets
$$
h_{i, j}(x)= -\sqrt{\frac{k''(0)}{2\pi}}\sum_{l\geq 0} B_l(x)_{i, j} =-\sqrt{\frac{k''(0)}{2\pi}}\sum_{k\in E} \E_i\Bigl[\sum_{l\geq 0} 1_{[-x, 0]}(S_{\tau_l}), X_{\tau_l}=k\Bigr] (A_+)_{k,j}.
$$
Notice that, for $x $ large enough, one gets $\displaystyle \E_i\Bigl[\sum_{l\geq 0} 1_{[-x, 0]}(S_{\tau_l}), X_{\tau_l}=k\Bigr] >0$ for any $i, k \in E$ since $S_{\tau_1}$ is finite $\P_i$-a.s. ;  so is $ h_{i, j}(x)$, since by \ref{eq28} at least one of the terms
$(A_+)_{k,j}$ is non negative.  We will see that this   property holds in fact for any $x\geq 0$.

First, one gets  the

\begin{lemma}\label{lem4.1.4}
For any $(i,j)\in E\times E$, we have
$$\sqrt{n}\ \P_i(m_n=0,X_n=j)=\sqrt{n}\ \P_i(T^*_->n,X_n=j)\longrightarrow\left(-\sqrt{\frac{k''(0)}{2\pi}}\right)(A_+)_{i,j},\quad\text{as }n\rightarrow+\In.$$
\end{lemma}
\begin{proof}
 Indeed, (\ref{eq27}) may be restated as follows
\begin{equation*}
 \begin{split}
   \sqrt{1-z}\left[I+\sum_{n=1}^{+\In}z^n\P_i(m_n=0,X_n=j)\right]
                                                                =&\sqrt{1-z}\,[I+\daP C_z(0)]\stackrel{z\rightarrow1}{\longrightarrow}-\sqrt{\frac{k''(0)}{2}}(A_+)_{i,j}.
 \end{split}
\end{equation*}
so that, by Lemma \ref{ELP&MP} (when $-(A_+)_{i,j}>0$),
\begin{equation}\label{eq33*}
\sqrt{n}\ \P_i(m_n=0,X_n=j)\stackrel{n\rightarrow+\In}{\longrightarrow}-\sqrt{\frac{k''(0)}{2\pi}}(A_+)_{i,j}.
\end{equation}
The same result holds when $-(A_+)_{i,j}=0$, by Corollary 1 in \cite{flajolet}.
\end{proof}

We will   use the following
\begin{lemma}\label{lem4.1.1}
For any $l \geq 0$, any $i, j \in E$ and  $x>0$ such that $h_{i,j}$ is discontinuous at $x$, we have
\begin{equation}\label{liminf}
\liminf_{n\rightarrow+\In}\sqrt{n}\,\P_i(S_{\t_l}\geq-x,\t_l\leq n,\t_{l+1}>n,X_n=j)\geq-\sqrt{\frac{k''(0)}{2\pi}}B_l(x)_{i,j}.
\end{equation}
\end{lemma}
\begin{proof} For any $0<\de<1$, we have
\begin{equation}\label{eq4.1.2}
 \P_i(S_{\t_l}\geq-x,\t_l\leq n,\t_{l+1}>n,X_n=j)\geq\P_i(S_{\t_l}\geq-x,\t_l\leq\de n,\t_{l+1}>n,X_n=j).
\end{equation}
From Markov property, we have
\begin{equation}\label{eq99*}
 \begin{split}
  \P_i(S_{\t_l}\geq-x,\t_l\leq\de n,\t_{l+1}>n,X_n=j)=&\sum_{k\in E}\E_i\left[(S_{\t_l}\geq-x,\t_l\leq\de n,X_{\t_l}=k)\P_k(\t_1>n,X_{n-\t_l}=j)\right]\\
                                                     =&\sum_{k\in E\atop 0\leq p\leq\de n}\E_i\left[(S_p\geq-x,\t_l=p,X_p=k)\P_k(\t_1>n,X_{n-p}=j)\right].
 \end{split}
\end{equation}
In addition, one gets
$$\sqrt{n}\,\P_k(\t_1>n,X_{n-p}=j)=\sqrt{n}\,\P_k(\t_1>n-p,X_{n-p}=j)-\sqrt{n}\,\P_k(n-p<\t_1\leq n,X_{n-p}=j).$$
Since $0\leq\P_k(n-p<\t_1\leq n,X_{n-p}=j)\leq\P_k(n-p<\t_1\leq n)$, by Lemma \ref{lem4.1.4}, we hence obtain that
$$\sqrt{n}\,\P_k(n-p<\t_1\leq n)=\sqrt{n}\, \P_k(\t_1>n-p)-\sqrt{n}\,\P_k(\t_1>n)\rightarrow0,\quad\text{as }n\rightarrow+\In.$$
So we have $\lim_{n\rightarrow+\In}\sqrt{n}\,\P_k(n-p<\t_l\leq n,X_{n-p}=j)=0$. By lemma \ref{lem4.1.4}, we get
\begin{equation*}
  \lim_{n\rightarrow+\In}\sqrt{n}\,\P_k(\t_l>n,X_{n-p}=j)=-\sqrt{\frac{k''(0)}{2\pi}}\;(A_+)_{k,j}.
\end{equation*}
Using Fatou's lemma, the inequalities (\ref{eq4.1.2}) and (\ref{eq99*}), one concludes
\begin{eqnarray*}
 \liminf_{n\rightarrow+\In}\sqrt{n}\,\P_i(S_{\t_l}\geq-x,  &\t_l\leq n,&  \t_{l+1}>n,X_n=j)\\
  &\geq&\liminf_{n\rightarrow+\In}\sqrt{n}\,\P_i(S_{\t_l}\geq-x,\t_l\leq\de n,\t_{l+1}>n,X_n=j)\\
  &\geq&\left(-\sqrt{\frac{k''(0)}{2\pi}}\right)\sum_{k\in E}\P_i(S_{\t_l}\geq-x,X_{\t_l}=k)(A_+)_{k,j}\\
 & =&\left(-\sqrt{\frac{k''(0)}{2\pi}}\right)B_l(x)_{i,j}.
\end{eqnarray*}

\end{proof}


\begin{proof} [Proof of Theorem \ref{theoloc-min-repartition}]
From $(\ref{eq41*})$ and the extended continuity theorem (Thm 2a, \Rmnum{13}.1, W. Feller \cite{Fell}), for any $(i,j)\in E\times E$ and any $x>0$ such that $h_{i,j}(\cdot)$ is continuous at $x$, one gets
$$\lim_{n\rightarrow+\In}\sqrt{n}\,\P_i(m_n\geq-x,X_n=j)=h_{i,j}(x) ; $$
By Lemma \ref{lem4.1.4}, the same result holds for  $x=0$.\\

Now, fix $x>0$ such that $h_{i,j}(\cdot)$ is discontinuous at $x$. The map  $x\mapsto h_{i,j}(x)$  being  increasing and right-continuous  on $\R^*_+$, the set of its points of discontinuity is countable and there thus exists a sequence $(\var_k)_{k\geq 1}$  of non negative reals converging towards $0$ and such that such
  $h_{i, j}$ is continuous at $x+\var_k$  for any $k\geq 1$; consequently, for any $k\geq 1$ one gets
$$\sqrt{n}\,\P_i(m_n\geq-x,X_n=j)\leq\sqrt{n}\,\P_i(m_n\geq-x-\var_k,X_n=j)$$
and so
$$\limsup_{n\rightarrow+\In}\sqrt{n}\,\P_i(m_n\geq-x,X_n=j)\leq h_{i,j}(x+\var_k).$$
The map $h_{i,j}$ being right continuous, one gets
\begin{equation}\label{eq42*}
\limsup_{n\rightarrow+\In}\sqrt{n}\,\P_i(m_n\geq-x,X_n=j)\leq h_{i,j}(x).
\end{equation}
On the other hand,  for any $N\leq n$ and $0\leq l<N$, one gets
$$\P_i(m_n\geq-x,X_n=j)\geq\sum_{l=0}^{N}\P_i(S_{\t_l}\geq-x,\t_l\leq n,\t_{l+1}>n,X_n=j)$$
which readily implies, by Lemma \ref{lem4.1.1}
\begin{equation}\label{eq43*}
 \begin{split}
  \liminf_{n\rightarrow+\In}\sqrt{n}\,\P_i(m_n\geq-x,X_n=j)&\geq\sum_{l=0}^{N}  \liminf_{n\rightarrow+\In}\sqrt{n}\P_i(S_{\t_l}\geq-x,\t_l\leq n,\t_{l+1}>n,X_n=j)
  \\
  &\geq \left(-\sqrt{\frac{k''(0)}{2\pi}}\right)\sum_{l=0}^N B_{l}(x)_{i,j}\\
   &=\left(-\sqrt{\frac{k''(0)}{2\pi}}\right)\sum_{l=0}^N\left(\sum_{k\in E}\P_i(S_{\t_l}\geq-x,X_{\t_l}=k)(A_+)_{k,j}\right)\\
   &\stackrel{N\rightarrow+\In}{\longrightarrow}\left(-\sqrt{\frac{k''(0)}{2\pi}}\right)\left[\left(I+\daN^* B^*_1(1_{[-x,0]})\right)A_+\right]_{i,j}=h_{i,j}(x).
 \end{split}
\end{equation}
Combining (\ref{eq42*}) and (\ref{eq43*}), one gets the expected conclusion at $x$.

Now we are going to prove that for any $ j \in E $, the function $ (x, i) \mapsto h_{i,j}(x)$ is   harmonic  with respect to  $(S_n, X_n)$ and  positive on $\R\times E$. One gets
\begin{equation}\label{eq101*}
\sqrt{n+1}\,\P_i(m_{n+1}\geq -x,X_{n+1}=j)=\sqrt{\frac{n+1}{n}}\sum_{i_1\in E}\int p_{i,i_1}\sqrt{n}\P_{i_1}(m_n\geq -x-y_1,X_{n+1}=j)F(i,i_1,\d y_1),
\end{equation}
with $\displaystyle\E_i(x+|Y_1|)=x+\sum_{j\in E}p_{i,j}\int_\R|u|F(i,j,\d u)<+\In.$ We now need the
\begin{lemma}\label{lemm1}
 There exists a constant $C>0$  such that
 \begin{equation}\label{eq33}
\forall i, j \in E,\ \forall x \geq 0, \qquad  \s{n} \  \P_i(m_n\geq-x,X_n=j)\leq C(x+1).
 \end{equation}
\end{lemma}
By the dominated convergence theorem, tending $n\rightarrow+\In$ in (\ref{eq101*}), one thus gets
\begin{equation}\label{eq35}
 \forall x \geq 0,\quad h_{i,j}(x)=\sum_{i_1\in E}\int p_{i,i_1}h_{i_1,j}(y_1+x)F(i,i_1,\d y_1)=\E_i[h_{X_1,j}(x+Y_1)],
\end{equation}
which means  that $(x, i)\mapsto h_{i,j}$ is harmonic for $( S_n, X_n)$ on $\R^+\times E$.

By equality (\ref{eq000}) of Theorem \ref{theoloc-min-laplace}, for $\l>0$, one gets $\displaystyle \lim_{\tau\rightarrow0^+}\frac{H_{i, j}(\tau\l)}{H_{i, j}(\tau)}=\frac{1}{\l} $
and the classical  Tauberian theorem (see for instance Thm 1, \Rmnum{13}.5, W. Feller \cite{Fell}), we get
\begin{equation}\label{eq37}
h_{i,j}(x)=\mu_{i,j}([-x,0])\sim\frac{H_{i,j}(1/x)}{\Gamma(2)}\sim\sqrt{\frac{2}{k''(0)}}\ \nu_j\ x \quad {\rm as}\quad x\rightarrow+\In.
\end{equation}

At last, assume that there exists $(i_0,j_0)\in E\times E$  such that $h_{i_0,j_0}(0)=0$. Iterating Formula (\ref{eq35}), one gets   for any $n\geq1$,
$$
h_{i_0,j_0}(0)=\E_{i_0}[h_{X_n,j_0}(S_n)],
$$
so that
\begin{equation}\label{eq38}
 h_{X_n,j_0}(S_n)=0 \qquad \P_{i_0} \text{ a.s.}
\end{equation}
By (\ref{eq37}), there exists $M_{j_0}\geq1$, such that for $x\geq M_{j_0}$,
\begin{equation}\label{eq39}
\inf_{i\in E}h_{i,j_0}(x)\geq\frac{1}{2}\sqrt{\frac{2}{k''(0)}}\ \nu_{j_0}>0
\end{equation}
and the central limit theorem for Markov chains (\cite{guivarc'h-hardy}) implies that for any $i\in E$,
$$
\P_i\left(\frac{ S_n}{\sqrt{n}}\geq M_{j_0}\right)\stackrel{n\rightarrow+\In}{\longrightarrow}\frac{1}{\sqrt{\pi k''(0)}}\int_{M_{j_0}}^{+\In}e^{-\frac{x^2}{2k''(0)}}\d x:=\a(M_{j_0})>0.
$$
Setting $\displaystyle B_n=\Bigl\{\omega : \frac{ S_n(\omega)}{\sqrt{n}}\geq M_{j_0}\Bigr\}$ and $\displaystyle B=\limsup_{n\to +\infty}B_n$, then for all $i\in E$, one thus may write
$$\P_i(B)=\lim_{m\rightarrow+\In}\P_i\Bigl(\bigcup_{n\geq m}B_n\Bigr)\geq\lim_{m\rightarrow+\In}\P_i(B_m)=\a(M_{j_0})>0.$$

For all $\omega \in B$, one gets $\displaystyle \limsup_{x\rightarrow+\In}\left[ S_n(\omega)\right]=+\In$ and so, by
(\ref{eq39}), one obtains
$$
\limsup_{n\rightarrow+\In}\inf_{i\in E}\left[h_{i,j_0}(x_0+S_n)1_B\right]\geq\frac{1}{2}\sqrt{\frac{2}{k''(0)}}\,\nu_{j_0}>0 \quad \P_{i_0}\text{-a.s.}
$$
This contradicts (\ref{eq38}) since $\P_i(B)>0$, for any $i\in E$.
Then, for any $i, j \in E$ and $x\geq 0$ one gets $h_{i, j}(x)\geq
h_{i, j}(0)>0$.
\end{proof}

It remains to prove Lemma \ref{lemm1}; we will use the two following facts, whose proofs may be found in  \cite{koz}  :

\begin{fact}[\cite{koz}]\label{lem4.6}
 Let $c, \nu\in \R^*_+$ and $(a_n)_{n\geq0}$ be a monotone sequence of non negative reals such that $\displaystyle \sum_{n=0}^{+\In}a_ns^n\leq c(1-s)^{-\nu}$
 for any $s\in[0,1[$.
 Then
 $$\forall n \geq 2, \qquad a_n\leq c\e(1-\e^{-1})^{-\nu}2^{1+\nu}n^{\nu-1}.$$
\end{fact}

 \begin{fact}[\cite{koz}]\label{lem4.5}
 Let $H$ be a non-decreasing function on $\R^+$ such that $H(0)=0$ and  the integral
  $\displaystyle \widetilde{H}(\l):=\int_{0}^{+\In}\e^{-\l x}\d H(x)$ does exist for any $\l>0$.
  If there exist   $\delta, \ga >0$ such that
  $$\forall \l\in]0,\delta],\quad \widetilde{H}(\l)\leq c\l^{-\ga},$$
  then, for all $x\geq\delta^{-1}$, one gets $H(x)\leq c\ \e\  x^\ga.$
\end{fact}

\begin{proof}[Proof of Lemma \ref{lemm1}]
Taking into account (\ref{eq500}), we get for any $i\in E$,
\begin{equation*}
  \lim_{\l\rightarrow0}\l\sum_{j\in E}H_{i,j}(\l)=\sqrt{\frac{2}{k''(0)}}>0,
\end{equation*}
which implies that there exist two constants $\delta>0$ and $c>0$ such that for any $\l \in ]0, \delta]$ and $s\in ]0, 1[$,
$$\sup_{i\in E}\sum_{n=0}^{+\In}s^n\E_i(\e^{\l m_n})\leq c\l^{-1}(1-s)^{-1/2}.$$
For $\l>0$, the sequence $\Bigl(\E(\e^{\l m_n})\Bigr)_{n\geq0}$ is decreasing with respect to $n$ and the Fact $\ref{lem4.6}$ with $\nu=1/2$ leads to
$$\forall i \in E,\ \forall n \geq 2,\ \forall \l \in ]0, \delta],  \quad  \sqrt{n}\ \E_i(\e^{\l m_n})\leq c\e(1-\e^{-1})^{-1/2}2^{3/2}\l^{-1}.$$
Applying now Fact \ref{lem4.5} with $\ga =1$,  we get, for all $x\geq\delta^{-1}>0$, $n\geq2$ and $i, j\in E$,
$$\sqrt{n}\ \P_i(m_n\geq-x,X_n=j)\leq\ \sqrt{n}\ \P_i(m_n\geq-x)\leq c_1 x,$$
where $c_1=c\e^2(1-\e^{-1})^{-1/2}2^{3/2}$.

On the other hand,  for $0\leq x<\delta^{-1}$, one gets
$$\sqrt{n}\ \P_i(m_n\geq-x,X_n=j)\leq\sqrt{n}\ \P_i(m_n\geq-\delta^{-1},X_n=j)\stackrel{n\rightarrow+\In}{\longrightarrow}h_{i,j}(\delta^{-1})$$
and one thus may write, , for any $i,j\in E$ and $x\geq0$,
$$\sqrt{n}\ \P_i(m_n\geq-x,X_n=j)\leq c_1x+c_2$$
where $\displaystyle c_2=\sup_{\stackrel{n\geq 1}{i,j\in E}}\P_i(m_n\geq-\delta^{-1},X_n=j)$.
\end{proof}

We end this section with the following elementary consequence of the above :
\begin{fact}\label{ineqharmonic} There exists a constant $c\geq 1$ such that,
for any $i, j \in E$ and $x \geq 0$ one gets
$$
{x+1\over c} \leq h_{i, j}(x) \leq c(x+1)
$$
\end{fact}
\begin{proof}
By (\ref{eq37}), there exists $c_1>0$ and $x_1\geq 0$ such that $
{x+1\over c_1} \leq h_{i, j}(x) \leq c_1(x+1)
$
for $x \geq x_1$.
For $0\leq x\leq x_1$ one thus gets
$$
{h_{i, j}(0)\over c_1 h_{i, j}(x_1)} (1+x) \leq h_{i, j}(0) \ \leq h_{i, j}(x) \leq   h_{i, j}(x_1) \leq c_1(x+1)(x_1+1).
$$
and one set $c:= \max(c_1(x_1+1), c_1 {h_{i, j}(x_1)\over h_{i, j}(0)})$.
\end{proof}

\subsection{Proof of Theorem \ref{thm3}}
\begin{proof}
By the Markov property and Fubini's theorem, we have, for $0<\varepsilon<\l$,
\begin{equation*}
   \begin{split}
    &\sum_{n=0}^{+\In}z^n\E_i(\e^{\l m_n-\varepsilon
    S_n},X_n=j)\\
               =&\sum_{n=0}^{+\In}z^n[\delta_{i,j}+\sum_{k=1}^{n}\E_i(\e^{\l S_k-\varepsilon S_n},S_0>S_k,\cdots,S_{k-1}>S_k, S_{k+1}\geq S_k, \cdots,S_n\geq S_k,X_n=j)]\\
                =&\sum_{n=0}^{+\In}z^n\Big\{\delta_{i,j}+\sum_{k=1}^{n}\sum_{l\in
    E}\E_i\Big[\e^{(\l-\varepsilon)S_k},S_0>S_k,\cdots,S_{k-1}>S_k,X_k=l\Big]\times
   \\
   &\qquad\qquad\qquad\qquad\qquad\qquad\qquad\qquad\qquad\qquad\qquad\E_l[\e^{-\varepsilon S_{n-k}},S_1\geq0,\cdots,S_{n-k}\geq0,X_{n-k}=j]\Big\}\\
                =&\sum_{l\in E}\Big[\sum_{k=0}^{+\In}z^k\E_i(\e^{(\l-\varepsilon)S_k};S_1>S_k,\cdots,S_{k-1}>S_k,S_k<0,X_k=l)\Big]\times\\
   &\qquad\qquad\qquad\qquad\qquad\qquad\qquad\qquad\qquad\qquad\qquad\Big[\sum_{p=0}^{+\In}z^p\E_l(\e^{-\varepsilon S_p};S_1\geq0,\cdots,S_p\geq0,X_p=j)\Big]\\
                =&\Big\{(I+\daN^* B^*_z(\l-\varepsilon))(I+\daP C_z(-\varepsilon))\Big\}_{i,j}.
   \end{split}
  \end{equation*}
So by the first assertion of Theorem \ref{NOsingularite}, letting $z \to 1$, one obtains
$$\sum_{n=0}^{+\In}\E_i(\e^{\l m_n-\varepsilon
    S_n},X_n=j)=\{(I+\daN^* B^*_1(\l-\varepsilon))(I+\daP C_1(-\varepsilon))\}_{i,j}<+\In.$$
\end{proof}

\section{Appendix}
  \subsection{Absolutely continuous components for $k$ times convolution of a matrix of positive measures on $\R$}\label{appendixA}

We use here the Notations \ref{notation2.2.1} and we prove the

\begin{lemma}\label{lem2.6}
Assume that $M(\d x)=(\mu_{i,j}(\d x))_{1\leq i,j\leq N}$ is a matrix of
positive measures on $\R$. If the following two conditions hold
simultaneously:
\begin{enumerate}
 \item\label{con1} there exist $(i_0,j_0)\in \{1,\cdots,N\}^2$ and $n_0\geq1$, such that $\mu_{i_0,j_0}^{(n_0)}(\d x)$ has an absolutely continuous component;
  \item\label{con2} there exists $n_1\geq1$, such that $M^{\bullet n_1}(\R)>0$,
\end{enumerate}
then for any $k\geq(n_0+1)n_1n_0$, one gets $M^k(\R)>0$ and there exists at least one
absolutely continuous component term in $M^{\bullet k}$.
\end{lemma}
\begin{proof} There are two cases to consider.
 \begin{description}
 \item[Case 1] \underline{$i_0=j_0$}. The matrix
 $M^{\bullet n_0}$ has thus an absolutely continuous component term in its
 diagonal. Since $M^{\bullet n_0n_1}=(M^{\bullet n_0})^{n_1}$, it is clear that
 there also exists an absolutely continuous component term on the
 diagonal of the matrix $M^{\bullet n_0n_1}$. Moreover, one gets $M^{\bullet n_1}(\R)>0$, so
that $M^{\bullet n_0n_1}(\R)=(M^{\bullet n_1}(\R))^{n_0}>0$.\\
 Consequently the matrix $M^{\bullet n_0n_1}$ has an absolutely continuous component term on its
 diagonal and $M^{\bullet n_0n_1}(\R)>0$. This implies that for any
 $k\geq(n_1+1)n_0n_1>n_0n_1$, the matrix $M^{\bullet k}$ has at least one
 absolutely continuous component term and $M^{\bullet k}(\R)>0$.

 \item[Case 2]\underline{$i_0\neq j_0$}. Set
 $n'_1=(n_1+1)n_0$. The positivity of $M^{\bullet n_0n'_1}(\R)$ can be obtained easily   using the same argument as in Case 1. Remark that
 $$\mu^{(n'_1)}_{j_0,j_0}(\d x)=\sum_{l=1}^N\mu^{(n_0n_1)}_{j_0,l}(\d x)\ast\mu^{(n_0)}_{l,j_0}(\d x).$$
 Since $M^{\bullet n_0n_1}(\R)>0$ and $\left(M^{\bullet n_0}(\d x)\right)_{i_0,j_0}$ has an absolutely continuous component term, so has the measure
 $\mu^{(n'_1)}_{j_0,j_0}(\d x)$. We
 are therefore in the first case and conclude easily.
 \end{description}
\end{proof}
In particular, we get the following lemma:

\begin{lemma}\label{lem2.7}
 If the hypotheses of lemma \ref{lem2.6} are valid, there exists
 $k_1\geq1$  such that $M^{\bullet k_1}(\R)>0$ and all the terms of
 $M^{\bullet k_1}$ have absolutely continuous components.
\end{lemma}
\begin{proof}
 Take $k_1=4k_0$ with $k_0=(n_0+1)n_1n_0$. The positivity of $M^{\bullet k_1}(\R)$ is an immediate
 consequence of lemma \ref{lem2.6}. In addition,
 $$M^{\bullet 2k_0}=M^{\bullet k_0}M^{\bullet k_0}.$$
By Lemma \ref{lem2.6}, one has $M^{\bullet k_0}(\R)>0$ and $M^{\bullet k_0}$ has an
 absolutely continuous component term $\mu^{(k_0)}_{i'_0,j'_0}$. So according to the above equality, we see that every term of
 $M^{\bullet 2k_0}(e_{i'_0})$ and $M^{\bullet 2k_0}(e_{j'_0})$ has an absolutely continuous component. It is thus clear that all the terms of
 the matrix $M^{\bullet k_1}=M^{\bullet 4k_0}$ have an absolutely continuous component.
\end{proof}

\subsection{Proof of Theorem \ref{thm55}}\label{appendixB}

\begin{proof}[Proof of Theorem \ref{thm55}]

\begin{enumerate}
 \item The first assertion is a direct consequence of the perturbation theorem (see Theorem 9 of Chapter 7 in \cite{Schwartz} for instance).
 \item To prove the second assertion, we will use the following lemmas:
   \begin{lemma}\label{rayonainfini}
     There exist  $\a_1>0$, $\b_1>0$ and $\chi_1\in ]0, 1[$  such that $r(P(\l))\leq\chi_1$ for any $\l\in \mathbb C$ satisfying  $|\Re\ \l|\leq\a_1$ and $|\Im\ \l|\geq\b_1$.
     \end{lemma}
     \begin{lemma}\label{royoncarre}
      For any   $0< a<b$, there exist  $\a_{a, b}>0$  and $\chi_{a, b}\in ]0, 1[$ such  that  $r(P(\l))\leq\chi_{a, b}$  for any $\l\in \mathbb C$ satisfying $|\Re\ \l|\leq\a_{a, b}$ and $a\leq|\Im\ \l|\leq b$.
     \end{lemma} Theorem   \ref{thm55} can thus be proved easily. Indeed, it is sufficient to fix $a,b$ in Lemma \ref{royoncarre} in the following way : $a=\a_0$, $b=\b_1$ given by Lemma  \ref{rayonainfini} and   $\a'_0=\inf(\a_0,\a_1,\a_{a, b}), \chi=\inf(\chi_1,\chi_{a, b})$.  \end{enumerate}
\end{proof}

It remains to prove Lemma \ref{rayonainfini} and Lemma \ref{royoncarre}. We first need the following fact :

\begin{fact}\label{LemApp3}
      Fix $\ga>0$ and let  $f:\R\mapsto\R$  be such that  the function $y\mapsto\e^{\ga \vert y\vert}f(y)$ belongs to $\mathbb{L}^1(\R,\d x)$. Then
$$\lim_{|t|\rightarrow+\In\atop t\in\R}\sup_{|a|\leq\ga}\Big|\int_\R\e^{(a+it)y}f(y)\d y\Big|=0.$$
    \end{fact}

 \begin{proof}
   For any $\varepsilon>0$, there exists a function $f_\varepsilon\in C^1$ and with compact support $\subset [-M, M]$ and such that
\begin{equation}\label{eqApp3}
\int_\R\e^{\ga|y|}|f(y)-f_\varepsilon(y)|\d y<\varepsilon.
\end{equation}
For any $a\in[-\ga,\ga]$, one has
 \begin{equation*}
   \begin{split}
     \Big|\int\e^{(a+it)y}f(y)\d y\Big|&\leq\Big|\int_\R\e^{(a+it)y}\Big(f(y)-f_\varepsilon(y)\Big)\d y\Big|+\Big|\int_\R\e^{(a+it)y}f_\varepsilon(y)\d y\Big|\\
       &\leq\int_\R\e^{\ga|y|}|f(y)-f_\varepsilon(y)|\d y+\frac{\e^{\ga M}}{|t|}\int_\R|f'_\varepsilon(y)|\d y
   \end{split}
 \end{equation*}
Using (\ref{eqApp3}) and letting $t\to +\infty$, one can obtain the expected result.
  \end{proof}

\begin{proof}[Proof of Lemma \ref{rayonainfini}]
 Set $M=\Big(p_{i,j}F(i,j,\d x)\Big)_{i,j}$. By Lemma \ref{lem2.7},  there exists $k_1\geq0$ such that all the terms of the matrice $M^{\bullet k_1}$ have absolutely continuous components. Using the fact that
$$M^{\bullet k_1}_{i,j}(\d x)=\varphi_{k_1,i,j}(x)\d x+\theta_{k_1,i,j}(\d x),$$
where for any $(i,j)\in E\times E$,
\begin{description}
 \item [$\bullet$] the function $\varphi_{k,i,j}$ is strictly
positive, belongs to $\mathbb{L}^1(\R,\d x)$ and satisfies
$$0<\int\varphi_{k,i,j}(x)\d x\leq1;$$
 \item[$\bullet$] $\theta_{k,i,j}(\d
x)$ is a singulary measure with respect to the Lebesgue measure such
that
$$0\leq\int\theta_{k,i,j}(\d x)<1.$$
\end{description}
Recall that the matrice containing the singulary measures $\theta_{k_1,i,j}$ is denoted by $\Theta_{k_1}(\d x)$ and its relative Laplace transform term by term is denoted by $\L(\Theta_k)(\l)$, for $|\Re\ \l|\leq\a_0$.\\
By Lemma \ref{LemApp3}, we have
\begin{equation*}
  \begin{split}
    \limsup_{|t|\rightarrow+\In\atop t\in\R}\sup_{|a|\leq\a_0}\|P^{k_1}(a+it)\|&\leq\lim_{|t|\rightarrow\In\atop t\in\R}\sup_{|a|\leq\a_0}\|\L(H_{k_1})(a+it)\|+\limsup_{|t|\rightarrow+\In\atop t\in\R}\sup_{|a|\leq\a_0}\|\L(\Theta_{k_1})(a+it)\|\\
&\leq\sup_{|a|\leq\a_0}\|\L(\Theta_{k_1})(a)\|.
  \end{split}
\end{equation*}
Moreover,  $\|\L(\Theta_{k_1})(0)\|=1-\delta$ with $\delta\in ]0, 1[$ ;  by  continuity of the map  $x\in\R\mapsto\L(\Theta_{k_1})(x)$,   there thus exists a real number $\a_1$ such that
$$\sup_{|a|\leq\a_1}\|\L(\Theta_{k_1})(a)\|\leq1-\delta/2<1.$$
Set $\chi_1=1-\delta/4$ and choose $\b_1>0$ such that for any $\l\ \in\C$ satisfying $|\Re\ \l|\leq\a_1$ and $|\Im\ \l|\geq\b_1$ one gets $\|P^{k_1}(\l)\|\leq\chi_1, $
which implies $r(P(\l))\leq\chi_1^{1/k_1}<1.$
\end{proof}

\begin{proof}[Proof of Lemma \ref{royoncarre}]
Fix  $\l\in\C$ s.t. $|\Re\ \l|\leq\a_1$ and $|\Im\ \l|\in [a,b]$. Since for any $i,j\in E$ the measure   $P^{\bullet k_1}_{i,j}$  has an absolutely continuous component,  one gets
$$|P^{k_1}_{i,j}(\l)| < P^{k_1}_{i,j}(\Re\ \l),$$
i.e.
$|P^{k_1}_{i,j}(\l)| \leq \rho_\l  P^{k_1}_{i,j}(\Re\ \l)$ with $0<\rho_\l<1$ ; by continuity of the map $\l\mapsto |P^{k_1}_{i,j}(\l)|$,  one gets $\displaystyle \rho_{a, b}:= \sup_{ \stackrel{|\Re\ \l|\leq\a_1}{|\Im\ \l|\in [a,b]}}\rho_\l \in ]0, 1[$. There thus exists $0<\rho<1$ such that
$$|P^{k_1}_{i,j}(\l)|\leq\rho P^{k_1}_{i,j}(\Re\ \l).$$
Therefore,  for any $\l$ such that ${|\Re\ \l|\leq\a_1}$ and $|\Im\ \l|\in [a,b]$, one gets
\begin{equation}\label{eqApp1}
r(P(\l))\leq\rho_{a, b}^{1/k_1}k(\Re\ \l).
\end{equation}
But,   one gets
$$|k(\Re\ \l)-1|\leq|\Re\ \l|\sup_{-\a_1\leq u\leq\a_1}|k'(u)|\leq\a_1 M_{\a_1},$$
whith $ \displaystyle M_{\a_1}=\sup_{-\a_1\leq u\leq\a_1}|k'(u)| <+\infty$. Finally, for $\a_1$  small enough, one gets
$$
\chi_{a, b} := \sup_{ \stackrel{|\Re\ \l|\leq\a_1}{|\Im\ \l|\in [a,b]}}r(P(\l))\   \in \ ]0, 1[.
$$
 \end{proof}

\nocite{boro} \nocite{bra1} \nocite{bra2}


\begin{thebibliography}{99}
\bibitem{boro} Borovkov, K. A. (1980). Stability theorems and estimates of the rate of convergence of the components of factorizations for walks defined on Markov chains.
Theory of Probability ans its applications 25, no. 2, 325-334.
\bibitem{bra1} Bratiichuk, M. S. (1997). Properties of operators of additive Markov processes \Rmnum{1}. Theor. Probab. and Math. Statist., no. 55, 19-28.
\bibitem{bra2} Bratiichuk, M. S. (1998). Properties of operators of additive Markov processes \Rmnum{2}. Theor. Probab. and Math. Statist., no. 57, 1-10.
\bibitem{Schwartz} Dunford, N., and Schwartz, J. (1988). Linear Operators, Part II Spectral Theory, Self Adjoint Operators in Hilbert Space, Johy Wiley \& Sons, INC., USA.
\bibitem{Fell} Feller, W. (1968). An introduction to probability theory and its applications, Johy Wiley \& Sons, INC., USA.
\bibitem{flajolet} Flajolet, P., and Odlyzko, A. (1990). Singularity analysis of generating functions. SIAM J. Discrete Math. 3, 216-240.
\bibitem{geig1} Geiger, J. and Kersting, G. (2000). The survival probability of a critical branching process in random environment. Ther. Verojatnost. i Primenen 45, 607-615.
\bibitem{guivarc'h-hardy} Guivarc'h, Y. and Hardy, J. (1988). Th\'eor\`emes limites pour une classe de cha\^{i}nes de Markov et applications aux diff\'eomorphismes d'Anosov. Ann. Inst. Henri Poincar\'{e} 24, no. 1, 73-98.
\bibitem{quiv} Guivarc'h, Y., Le Page, E., and Liu, Q. (2003). Normalisation d'un processus de branchement critique dans un environnement al\'{e}atoire. C. R. Acad. Sci. Paris. 337, no.9, 603-608.
\bibitem{koz} Kozlov, M. V. (1976). On the asymptotic behavior of the probability of non-extinction for critical branching processes in a random environment. Theory Probab. 21, no. 4, 791-804.
\bibitem{emile} Le Page, E., and Peign\'{e}, M. (1997). A local limit theorem on the semi-direct product of $\mathbb{R}^{*+}$ and $\mathbb{R}^d$. Ann. Inst. Henri Poincar\'{e} 33, 223-252.
\bibitem{ye} Le Page, E., and Ye, Y. (2010). The survival probability of a critical branching process in Markovian random environment. C. R. Acad. Sci. Paris. 348, no. 5, 301-304.
\bibitem{pres2} Presman, E. L. (1969). Factorization methods and a boundary value problem for sum of random variables defined on a Markov chain. Math. USSR-Izv. 3, no. 4, 815-852.
\bibitem{Widd} Widder, D. V. (1946). The Laplace transform, Princeton University Press.

\end{thebibliography}
\end{document}